\documentclass{article}
\usepackage{amsmath}
\usepackage{amsfonts}
\usepackage{amscd}
\usepackage{amsthm}
\usepackage{xcolor}
\usepackage{graphicx}

\usepackage{geometry}
\allowdisplaybreaks

\newtheorem{theorem}{Theorem}[section]
\newtheorem{lemma}[theorem]{Lemma}
\newtheorem{cor}[theorem]{Corollary}

\newtheorem{remark}[theorem]{Remark}
\theoremstyle{definition}
\newtheorem{definition}[theorem]{Definition}

\newtheorem{innercustomthm}{Theorem}
\newenvironment{customthm}[1]
  {\renewcommand\theinnercustomthm{#1}\innercustomthm}
  {\endinnercustomthm}

\newcommand{\abs}[1]{\left\lvert#1\right\rvert}
\newcommand{\norm}[1]{\left\lVert#1\right\rVert}
\newcommand{\set}[1]{\left\{#1\right\}}
\newcommand{\interval}[1]{\left[#1\right]}
\newcommand{\expr}[1]{\left(#1\right)}
\newcommand{\id}{\operatorname{id}}
\newcommand\numberthis{\addtocounter{equation}{1}\tag{\theequation}}
\newcommand\IS{\mathit{IS}}
\newcommand\FS{\mathit{FS}}
\newcommand\HS{\mathit{HS}}
\newcommand\VS{\mathit{VS}}

\newcommand{\dom }{\,{\rm dom}\,}

\newcommand{\sgn }{{\rm sgn}\,}

\newcommand{\im }{\mbox{im} \,}

\newcommand{\Id }{\mbox{Id} \,}

\newcommand{\comment}[1]{\mbox{}}

\def\qed{{\hfill{\vrule height5pt width3pt depth0pt}\medskip}}

\author{
  Jacek Kubica\\
  \texttt{jacek.kubica@ii.uj.edu.pl}
  \and
  Piotr Zgliczyński\\
  \texttt{umzglicz@cyf-kr.edu.pl}
  \and
  Piotr Kalita\\
  \texttt{piotr.kalita@ii.uj.edu.pl}
}
\title{ Pitchfork bifurcation and heteroclinic connections in the Kuramoto--Sivashinsky PDE}

\begin{document}

\maketitle

\begin{abstract}
  \textcolor{black}
  {
 We present a method for the complete analysis of the dynamics of dissipative Partial Differential Equations (PDEs) undergoing a pitchfork bifurcation. We apply our technique to the Kuramoto--Sivashinsky PDE on the line to obtain a computer-assisted proof of the creation of two symmetric branches of non-symmetric fixed points and heteroclinic connections between the symmetric fixed point and the new ones. The range of parameters is given explicitly and is large enough to allow for the rigorous continuation of the fixed points and heteroclinic connections created during the bifurcation.
  }
\end{abstract}


{\let\thefootnote\relax\footnote{{The work of all three authors was supported by National Science Center (NCN) of Poland under project No. UMO-2016/22/A/ST1/00077. The work of Piotr Kalita was moreover supported by FAPESP, Brazil grant 2020/14075-6 and Ministerio de Ciencia e Innovaci\'{o}n of Kingdom of Spain under	project No. PID2021-122991NB-C21.}}}

\section*{Introduction}

\textcolor{black}
{
 We present a method for the complete analysis of the dynamics of dissipative PDEs undergoing a pitchfork bifurcation. By this, we mean the following. Assume that we have a parameter-dependent system with a 'reflection' symmetry containing a symmetric fixed point, where, as the parameter changes, a single eigenvalue passes (increases) through zero. At the bifurcation point, two symmetric branches of non-symmetric fixed points are created, along with heteroclinic connections between the symmetric fixed point and the new ones. We also obtain an explicit parameter range over which the new fixed points and heteroclinic connections exist. This range is large enough to allow for rigorous continuation of the fixed points and heteroclinic connections created during the bifurcation, away from the bifurcation locus.
}

 The dissipative PDE to which we apply
our \textcolor{black}{technique} is the Kuramoto--Sivashinsky equation \cite{KT,S}
(in the sequel we will refer to it as the KS equation)
\begin{equation}\label{eq_ks_pde}
u_t = -\mu u_{xxxx} - u_{xx} + (u^2)_x, \qquad \mu>0
\end{equation}
where $x \in \mathbb{R}$, $u(t,x) \in \mathbb{R}$. The equation is equipped with the
odd and periodic boundary conditions
\begin{equation} \label{eq_ks_bc}
  u(t,x)=-u(t,-x), \qquad u(t,x) = u(t,x+ 2 \pi).
\end{equation}

\par
Expanding a solution of the KS equation in the sine Fourier basis
$$u(t, x) = \sum_{k = 1}^\infty -2a_k \sin(kx),$$
the KS equation
becomes an infinite-dimensional ladder of ODEs (see \cite{ZM})
\begin{align} \label{eqn_ks_fourier_with_mu_intro}
  \begin{split}
    \dot{a_k}&=k^2(1-\mu k^2)a_k - k \sum_{i=1}^{k-1} a_i a_{k-i} + 2 k \sum_{i=1}^{\infty} a_i a_{k+i}, \quad k=1,2,\dots
  \end{split}
  \end{align}
We denote the right-hand side by $F$.
We see that  \textcolor{black}{as $\mu$ decreases through $1$,} the origin loses stability and because of the symmetry $F\expr{R\expr{a}} = R\expr{F\expr{a}}$, for $R$ given by
\begin{align*}
  a_{2k} &\mapsto a_{2k}, \\
  a_{2k + 1} &\mapsto -a_{2k + 1},
\end{align*}
we expect that the corresponding bifurcation is a pitchfork bifurcation,
also known as the symmetry-breaking bifurcation (see \cite{ZBif}[Section 6.1] for more detailed discussion).

Our first main result is the analytical proof of  \textcolor{black}{a}  pitchfork bifurcation from zero in the KS equation
when $\mu = 1$,
as stated in the Theorem \ref{thm_bif}. This proof is done in Section \ref{sec_proof_of_bifurcation}.
Essentially (omitting some subtleties) what we prove is that there exists a parameter range $\mu_+ < 1 < \mu_-$ and
the "big" self-consistent bounds $V$ (here by "big" we mean the same for every $\mu$ in the given range)
such that
\begin{enumerate}
  \item[(i)] For $\mu \in \interval{\mu_+, \mu_-} \setminus \set{1}$ the origin is a hyperbolic fixed point.
  \item[(ii)] For $\mu \in [1, \mu_-]$ the origin is an attracting fixed point. Moreover, $\set{0}$ is the
        maximal invariant set in $V$.
  \item[(iii)] For $\mu < 1$ there exist two non-zero hyperbolic fixed points $u_{\pm}^\mu \in V$. Moreover, there exist heteroclinic
        connections from the origin to those points and the maximal invariant set in $V$ is the set
        consisting of those fixed points, the origin and the mentioned heteroclinic connections.
\end{enumerate}

Let us stress that because we have the complete characterization of the maximal invariant sets,
hence we obtain the complete description of the dynamics near the bifurcation. This should be
contrasted with the other work on the bifurcations \cite{ZBif, AriKo, BQ} where the dynamics
is not considered.

\par
Let us  discuss briefly problems and techniques related to establishing claims (i),(ii),(iii).

\par
To handle the issue of infinite dimension in system (\ref{eqn_ks_fourier_with_mu_intro}) and  its relations with finite dimensional Galerkin projections  we use the approach based on the self-consistent bounds developed in \cite{ZM, ZNS, ZA}.
This allows us to apply finite-dimensional geometric tools to dissipative systems.  The method is described in Sections \ref{sec:method}.

A usual in the analysis of bifurcations the system has to be brought first to a suitable normal form (for the exposition of the normal form theory, see for example \cite[Chapter 5]{Arn}).
Transformation of the KS equation
to the normal form
is discussed in detail in Section \ref{sec_proof_pitchfork_ks} and
Appendix \ref{sec_app_normal_form}.

\par
\textcolor{black}{The} notion of hyperbolicity that we use differs from the functional-analytic one based on the spectrum. Connecting
such definitions with what they mean for the dynamics is usually complicated in infinite dimensional dynamical
system (see \cite{FSV, SS} for the approach connecting spectrum with the dynamics). Instead, we focus on the
expected behavior of the solutions of the system near the hyperbolic fixed points. We formulate it precisely
in Section \ref{sec_manifolds}, in Definition \ref{def_hyperbolicity}. We also prove that our definition
implies the standard one in case of the finite-dimensional system. 

To describe the local behavior near fixed points we use the logarithmic norms and the cone conditions.
The logarithmic norms and their basic properties are recalled
in Appendix \ref{app_lognorms} and their adaptation in context of the self-consistent bounds is described in Section \ref{sec:cov-proj}
based on \cite{ZNS}.
The cone conditions and their relation with hyperbolicity are described also in Section \ref{sec_manifolds} and
their verification in the context of the self-consistent bounds in Section \ref{sec_cc_ver}.
\par
Claim (iii) describes the dynamics after the bifurcation. Its proof consists of three parts.
First, we establish that the unstable manifold of the source is a graph over the bifurcation direction
(which is an unstable direction near the origin). To prove that, we verify the cone conditions.
The relation of the cone conditions to the invariant manifolds is also explained in Section \ref{sec_manifolds}
and their adaptation to the context of the self-consistent bounds is described in Section \ref{sec_cc_ver}.
Next, by a straightforward computation we check that some point on the unstable manifold is transported to the
attracting region of a fixed point born after the bifurcations.

\par
In Section \ref{sec_pitchfork_models} we discuss two simple ODE models of  \textcolor{black}{a}  pitchfork bifurcation
to demonstrate our techniques and discuss the difficulties we need to overcome, first without
any unstable directions and then with one unstable direction.
\par
\textcolor{black}{The} proof of the Theorem \ref{thm_bif} is purely analytical and does not need any computer assistance. Its shortcoming
is that it does not give
us the explicit range $\interval{\mu_+, \mu_-}$. This is why one of the goals in Section \ref{sec_gen_bif_thm} is
to extract from the proof of Theorem \ref{thm_bif} the inequalities used in it in a form which can be
verified rigorously on the computer using the interval arithmetic. We also go a step further and give additional inequalities
which allow us to prove  \textcolor{black}{a} pitchfork bifurcation when some directions are unstable. This constitutes the second main result
of this paper.
Below by  \textcolor{black}{a}  pitchfork bifurcation on a given parameter range we mean as before the bifurcation with
the full description of the dynamics.

\begin{customthm}{\ref{thm_comp_assisted_bif1}}
   Pitchfork bifurcation occurs in the KS equation for the
  parameter $\mu$ in the range $\interval{0.99, 1.01}$.
\end{customthm}
\begin{customthm}{\ref{thm_comp_assisted_bif2}}
   Pitchfork bifurcation occurs in the KS equation
  for the parameter $\mu$ in the range $\interval{0.25 - 0.0002, 0.26}$.
\end{customthm}

Let us mention that Theorem \ref{thm_comp_assisted_bif2} does not seem to give
any additional information compared to Theorem \ref{thm_comp_assisted_bif1}, because
if $u$ is the solution of the KS equation for some $\mu > 0$, then
$\tilde{u}(t, x) := ku(k^2t, kx)$ is the solution of the KS equation for $\frac{\mu}{k^2}$
(See \cite[Lemma 6.1]{ZBif}). It would still be enough for us to test our method,
nevertheless this transformation would give us
a complete description of the dynamics only on the invariant subspace of even modes, while
approach in this paper gives us complete description near the bifurcation on the whole space.
\par
\textcolor{black}{Our} final main results are the computer assisted proofs of the heteroclinic connection at
the end of the parameter range from Theorem \ref{thm_comp_assisted_bif1} and another
away from the bifurcation.
\begin{customthm}{\ref{thm_comp_assisted_connections}}
For system (\ref{eq_ks_pde}) with $\mu \in \set{0.99, 0.75}$ there exists  a heteroclinic
connection between two fixed points: the unstable zero solution
and the attracting fixed point.
\end{customthm}
This result can be seen as the continuation of Theorem \ref{thm_comp_assisted_bif1}.
We chose one parameter at the end of the parameter range in the mentioned theorem ($\mu = 0.99$,
correspondingly $\lambda_1(\mu) = 0.01$)
and the other quite away from it ($\mu = 0.75$,
correspondingly $\lambda_1(\mu) = 0.25$). This gives us assurance that given enough
computation time it should pose no problem to continue the heteroclinic connection for
every intermediate value of $\mu$.
\textcolor{black}{Let us mention that rigorous construction of heteroclinic connections between equilibria for dynamical systems driven by dissipative PDEs with the use of computer assisted methods was done, for example, in \cite{CyrWan} for the Ohta--Kawasaki model or, more recently, for complex-valued heat equation with quadratic nonlinearity in \cite{JLT}. However, none of those heteroclinic connections are associated with a bifurcation.}
\par
The method of the proof of Theorem \ref{thm_comp_assisted_connections} is similar to the proof of the claim (iii) mentioned above.
In Appendix \ref{sec_cc_and_ln_ver} we provide means of obtaining rigorous estimate of the unstable manifold of the source point
and of the basin of attraction of the target point. Then we "connect" those estimates using the rigorous integration
based on the algorithm for the dissipative PDEs described in \cite{ZPer2}. Data regarding
the obtained fixed points and the integration is in Section \ref{sec_data_heteroclinic}.
\par
\par
In computations the CAPD \cite{Capd} library is used. Code is available at \cite{github}. 
\textcolor{black}{None of the computations time exceeded 20 minutes on regular PC. The longest computation
was necessary to prove the heteroclinic connection right after the bifurcation.}

% \section{Topological setup}

\section{Cone conditions, unstable manifold and hyperbolicity for the local semiflows in  Hilbert space\textcolor{black}{s}}
\label{sec_manifolds}

Let $H$ be a Hilbert space,
$X_n$ be an $n$-dimensional subspace of $H$ and set $Y_n := X_n^\bot$. We denote by $P_n:H\to X_n$ and $P_n^\perp :H \to Y_n$ respectively the orthoprojections onto $X_n$ and $Y_n$.

\begin{definition}
    We say that a partial (with respect to
    the first variable) map $\varphi$ from $[0, \infty) \times H$ to $H$
    is a \emph{local semiflow} if
    \begin{itemize}
      \item $\varphi$ is continuous,
      \item for every $z \in H$ there exists $t_{\max}(z) \in \interval{0, \infty}$
            such that $\varphi(t, z)$ exists only for $t \in [0, t_{\max}(z))$ (or for $t \in [0, \infty)$ if $t_{\max} = \infty$).
    \item for every $z \in H$ we have $\varphi(0, z) = z$,
      \item for every $z \in H$ and for $s, t > 0$ such that $s + t \leq t_{\max}(z)$ we have $\varphi(t, \varphi(s, z)) = \varphi(t + s, z)$.
    \end{itemize}
    If $\varphi$ is a full map (i.e. $t_{\max}(z) = \infty$ for all $z \in H$), then we simply say that it is a \emph{semiflow}.
\end{definition}

\textcolor{black}
{
The notion of isolating cuboid introduced below is a   special case of the
isolating block  known from the Conley index theory, see \cite{MM}.
}

\begin{definition}
Let $c_u:X_n \to \mathbb{R}^n$ be a homeomorphism.
Consider a set $N \subset H$ given by
$$
N := c_u^{-1} (\overline{B}(0,1)) \oplus T
$$
and for $\varepsilon > 0$ denote
$$
N^\varepsilon := c_u^{-1}(\overline{B}(0,1+\varepsilon)) \oplus T,
$$
where $B(0,1) \subset \mathbb{R}^n$ is an open ball and
where $T \subset Y_n$ is a compact set.
\par
Assume that  we have a local semiflow $\varphi$  \textcolor{black}{and let $\varepsilon > 0$}.
We say that $N$ is an $\varepsilon$-isolating cuboid for  $\varphi$  if there exists a time $t(\varepsilon) > 0$
such that for every $s \in (0, t(\varepsilon)]$ we have
\begin{itemize}
    \item[(I0)] if for $z \in N^\varepsilon$ we have $t_{\max}\expr{z} < \infty$,
                then there exists $t > 0$ such that $\varphi(t, z) \not \in N^\varepsilon$
                (in other words, all points that do not leave $N^\varepsilon$ \textcolor{black}{forward in time} 
                have \textcolor{black}{a}
                 full forward trajectory),
    \item[(I1)] $\varphi\expr{s, N} \subset N^\varepsilon$,
    \item[(I2)] \textcolor{black}{
       $\varphi\expr{s, N^\varepsilon \setminus N} \cap N = \emptyset$.}
\end{itemize}
\textcolor{black}{
In the sequel when we want to clearly point to which coordinate directions span $X_n$,  we say that $N$ is
\emph{$\varepsilon$-isolating cuboid with the unstable directions $x_1, \dots, x_n$}.
}
\end{definition}
\begin{remark}
    In our approach to the Kuramoto--Sivashinsky equation we use $T$ of the form
    $$
    T = \prod_{k = n + 1}^\infty I_k,
    $$
    where $I_k := \interval{a_k^-, a_k^+}$. In this case we call sets
    $$N^+ := c_u^{-1}(B(0,1)) \times \set{z \in T \mid \exists k \geq n + 1 \; z_k \in \partial I_k}$$
    the
    \emph{entry set} (since by (I1) every point from this set flows into $N$)
    and $N^- := c_u^{-1}(\partial B(0, 1)) \times T$
    the \emph{exit set} (since by (I2) every point from it flows out of $N$).
\end{remark}

\begin{definition}
    Let $\varphi$ be a local semiflow on $N \subset H$ and let $\bar{z} \in N$.
    For an interval $I \subset \mathbb{R}$ such that $0 \in I$ we say that a function
    $z: I \to N$ is a \emph{trajectory} through $\bar{z}$ in $N$ if
    \begin{itemize}
        \item $z(0) = \bar{z}$,
        \item for all $t < 0$ such that $t \in I$ we have $\varphi(-t, z(t)) = \bar{z}$,
        \item for all $t > 0$ such that $t \in I$ we have $\varphi(t, \bar{z}) = z(t)$.
    \end{itemize}
    We say that a trajectory $z(t)$ is a \emph{full backward} (\emph{full forward}) trajectory through $\bar{z}$ in $N$
    if $(-\infty, 0] \subset I$ ($[0, \infty) \subset I$). We say that it is a
    \emph{full trajectory} through $\bar{z}$ trajectory if it is a full backward trajectory and a full forward trajectory through $\bar{z}$.
\end{definition}

\begin{remark}
    As $\varphi$ is only a semiflow, even if a backward trajectory
    exists on some time interval, it might not be unique.
\end{remark}

\begin{definition}
    Let $\varphi$ be a local semiflow on $N \subset H$ and assume that $z_0 \in N$
    is its fixed point.
    We define \emph{unstable manifold} $W^u_N\expr{z_0}$ of $z_0 \in N$
    as a set of such points $\overline{z} \in N$ that
    \begin{itemize}
        \item $z$ has  \textcolor{black}{a} full backward trajectory in $N$,
        \item for every full backward trajectory $z(t)$ of $\overline{z}$ in $N$ we have $\lim_{t \to -\infty} z(t) = z_0$.
    \end{itemize}
\end{definition}

\begin{definition}
    Let $\varphi$ be a local semiflow on $N \subset H$ and assume that $z_0 \in N$
    is its fixed point.
    We define \emph{stable manifold} $W^s_N\expr{z_0}$ of $z_0 \in N$
    as a set of such points $\overline{z} \in N$ for which the forward trajectory is full
    in $N$ and we have $\lim_{t \to \infty} z(t) = z_0$.
\end{definition}

\par
Our goal is to give conditions under which the (un)stable manifold in $N$
is a graph over the (un)stable directions. To this end we define the cone conditions.
Let us discuss some notation first.
Let $\mathcal{Q}: H = X \oplus Y \ni (x, y) \mapsto \alpha(x) - \beta(y) \in \mathbb{R}$,
where $\alpha, \beta$ are \textcolor{black}{continuous} positive defined quadratic forms.
The positive and negative cones on the set $N^\varepsilon$
that we work with here are given by
$Q^+ := \set{z \in N^{\varepsilon} \mid \mathcal{Q}(z) > 0}, Q^- := N^\varepsilon \setminus Q^+$.
In our later uses we will simply use  $\mathcal{Q}(x, y) = \norm{x}^2 - \norm{y}^2$,
but our results do not depend on it.
\par
Throughout this section we keep $\mathcal{Q}$ fixed. By the
cone conditions we mean the forward-invariance of the positive
cones $z + Q^+$ for every $z \in N^\varepsilon$ together with
expansion (contraction) in the positive (negative, i.e. in $z + Q^-$) cones
expressed in terms of $\mathcal{Q}$. More precisely,
we state the following definition.

% In the proof of the stable and unstable manifold theorems, we need a stronger notion of the cone conditions.
% Additional property we demand from $\mathcal{Q}$ in the definition below imply that
% $\mathcal{Q}$ is a two-point Liapunov function.

\begin{definition}
Let $N$ be an $\varepsilon$-isolating cuboid for a continuous local
semiflow $\varphi$.
% $$
% \mathcal{Q}:H \to \mathbb{R}
% $$
% be a map defined by $\mathcal{Q}(u) = \alpha(c_u(P_n u)) - \beta(P_n^\bot u)$,
% where $\alpha:\mathbb{R}^n\to \mathbb{R}$ and $\beta:Y_n \to \mathbb{R}$ are continuous positively defined quadratic forms.
We say that the semiflow $\varphi$ satisfies the \emph{cone conditions}
% (the \emph{weak cone conditions})
on $N$ if for every $z_1, z_2 \in N^\varepsilon$ such that $z_1\neq z_2$
% ($\mathcal{Q}\expr{z_1 - z_2} > 0$)
the derivative $\frac{d}{dt} \mathcal{Q}\expr{\varphi\expr{t, z_1} - \varphi\expr{t, z_2}}|_{t=0}$ exists and there holds
$$
\frac{d}{dt} \mathcal{Q}\expr{\varphi\expr{t, z_1} - \varphi\expr{t, z_2}}|_{t=0} > 0.
$$
\end{definition}

Treating $\mathcal{Q}$ as a 'metric' the above definition immediately implies the expansion in the positive cone and contraction in the negative one. Later, when discussing the hyperbolicity, we introduce the strong
cone conditions which give us expansion and contraction with respect to the norm,
but we do not need them to prove the stable and unstable manifold theorems.
\begin{figure}
\begin{center}
\includegraphics[scale=0.4]{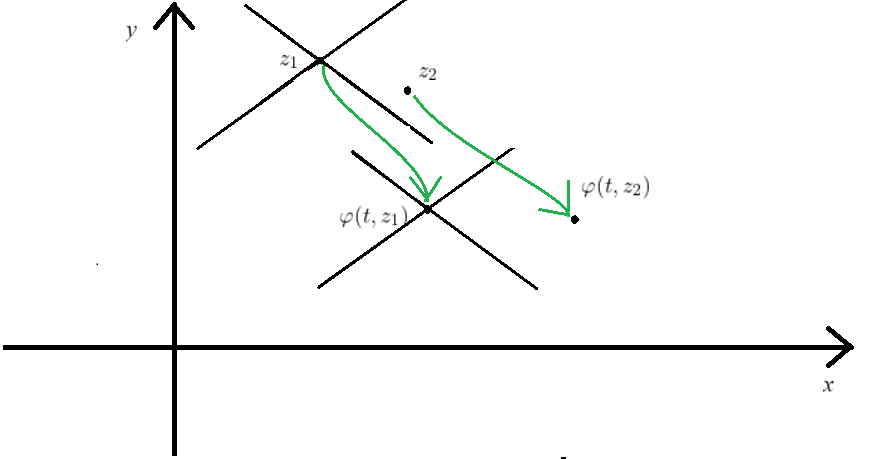}
\caption{\textcolor{black}{Illustration of the property of the positive cones invariance given in Definition \ref{def_conecond} and of the cone conditions. The  positive cones are given by $\mathcal{Q} > 0$. The point $z_2$ belongs to the positive cone of $z_1$ and thus $\varphi(t,z_2)$ must belong to the positive cone of $\varphi(t,z_1)$.                     Note that                        although the picture shows norm expansion and contraction respectively
                           in $x$ and $y$ direction, to demonstrate these properties we use the strong cone conditions which are
                           introduced later in this section rather than just the cone conditions.}}
\end{center}
\end{figure}
\par
Now we state formally our remark that the cone conditions give us the invariance of the positive cones.

\begin{definition} \label{def_conecond}
    Let $N$ be an $\varepsilon$-isolating cuboid for a continuous local
    semiflow $\varphi$. We say that $\varphi$ has
    \emph{invariant positive cones} on $N$ if for $z_1, z_2 \in N$ such that $z_1 - z_2 \in Q^+$
        we have for every $s \in \interval{0, t(\varepsilon)}$ that $\varphi(s, z_1) - \varphi(s, z_2) \in Q^+$.
\end{definition}

\begin{lemma}
    Let $N$ be an $\varepsilon$-isolating cuboid for a continuous local
semiflow $\varphi$. If the cone conditions are satisfied on this set, then the flow
has invariant positive cones.
\end{lemma}
\begin{proof}
    Assume that the cone conditions are satisfied.
    Let $z_1, z_2 \in N$ be such that $z_1 - z_2 \in Q^+$.
    Let $s < t(\varepsilon)$. By (I1) we have $\varphi(s, z_1) - \varphi(s, z_2) \in N^\varepsilon$, thus
    by the cone conditions
\begin{align*}
    \mathcal{Q}\expr{\varphi(s, z_1) - \varphi(s, z_2)} > \mathcal{Q}\expr{z_1 - z_2} > 0.
\end{align*}
\end{proof}

Now we introduce the notion of the disks, which simply are the Lipschitz functions
over the unstable or stable directions expressed in the terms of $\mathcal{Q}$.
\begin{definition}
     Let $N \textcolor{black}{= P_n N \oplus P_n^\perp N} \subset H$.
    We say that a continuous function $h: P_n N \to N$ is a \emph{horizontal disk in $N$} if
    for every $x \in P_n N$ we have $P_n h(x) = x$ and
    for every $x_1, x_2 \in P_n N$ such that $x_1 \neq x_2$ we have $\mathcal{Q}\expr{h\expr{x_1} - h\expr{x_2}} > 0$.
\end{definition}

\begin{definition}
     Let $N \textcolor{black}{= P_n N \oplus P_n^\perp N} \subset H$.
    We say that a continuous function $v: P_n^\bot N \to N$ is a \emph{vertical disk in $N$} if
    for every $y \in P_n N$ we have $P_n^\bot v(y) = y$ and
    for every $y_1, y_2 \in P_n^\bot N$ such that $y_1 \neq y_2$ we have $\mathcal{Q}\expr{v\expr{y_1} - v\expr{y_2}} < 0$.
\end{definition}

% \begin{definition}
%     Let $\varphi$ be a local semiflow on $N$. We say that
%     $S \subset N$ is \emph{forward-invariant} (\emph{backward-invariant})
%     \emph{relative to $N$} if for every $z \in S$ and for every $t > 0$ such that
%     $\varphi(\interval{0, t}, z) \subset N$ ($t < 0$ such that $\varphi(\interval{-t, 0}, z) \subset N$)
%     we have  $\varphi(\interval{0, t}, z) \subset S$ ($\varphi(\interval{-t, 0}, z) \subset S$).
% \end{definition}

In the proofs of the theorems in the proceeding part of this section we assume that $c = \id$,
as the arguments will be easy to generalize. Consequently, in those proofs we will
have $P_n N = \overline{B}(0, 1)$.
\par
We now show that if $N$ is an $\varepsilon$-isolating cuboid for the
local semiflow $\varphi$ for which the positive cones are invariant and $h$
is a horizontal disk in $N$,
then for small times $t > 0$ the set $\varphi(t, h(P_n N))$ contains the
image of another horizontal disk in $N$.

\begin{lemma} \label{disk_image_lemma}
    Let $N$ be an $\varepsilon$-isolating cuboid for a continuous local
    semiflow $\varphi$ and let $h: P_n N \to N$ be a horizontal disk in $N$. Assume that the positive cones are invariant for $\varphi$.
    % \begin{itemize}
    %     \item[(CF)] For every $z_1, z_2 \in N$,
    %     if $z_1 - z_2 \in Q^+ := \set{z \in N \mid \mathcal{Q}\expr{z} > 0}$, then for every $s \in \interval{0, t(\varepsilon)}$ we have
    %     $\varphi(s, z_1) - \varphi(s, z_2) \in Q^+$.
    % \end{itemize}
    Then for every $s \in (0, t(\varepsilon)]$ there exists a horizontal disk $h^*: P_n N \to N$
    such that for every $x \in P_n N$ there exists $\overline{x}$ such that for all $\tau \in \interval{0, s}$
    we have $\varphi\expr{\tau, h\expr{\overline x}} \in N$ and $h^*(x) = \varphi\expr{s, h\expr{\overline{x}}}$.
\end{lemma}
\begin{proof}
    % For every $\overline{x}\in \overline B(0,1)$ and $s\in (0,t(\varepsilon)]$ we decompose
    %  $\varphi\expr{s,  \overline {x} + h\expr{\overline{x}}} = w + z$  with $z \in X_n^\bot$ and $w \in X_n$. By (I1), $z\in T$ and $w \in \overline{B}(0,1 + \varepsilon)$.
    %  \par
    We fix $s \in \interval{0, t(\varepsilon)}$.
    \par
     First, we need to show that for every $x \in \overline B(0, 1)$ there exists a unique
     $\overline{x} \in \overline{B}(0, 1)$ such that for some $y \in T$ we have
     $\varphi\expr{s, h\expr{\overline{x}}} = x + y$.
     \par
     To prove the uniqueness, fix $s \in (0, t(\varepsilon)]$ and
     take $\overline{x}_1, \overline{x}_2 \in \overline{B}(0, 1)$ such that $\overline{x}_1 \neq \overline{x}_2$.
     Since $\mathcal{Q}\expr{h\expr{\overline{x_1}} - h\expr{\overline{x_2}}} > 0$, so by the invariance
     of the positive cones we have \\
     $\mathcal{Q}\expr{\varphi\expr{s, h\expr{\overline{x_1}}} - \varphi\expr{s, h\expr{\overline{x_2}}}} > 0$ and the claim follows.
     \par
     To prove the existence, we use the local Brouwer degree,
     whose basic properties we recall in Appendix \ref{app_brouwer_degree}.
     % By (I2) for $x \in B(0, 1)$ we have $\varphi\expr{s, x + h\expr{x}} \in \overline{B}(0,1+\varepsilon)$, thus we can consider a map $\Phi_s:\overline B(0,1 ) \to \overline{B}(0,1+\varepsilon)$ defined by
     For $\tau \in \interval{0, s}$ we define  map $\Phi_{\tau}:\overline B(0,1 ) \to \mathbb{R}^n$ by
     % $$
     % \overline{B}(0, 1)\ni x \mapsto  P_n \varphi\expr{s, x + h\expr{x}} \in \overline{B}(0,1+\varepsilon).
     % $$
     $$
     \overline{B}(0, 1)\ni \overline x \mapsto  P_n \varphi\expr{\tau, h\expr{\overline x}} \in \mathbb{R}^n.
     $$
     Fix $x \in B(0, 1).$ We will show that $\deg(\Phi_s, B(0, 1), x) \neq 0$,
     which implies that there exists $\overline{x} \in B(0, 1)$ such that $\Phi_s(\overline{x}) = x$.
     Define the homotopy $H$ by
     $$
     H_r(\overline x) = \Phi_{rs}(\overline x) \quad \textrm{for}\quad r\in [0,1].
     $$
     Let $\tau\in [0, s]$. Since by (I2) for $\overline x \in \partial B(0, 1)$ we have
     $\Phi_{\tau}(\overline x) \not\in B(0, 1)$,
     it follows that for each $r \in \interval{0, 1}$ we have
     \begin{equation} \label{eqn_we_can_use_homotopy}
     \overline{B}(0,1) \cap H_r\expr{\partial B\expr{0, 1}} = \emptyset.
     \end{equation}
     Since $H_0 = \id, H_1 = \Phi_s$ and since (\ref{eqn_we_can_use_homotopy}) holds,
     by the homotopy property of the local Brouwer degree (Theorem \ref{thm_homotopy_property}) we have
     $$ \mathrm{deg}(\Phi_s, B(0,1), x) = \mathrm{deg}(\id, B(0,1), x) = 1.$$ In consequence, we get
     the needed existence for $x \in B(0, 1)$ by the existence property
     of the Brouwer degree (Theorem \ref{thm_solution_property}).
     \par
     For $x \in \partial B(0, 1)$ the existence of $\overline{x} \in \overline{B}(0, 1)$ such
     that $\Phi_s(\overline x) = x$ follows easily by the continuity of $\varphi$ and the compactness of $N^{\varepsilon}$.
     \par
     We have proved that there exists a function $h^*: \overline B(0,1) \to N$ such that
     for all $x \in \overline{B}(0, 1)$ there exists $\overline{x} \in \overline{B}(0, 1)$ such that
     $\varphi(s, h(\overline{x})) = h^*(x)$.
     It remains to prove that
     $h^*$ is a horizontal disk, but it immediately follows by the invariance of the positive cones.
    %  By the cone conditions and by (I1) and (I2), for $\tau \in \interval{0, s}$ and $x_1, x_2 \in \overline{B}(0, 1)$, $x_1 \neq x_2$, we have
    %  $$
    %  \mathcal{Q}\expr{\varphi\expr{s, h\expr{x_1}} - \varphi\expr{s, h\expr{x_2}} > \mathcal{Q}\expr{h\expr{x_1} - h\expr{x_2}}} > 0,
    %  $$
    %  which concludes the proof.
\end{proof}

\begin{lemma} \label{lemma_fixed_point}
    If $N$ is an $\varepsilon$-isolating cuboid for a continuous local
    semiflow $\varphi$ which satisfies the cone conditions on $N$,
    then $\varphi$ has a unique fixed point $z_0 \in N$. Moreover, for any
    horizontal disk $h$ in $N$ there exists $z \in h(P_n N)$ such that
    $\lim_{t \to \infty} \varphi(t, z) \in N$ and for this point
    $\lim_{t \to \infty} \varphi(t, z) = z_0$.
\end{lemma}

\begin{proof}
    We will first show that there exists $z \in N$ such that $\varphi(t, z) \in N$
    for all $t > 0$. Consider a horizontal disk $h: \overline B(0, 1) \to N$. Observe that by Lemma \ref{disk_image_lemma} it follows easily that
    for all $t > 0$ we have
     $$B(t) = \{ \overline{x} \in \overline B(0,1) :  \varphi(s, h(\overline{x})) \in N \ \mathrm{for}\ s\in [0,t]\} \neq \emptyset.$$
     Since $B(t) \subset \overline{B}(0, 1)$ is a compact set for every $t > 0$ and $B(s) \subset B(t)$ for $s \geq t$, we have
     $$
     \bigcap_{t\geq 0} B(t) \neq \emptyset
     $$
     and the claim follows.
     \par
     Thus we can pick $z \in N$ which has a full forward trajectory in $N$.
     Now consider $z_0$ in the $\omega$-limit set of $z$. We want to show that for each $t > 0$ we have $\varphi(t, z_0) = z_0$. Let us fix $t > 0$ and define
     $$\xi_t: \expr{0, \infty} \ni s \mapsto \mathcal{Q}\expr{\varphi\expr{s, z} - \varphi\expr{s, \varphi \expr{t, z}}} \in \mathbb{R}.$$
     If $z$ is not a fixed point of $\varphi$ (in case it is, the proof of existence is already finished), then function $\xi_t$ is increasing by the cone conditions and since $N$ is compact, $\xi_t$ is also bounded.
     This means that $\lim_{s \to \infty} \xi_t(s) = A$ for some $A \in \mathbb{R}$. Since $z_0 \in \omega(z)$,
     there exists a sequence
     $\set{t_n}_{n \in \mathbb{N}}$ such that $\varphi(t_n, z) \to z_0$.
     We have $\lim_{n \to \infty} \xi_t(t_n) = A$, so $\mathcal{Q}(z_0 - \varphi(t, z_0)) = A$.
     If we had $\varphi(t, z_0) \neq z_0$, then for $r > 0$ we would have by the cone conditions
     $$
     \mathcal{Q}(\varphi\expr{r, z_0} - \varphi\expr{t + r, z_0}) > A,
     $$
     but on the other hand it would hold that $\mathcal{Q}(\varphi\expr{r, z_0} - \varphi\expr{t + r, z_0}) = \lim_{n \to \infty} \xi_t(t_n + r) = A,$ a contradiction.
     \par
     That $z_0$ is a unique fixed point in $N$ is an immediate consequence of the cone conditions.
\end{proof}

We almost immediately get the stable manifold theorem.

\begin{theorem} \label{thm_stable_manifold}
    If $N$ is an $\varepsilon$-isolating cuboid for a continuous local
    semiflow $\varphi$ which satisfies the cone conditions on $N$,
    then $\varphi$ has a unique fixed point $z_0 \in N$ and there exists a vertical disk
    $w_{s} : P_n^\bot N \to N$ such that
    $W^s_N\expr{z_0} = w_s(P_n^\bot N).$
    Moreover, if a point $z \in N$ has  \textcolor{black}{a}  full forward trajectory in $N$,
    then $z \in W^s_N\expr{z_0}$.
\end{theorem}
\begin{proof}
    Consider $y \in P_n^\bot N$ and consider a horizontal disk $h: P_n \to N$
    given by $h(x) = (x, y)$. By Lemma \ref{lemma_fixed_point} we get that there
    exists unique $x^*$ such that $\varphi\expr{t, h\expr{x^*} = \expr{x^*, y}} \to z_0$.
    So we can define $w_s: P_n^\bot N \ni y \mapsto (x^*, y)$. It remains to show
    that $w_s$ is a vertical disk. Let $y_1, y_2 \in P_n^\bot N$ and assume that
    $\mathcal{Q}\expr{w_s\expr{y_1} - w_s\expr{y_2}} > 0$.
    It implies that there exists a horizontal disk passing both through $w_s\expr{y_1}$ and
    $w_s\expr{y_2}$, which is a contradiction with Lemma \ref{lemma_fixed_point},
    because both of those points converge to $z_0$.
\end{proof}

Proof of the unstable manifold theorem requires a bit more additional work.

\begin{theorem} \label{thm_unstable_manifold}
    If $N$ is an $\varepsilon$-isolating cuboid for a continuous local
    semiflow $\varphi$ which satisfies the cone conditions on $N$,
    then $\varphi$ has a unique fixed point $z_0 \in N$ and there exists a horizontal disk
    $w_{u} : P_n N \to N$ such that
    $W^u_N\expr{z_0} = w_u(P_n N).$
    Moreover, if a point $z \in N$ has a full backward trajectory in $N$,
    then $z \in W^u_N\expr{z_0}$.
\end{theorem}
\begin{proof}
    Consider a horizontal disk $h: \overline{B}(0, 1) \to N$. Fix $x \in \overline{B}(0, 1)$.
    Due to Lemma \ref{disk_image_lemma}, there exists a sequence $\set{a_k}_{k \in \mathbb{N}} \subset N$ such that
    for all $k \in \mathbb{N}$ we have $P_n a_k = x$ and there exists $\overline{x}_k$ for which
    $a_k = \varphi(k t(\varepsilon), h(\overline{x}_k))$.
    We will apply the Cantor's diagonal
    argument to prove that there is a convergent subsequence of $\set{a_k}_{k \in \mathbb{N}}$ such that its limit has an infinite backward orbit.
    \par
    Consider families of sequences $\set{k_i^l}_{i \in \mathbb{N}}$, $\set{b_{k_i^l}^l}_{i \in \mathbb{N}}$,
    where $l \in \mathbb{N}$,
    such that
    \begin{itemize}
        \item[(i)] for each $l_1 > l_2$ the sequence $\set{k_i^{l_1}}_{i \in \mathbb{N}}$ is a subsequence of $\set{k_i^{l_2}}_{i \in \mathbb{N}}$,
        \item[(ii)] sequence $\set{a_{k_i^1}}_{i \in \mathbb{N}}$ is convergent,
        % \item[(iii)] for each $i \in \mathbb{N}$ and $t \in (0, lt(\varepsilon)]$ there exists $\overline{z} \in N$ such that
        % $\varphi(t, \overline{z}) = a_{k_i^l}$,
        \item[(iii)] for each $l \in \mathbb{N}$ the sequence
        $\set{b_{k_i^l}^l}_{i \in \mathbb{N}}$ is convergent and we have
        $\varphi(lt(\varepsilon), b_{k_i^l}^l) = a_{k_i^l}$.
    \end{itemize}
    Existence of such a family follows easily from the definition of $\set{a_n}_{n \in \mathbb{N}}$ and the compactness of $N$.
    \par
    % Now for $n \in \mathbb{N}$ define $w_i := a_{k_1^i}$. By (i) an (ii) it is a convergent
    % sequence, so we can denote $w_x := \lim_{i \to \infty} w_i$.
    We define $w_x := \lim_{i \to \infty} a_{k_i^1}$ (this limit exists by (ii)) and $o_n := \lim_{i \to \infty} b_{k_i^n}^n$ for $n \in \mathbb{N}$ (those limits exist by (iii)). For $n \in \mathbb{N}$, since by (iii) we have
    $\varphi(n t\expr{\varepsilon}, b_{k_i^n}^n) = a_{k_i^n}$, it follows by (i) that
    $\varphi\expr{nt\expr{\varepsilon}, o_n} = w_x$,
    thus $w_x$ has an infinite backward orbit in $N$. By a reasoning analogous to the one
    in the proof of Lemma \ref{lemma_fixed_point}, we can prove that $\lim_{n \to \infty} o_n = z_0$.
    \par
    Thus we have proved that for every $x \in \overline{B}(0, 1)$ there exists $w_x \in W^u_N\expr{z_0}$
    such that $x = P_n w_x$.
    To show that it is unique assume
    that $w_x' \in N$ is such that $P_n w_x' = x$ and that there exists a sequence $\set{o_{k}'}_{k \in \mathbb{N}}$ in $N$ such that $\varphi\expr{nt\expr{\varepsilon}, o_n'} = w_x'$ for $n \in \mathbb{N}$. By the cone conditions we have
    $$0 \geq \mathcal{Q}(w_x - w_x') \geq \lim_{k \to \infty} \mathcal{Q}(o_k - o_k') = \mathcal{Q}(z_0 - z_0) = 0,$$
    so $w_x' = w_x$.
    % To show that it is unique assume for the sake of contradiction
    % that there exists some other $w_x' \in N$ such that $P_n w_x' = x$ for which there exists a sequence $\set{o_{k}'}_{k \in \mathbb{N}}$ in $N$ such that $\varphi\expr{nt\expr{\varepsilon}, o_n'} = w_x'$ for $n \in \mathbb{N}$. We have
    % $$0 > \mathcal{Q}(w_x - w_x') > \lim_{k \to \infty} \mathcal{Q}(o_k - o_k') = \mathcal{Q}(z_0 - z_0) = 0,$$
    % a contradiction.
    \par
    Thus we can define a map $w_u: \overline{B}(0, 1) \ni x \mapsto w_x \in N$. It remains
    to show that $w$ is a horizontal disk, but the proof is analogous to the proof
    of uniqueness of $w_x$ for each $x \in \overline{B}(0, 1)$.
\end{proof}

Now we define what we consider to be a hyperbolic point for a differential equation
in the Hilbert space.
Notion of hyperbolicity we use differs from the functional-analytic one based on the spectrum. Connecting
such definitions with what they mean for the dynamics is usually complicated in infinite dimension \textcolor{black}{ (see \cite{FSV, SS} for the approach connecting spectrum with the dynamics)}. Instead, we focus on what is the
expected behavior of the solutions of the system near the hyperbolic fixed points.
% and then we will use a stronger notion of the cone conditions to give sufficient conditions for a point to be hyperbolic
\begin{definition}
    Let $\varphi$ be a local semiflow on $N$. We say that
    $S \subset N$ is \emph{forward-invariant} (\emph{backward-invariant})
    \emph{relative to $N$} if for every every orbit $z$ in $N$ with $z(0) \in S$ and for all $t > 0$
    such that
    $z\expr{\interval{0, t}} \subset N$ ($t < 0$ such that $z\expr{\interval{-t, 0}} \subset N$)
    we have  $z\expr{\interval{0, t}}  \subset S$ ($z\expr{\interval{-t, 0}} \subset S$).
\end{definition}

\begin{definition} \label{def_hyperbolicity}
    We say that a fixed point $z_0 \in N$ is \emph{hyperbolic}, if there exist constants
    $A, B, \omega > 0$ such that the sets
    $z_0 + Q^+, z_0 + Q^-$ are nonempty and respectively forward- and backward-invariant relative to $N$
    and there holds
    \begin{align}
        &\forall z \in z_0 + Q^+ \ \norm{P_n \expr{\varphi\expr{t, z} - z_0}} \geq A e^{\omega t} \norm{P_n \expr{z - z_0}}, \textnormal{if } \varphi\expr{\interval{0, t}, z} \in N, \label{eqn_hyp1} \\
        &\forall z \in z_0 + Q^- \ \norm{P_n^\perp  \expr{\varphi\expr{t, z} - z_0}} \leq B e^{-\omega t} \norm{P_n^\perp  \expr{z - z_0}}, \textnormal{if } \varphi\expr{\interval{0, t}, z} \in Q^-. \label{eqn_hyp2}
    \end{align}
\end{definition}

For our definition of the hyperbolicity to be reasonable, it must
imply the standard notion of hyperbolicity based on the eigenvalues
when the dimension is finite.
The following theorem shows that it is indeed true.

\begin{theorem} \label{thm_hyperbolic_fin_dim}
    Let $f = (f_x, f_y): \mathbb{R}^N \to \mathbb{R}^N$ be a $\mathcal{C}^2$ function. If in the equation
    \begin{align*}
        x' &= f_x\expr{x, y}, \\
        y' &= f_y\expr{x, y}
    \end{align*}
    the point $p$ is a hyperbolic fixed point in the sense of Definition \ref{def_hyperbolicity} (with respect to
    the cones given by $\mathcal{Q}(x, y) = \norm{x}^2 - \norm{y}^2$), then it is hyperbolic
    in the standard sense (i.e. all eigenvalues of the linearization at the fixed point have nonzero real part).
\end{theorem}
    \begin{proof}
    Without any loss of generality assume that $p = 0$.
    Let $z = (x, y) \in Q^+$. Fix $t > 0$.
    Denote $L = Df(0)$ and $S(z) = (S_x, S_y)(z) := \varphi(t, z)$.
    Since the solutions of $z' = Lz$ are of the form $z(t) = e^{Lt}z(0)$,
    we have that $z(t) = DS(0)z(0)$ and $x(t) = DS_x(0)z(0), y(t) = DS_y(0)z(0)$.
    Thus to provide bounds for the solutions of the linearized equation we
    only need to focus on $DS(0)$. Now, due to the identity
    $$S_x(z) = S_x(z) - S_x(0) = \int_{0}^{1} DS_x\expr{sz} \mathrm{d}s \; z,$$
    by (\ref{eqn_hyp1}) we have for $\alpha > 0$
    $$\norm{\int_{0}^{1} DS_x\expr{s\alpha z} \mathrm{d}s \alpha z} \geq Ae^{\omega t}\norm{\alpha z}.$$

    Cancelling $\alpha$ on both sides and using the fact that when $\alpha \to 0^+$
    $$\int_{0}^{1} DS_x\expr{s\alpha z} \mathrm{d}s \to DS_x(0),$$
    we get that $\norm{DS_x(0)z} \geq Ae^{\omega t} \norm{z}$, which proves
    that (\ref{eqn_hyp1}) is satisfied for the linearized equation.
    We can analogically prove that (\ref{eqn_hyp2}) is also satisfied.

    % Denoting $L = Df(0, 0)$, we get
    % \begin{align*}
    %     \norm{P_n \varphi\expr{t, z}} \leq \norm{P_n e^{tL}z} + \norm{P_n \expr{e^{tL}z - \varphi\expr{t, z}}}.
    % \end{align*}
    % Let $\epsilon > 0$. If $z \in Q_\delta^+$ is close enogh to $0$, we get by using the Taylor expansion that
    % \begin{align*}
    %     \norm{P_n\expr{e^{tL}z - \varphi\expr{t, z}}} =
    %     \norm{P_n \expr{e^{tL} - \frac{\partial \varphi}{\partial z}\expr{t, z}}z} + O(\norm{x}^2)
    %     < \epsilon \norm{x}.
    % \end{align*}
    % In consequence,
    % \begin{align*}
    %     \norm{P_n e^{tL}z} > \expr{Ae^{\omega t} - \epsilon}\norm{x}.
    % \end{align*}
    % We see that taking $\epsilon < \frac{A}{2}$, we get that the condition (\ref{eqn_hyp1}) holds for the linearized
    % equation. We can analogously prove that the condtion (\ref{eqn_hyp2}) also holds.
    \par
    For a linear equation, conditions (\ref{eqn_hyp1}, \ref{eqn_hyp2}) obviously imply that $Q^+$
    is forward-invariant.
    Thus by Lemma \ref{disk_image_lemma} we see that denoting $D_y := \set{\expr{x, y} \mid x \in \mathbb{R}^n}$
    the set $\varphi\expr{t, D_y}$ contains the image of a horizontal disk for each $y \in \mathbb{R}^{N-n}$ and $t > 0$.
    From this and (\ref{eqn_hyp1}, \ref{eqn_hyp2}) we can easily deduce that for each $y$ there exists $z \in D_y$
    such that $\lim_{t \to \infty} \varphi(t, z) = 0$, thus its stable manifold is a graph over the $y$ directions.
    Since the equation is linear, this means that the stable manifold is a linear subspace of the same dimension as $y$.
    \par
    Reversing time we also get that the unstable manifold is a linear subspace of the same dimension as $x$.
    Thus the stable and unstable subspaces span the entire $\mathbb{R}^N$ and it is thus well known from the standard
    theory that no eigenvalue of $Df(0)$ can be purely imaginary.
    % It is clear that for a linear equation conditions (\ref{eqn_hyp1}, \ref{eqn_hyp2}) imply
    % the cone conditions. Again, in the case of a linear equation the resulting unstable and stable manifolds
    % are linear subspaces which span the entire $\mathbb{R}^N$, thus every
    % invariant subspace must be a subspace of one of those manifolds. It is thus clear that we cannot
    % have a purely imaginary eigenvalue, as by the standard theory on its corresponding subspace the conditions
    % (\ref{eqn_hyp1}, \ref{eqn_hyp2}) would be violated.
\end{proof}

One could ask why we have not used Theorem \ref{thm_stable_manifold} in the proof above,
but seemingly we have reproved it. This is because we only proved the positive
cone forward-invariance for the linearization; this invariance allowed only
to prove that the stable manifold is a graph over $y$ (not a vertical disk),
but because of the linearity it already allows us to establish the dimension
of the stable subspace.
\par
To prove that a fixed point is hyperbolic, we need to use a stronger notion of the cone conditions.
They give us exponential expansion (contraction) in cones not only in terms of $\mathcal{Q}$
like the cone conditions, but also in terms of the norm.

\begin{definition}
Let $N \subset H$ be a $\varepsilon-$isolating cuboid and assume $\varphi$ is a local semiflow on $N$.
Let
$$
\mathcal{Q}:H \to \mathbb{R}
$$
be a map defined by $\mathcal{Q}(u) = \alpha(c_u(P_nu)) - \beta(P_n^\bot u)$,
where $\alpha:\mathbb{R}^n\to \mathbb{R}$ and $\beta:Y_n \to \mathbb{R}$
are continuous positively defined quadratic forms.
For $\delta > 0$ define
\begin{align*}
    \mathcal{Q}_\delta(z) &:= \mathcal{Q}(z) + \delta \norm {P_n z}^2, \\
    \mathcal{Q}^\delta(z) &:= \mathcal{Q}(z) - \delta \norm {P_n^\bot z}^2, \\
\end{align*}
We say that a semiflow $\varphi$ satisfies the \emph{strong cone conditions}
on $N$ if for every $z_1, z_2 \in N, z_1 \neq z_2$,
the derivative $\frac{d}{dt} \mathcal{Q}\expr{\varphi\expr{t, z_1} - \varphi\expr{t, z_2}}|_{t=0}$
exists and there exists $\lambda > 0$ such
that for $z_1 \neq z_2, \lambda_* \in \interval{-\lambda, \lambda}$ and for $\delta > 0$
small enough there holds
\begin{align}
    \frac{d}{dt} \mathcal{Q}_\delta\expr{\varphi\expr{t, z_1} - \varphi\expr{t, z_2}}|_{t=0}
        &> \lambda_* \mathcal{Q}_\delta(z_1 - z_2), \label{eqn_strong_cc_1} \\
    \frac{d}{dt} \mathcal{Q}^\delta\expr{\varphi\expr{t, z_1} - \varphi\expr{t, z_2}} |_{t=0}
        &> \lambda_* \mathcal{Q}^\delta(z_1 - z_2), \label{eqn_strong_cc_2}
    \end{align}
% \begin{align}
% \frac{d}{dt} \expr{ \mathcal{Q}\expr{\varphi\expr{t, z_1} - \varphi\expr{t, z_2}} \pm \delta \norm{P_n \varphi\expr{t, z_1} - P_n \varphi\expr{t, z_2}}^2}|_{t=0} &> \lambda_* Q(z_1 - z_2), \label{eqn_strong_cc_1} \\
% \frac{d}{dt} \expr{ \mathcal{Q}\expr{\varphi\expr{t, z_1} - \varphi\expr{t, z_2}} \pm \delta \norm{P_n^\bot \varphi\expr{t, z_1} - P_n^\bot \varphi\expr{t, z_2}}^2}|_{t=0} &> \lambda_* Q(z_1 - z_2), \label{eqn_strong_cc_2}
% \end{align}
% &\frac{d}{dt} \mathcal{Q}\expr{\varphi\expr{t, z_1} - \varphi\expr{t, z_2}}|_{t=0} > \lambda_* Q(z_1 - z_2), \label{eqn_strong_cc_1} \\
% &\frac{d}{dt} \expr{ \mathcal{Q}\expr{\varphi\expr{t, z_1} - \varphi\expr{t, z_2}} + \delta \norm{\pi_x \varphi\expr{t, z_1} - \pi_x \varphi\expr{t, z_2}}^2}|_{t=0} > 0, \label{eqn_strong_cc_2} \\
% &\frac{d}{dt} \expr{ \mathcal{Q}\expr{\varphi\expr{t, z_1} - \varphi\expr{t, z_2}} - \delta \norm{\pi_y \varphi\expr{t, z_1} - \pi_y \varphi\expr{t, z_2}}^2 } |_{t=0} > 0. \label{eqn_strong_cc_3}

% &\frac{d}{dt} \mathcal{Q}\expr{\varphi\expr{t, z_1} - \varphi\expr{t, z_2}}|_{t=0} > \lambda_* \norm{z_1 - z_2}^2, \label{eqn_strong_cc_1} \\
% &\alpha(z) \leq \alpha_* \norm{P_n z}^2, \label{eqn_strong_cc_2} \\
% &\beta(z) \leq \beta_* \norm{P_n^\perp  z}^2.
\end{definition}

\begin{theorem}\label{thm_hyperbolicity_characterization}
    Let $N$ be $\varepsilon$-isolating cuboid for a continuous local
    semiflow $\varphi$ and assume that $\varphi$ satisfies the strong cone conditions on $N$.
    Then there exists a unique fixed point $z_0 \in N$ and moreover this fixed point is hyperbolic.
\end{theorem}
\begin{proof}
    The fixed point $z_0$ exists by Lemma \ref{lemma_fixed_point}.
    Without loss of generality, we assume that $z_0 = 0 \in N$.
    \par
    The cones $Q^+, Q^-$  and respectively forward-
    and backward-invariant relative to $N$ respectively by
    (\ref{eqn_strong_cc_1}, \ref{eqn_strong_cc_2}) with $\delta = 0$.
    \par
    Now let $z \in Q^+$ and
    assume $t > 0$ is such that $\varphi\expr{\interval{0, t}, z} \subset N$
    (hence, by the relative forward-invariance, $\varphi\expr{\interval{0, t}, z} \subset Q^+$).
    Since $\alpha$ is a continuous
    positively defined quadratic form, there exists $\alpha_* > 0$ such that for all $z \in H$
    we have $\alpha(P_n z) \leq \alpha_* \norm{P_n z}^2$.
    Observe that $\mathcal{Q}(z) > 0$ implies $\mathcal{Q}_\delta(z) > \delta \norm{P_n z}^2$
    and that $\mathcal{Q}_\delta(z) = \alpha(P_n z)
    + \delta\norm{P_n z}^2 - \beta(P_n^\bot z) \leq \expr{\alpha_* + \delta}\norm{P_n z}^2$.
    Thus using (\ref{eqn_strong_cc_1}) with $\lambda_* = \lambda$
    we have for $\delta > 0$ small enough
    $\norm{P_n \varphi\expr{t, z}}^2 \geq \frac{1}{\alpha_* + \delta}
    \mathcal{Q}_\delta\expr{\varphi\expr{t, z}} > \frac{1}{\alpha_* + \delta} e^{\lambda t} \mathcal{Q}_\delta(z)
    \geq \frac{\delta}{\alpha_* + \delta} e^{\lambda t} \norm{P_n z}^2$.
    \par
    Now let $z \in Q^-$ and
    assume $t > 0$ is such that $\varphi\expr{\interval{0, t}, z} \subset Q^-$. Since $\beta$
    is a continuous
    positively defined quadratic form, there exists $\beta_* > 0$ such that for all $z \in H$ we have
    $\beta(P_n^\perp z) \leq \beta_* \norm{P_n^\perp z}^2$.
    Observe that $\mathcal{Q}(z) < 0$ implies $\mathcal{Q}^\delta(z) < -\delta \norm{P_n^\bot z}^2$
    and that $\mathcal{Q}^\delta(z) = \alpha(P_n z)
    - \beta(P_n^\bot z) - \delta\norm{P_n^\bot z}^2 \geq -\expr{\beta_* + \delta}\norm{P_n^\bot z}^2$.
    Thus using (\ref{eqn_strong_cc_2})
    with $\lambda_* = -\lambda$ we have
     $\norm{P_n^\perp  \varphi\expr{t, z}}^2 < -\frac{1}{\delta} \mathcal{Q^\delta}\expr{\varphi\expr{t, z}} <
     -\frac{1}{\delta} e^{-\lambda t} \mathcal{Q^\delta}(z) \leq
     \frac{1}{\delta} e^{-\lambda t} (\beta_* + \delta) \norm{P_n^\perp  z}^2$.
\end{proof}

\section{ Pitchfork bifurcation} \label{sec_pitchfork_models}

\subsection{Definition of \textcolor{black}{a} pitchfork bifurcation}

Consider a family of equations
\begin{equation}\label{eq_pde_mu} \frac{du}{dt} = F_{\nu}(u)\end{equation}
for $\nu \in \mathbb{R}$.
In this subsection we give the definition of what we
understand by \textcolor{black}{a} pitchfork bifurcation. In the next subsections, we discuss simple
ODE models to motivate the definition and 
%\textbf{cos takiego:}
 \textcolor{black}{to present our approach to proving the bifurcation in this simplified setting.}
.

\begin{definition}
    Let $H$ be a Hilbert space, $S \subset H$ and assume that $\varphi$ is a local semiflow.
    We say that $I \subset S$ is the \emph{maximal invariant set} in $S$ if
    $\overline z \in S$ belongs to $I$ if and only if
    \begin{itemize}
        \item there is a full solution passing through $\overline z$,
        \item for every full solution $z$ passing through $\overline{z}$
        we have $z(t) \in I$ for $t \in \mathbb{R}$.
    \end{itemize}
\end{definition}

\begin{definition}
\label{def_pitchfork}
Assume that $R$ is an $\varepsilon$-isolating cuboid  \textcolor{black}{for}
the local semiflow $\varphi_\nu$ associated with (\ref{eq_pde_mu})
for each $\nu \in \interval{\nu_-, \nu_+},$ where $\nu_- < 0 < \nu_+$.
We say that (\ref{eq_pde_mu}) undergoes  \textcolor{black}{a}  \emph{pitchfork bifurcation} on
the interval $\interval{\nu_-, \nu_+},$ where $\nu_- < 0 < \nu_+$, if
\begin{itemize}
    \item[(P1)] There exists $u_0 \in R$ such that $u_0$ is a fixed point for all $\nu \in \interval{\nu_-, \nu_+}$. Moreover, $u_0$ is hyperbolic in sense of Definition \ref{def_hyperbolicity} for $\nu \neq 0$.
    \item[(P2)] For $\nu \in \interval{\nu_-, 0}$ the set $\set{u_0}$ is the maximal invariant set in $R$.
    \item[(P3)] For $\nu \in (0, \nu_+]$ there exist hyperbolic fixed points $u_-^\nu, u_+^\nu \in R$
    different from $u_0$
    such that there exist heteroclinic connections from $u_0$ to $u_-^\nu$ and $u_+^\nu$.
    Moreover, the set consisting of $u_0$, of the points in the heteroclinic connections from $u_0$ to $u_\pm^\nu$
    and of $u_\pm^\nu$ is the maximal invariant set in $R$.
\end{itemize}

We refer to the set $R$ as to  \textcolor{black}{a}
\emph{set isolating all dynamics near the bifurcation}. Moreover, if $R$ is forward-invariant
(i.e. when there are no unstable directions),
then we call it  \textcolor{black}{a}  \emph{set trapping all dynamics near bifurcation}. The second
case corresponds to the case where all directions except the bifurcation direction
are stable.
\end{definition}

\subsection{Bifurcation model without unstable directions} \label{subsection_pitchfork_model}

We start off by considering a simple pitchfork bifurcation ODE model to illustrate what
we understand by the full description of a bifurcation and to present a technique
which will later be used in the proof of \textcolor{black}{a} pitchfork bifurcation in the Kuramoto--Sivashinsky
system.
\par
% In this section try to present things as self-contained as we can and only refer
% to more formal results and definitions presented later.
Throughout this
section  \textcolor{black}{whenever we write about a}  hyperbolic fixed point we mean hyperbolicity in the sense of Definition
\ref{def_hyperbolicity}. By Theorem \ref{thm_hyperbolic_fin_dim} it also implies
the standard hyperbolicity for ODEs (and in this model we could check
this hyperbolicity straightforwardly), but our goal here is to present techniques
which are later used in the infinite-dimensional case.
\par
Consider \textcolor{black}{the folowing planar} equation with parameter $\nu \in \mathbb{R}$.
\begin{align} \label{eq_model}
    \begin{split}
        x' &= x(\nu - x^2) + x^3y + xy^2, \\
        y' &= -y + x^2y + x^4,
    \end{split}
\end{align}
and denote by $\varphi_\nu$ the associated dynamical system and by $F_\nu = (f, g)_\nu$ the right-hand side of (\ref{eq_model}).
We see that as $\nu$ passes through $0$, the zero solution loses its stability, hence a bifurcation of the $(0, 0)$ fixed point
occurs. \textcolor{black}{ To understand}  the character of this bifurcation consider terms $(\tilde f, \tilde g)_\nu(x, y) = (x(\nu - x^2), -y)$.
Close to the origin, these terms dominate remaining terms of $F_\nu$ and thus they drive the entire dynamics \textcolor{black}{near the origin}.
This is why given the symmetry $(\tilde f, \tilde g)_\nu(-x, y) = (-\tilde f, \tilde g)_\nu(x, y)$ \textcolor{black}{ of (\ref{eq_model}) },
the bifurcation we expect is indeed \textcolor{black}{a} pitchfork bifurcation (also known as the symmetry breaking bifurcation).
\par
Thus for $\nu > 0$ we expect two new fixed points to be born. Approximately, those fixed points satisfy
$(\tilde f, \tilde g)_\nu(x, y) = 0$, which means that $x = \pm \sqrt{\nu}$ and $y = 0$.
We want to guess, for a fixed $\nu > 0$, a set which encompasses  all interesting local
dynamics after bifurcation. Obviously, it needs to contain the approximate fixed points
plus some spare space. This leads us to \textcolor{black}{the following} guess
$$B_\nu = \interval{-\sqrt{2\nu}, \sqrt{2\nu}} \times \interval{-\nu^\kappa, \nu^\kappa}, \kappa > 0.$$
Moreover, we would like this set to be forward-invariant, so that we can hope (and later prove)
that we do not leave anything interesting outside of it. Thus the vector field on the \textcolor{black}{boundary} 
should point in \textcolor{black}{the direction of the set interior}. This gives us the conditions
\begin{itemize}
    \item if $x = \pm \sqrt{2\nu}$, then, since $\nu - x^2$ is negative, we need to have $\abs{x(\nu - x^2)} > \abs{x^3y + xy^2}$,
    \item if $y = \pm \nu^\kappa$, then we need to have $\abs{y} > \abs{x^2y + x^4}$.
\end{itemize}
We can see that those conditions are easily satisfied for any $\textcolor{black}{2 > }\kappa > \frac{1}{2}$    if $x, y$ are close
enough to $0$. We pick $\kappa = 1$.
\par
As indicated before, we proceed to what we call \textbf{Step I -- existence of the set isolating the bifurcation}.
Namely, we show that there exists $\nu_+ > 0$ such that for $\nu \in (0, \nu_+]$ \textcolor{black}{the whole} interesting dynamics near the bifurcation happens
on $B_\nu$. We do it by showing that there exists a set $R$ independent of $\nu$ such that for
each $\nu \in (0, \nu_+]$ every point in $R$ flows eventually into $B_\nu$. More precisely, the following lemma holds.

\begin{lemma}[\textbf{Step I -- existence of the set isolating the bifurcation, situation after the bifurcation}] \label{lemma_big_isolation_model}
    There exists $\nu_+ > 0$ such that defining
    $R := B_{\nu_+}$ there holds that for each $\nu \in (0, \nu_+]$ there
    exists a time $T$ such that for all $t > T$ we have $\varphi\expr{t, R} \subset B_\nu$.
\end{lemma}
\begin{proof}
Fix $\nu_+ > 0, \nu \in (0, \nu_+]$. We will prove that if $\nu_+$ is small enough,
then for  any $\theta \in \interval{\nu, \nu_+}$ there exists a time $T > 0$ such that for all $t > T$
we have $\varphi_\nu(t, B_\theta) \subset B_{\bar\theta}$, where $\bar\theta := \max \set{\frac{\theta}{2}, \nu}$.
If $\theta$ is small and if
$(x, y) \in B_\theta \setminus B_{\max \set{\frac{\theta}{2}, \nu}}$,
we have
\begin{itemize}
    \item if $\sqrt{2}\theta^{\frac{1}{2}} \geq x \geq \sqrt{2}\bar\theta^{\frac{1}{2}}$, then
    $x' = x(\nu - x^2) + x^3y + xy^2 < \bar\theta(\nu - 2\bar\theta) + r(x, y)$,
    where $r(x, y) = O(\theta^{\frac{5}{2}})$, so $x' < 0$, because
    $\nu - 2\bar\theta < -\bar\theta \leq -\frac{1}{2} \theta < 0$,
    \item if $\sqrt{2}\theta^{\frac{1}{2}} \geq -x \geq \sqrt{2}\bar\theta^{\frac{1}{2}}$, then analogously
    $x' > 0$,
    \item if $\theta \geq y \geq \bar\theta$, then $y' = -y + x^2y + x^4 < 0$, as $x^2y + x^4 = O(\theta^2)$ and $-y \leq -\frac{\theta}{2} < 0$,
    \item if $\theta \geq -y \geq \bar\theta$, then analogously $y' > 0$.
\end{itemize}
\par
Define $\theta_0 = \nu_+$ and $\theta_k = \max\set{\frac{\theta_{k - 1}}{2}, \nu}$. For some
$n \in \mathbb{N}$ we have $\theta_n = \nu$.
By what we have shown above, if $\nu_+$ is small enough, then there exists a sequence of times $t_0, t_1, \dots, t_{n - 1} > 0$ such that
$\varphi\expr{t_i, B_{\theta_i}} \subset B_{\theta_{i + 1}}$, $i = 0, \ldots, n - 1$.
This completes the proof for $\nu \geq 0$
\end{proof}
We can easily see that using  \textcolor{black}{a} similar reasoning as above we can prove that for $\nu \leq 0$ the origin is attracting.
Indeed, if $\nu_+$ and $\abs{\nu_-}$ are small enough, then for $\nu \in \interval{\nu_-, 0}$
and for all $\theta \in \interval{0, \nu_+}$  by nearly identical calculations as above there exists
a time $T > 0$ such that for all $t > T$
we have $\varphi_\nu(t, B_\theta) \subset B_{\frac{\theta}{2}}$. This means that every
point in $R$ flows into $B_\xi$ for arbitrarily small $\xi$. This gives us the following lemma.

\begin{lemma}[\textbf{Step I -- existence of the set isolating the bifurcation, situation before the bifurcation}] \label{lemma_model_origin_stable}
    There exist $\nu_- < 0 < \nu_+$ such that, setting $R := B_{\nu_+}$, for all $\nu \in \interval{\nu_-, 0}$
    and all $\xi \in (0, \nu_+]$ there exists $T > 0$ such that for $t > T$ we have $\varphi_\nu\expr{t, R} \subset B_\xi$.
    Consequently, for $\nu \in \interval{\nu_-, 0}$ the origin $(0, 0)$ is an attracting fixed point.
\end{lemma}

We now proceed to the \textbf{Step II -- hyperbolicity of the origin before the bifurcation}.
To describe the behavior before the bifurcation completely, we wish to establish that for $\nu < 0$ the origin
is a hyperbolic fixed point. We prove it using the logarithmic norms described in Appendix \ref{app_lognorms}.
Let us note that if we prove that for the logarithmic norm $\mu_2$ given by the $l_2$ norm
we have for some $l < 0$ that $\mu_2\expr{DF(z)} \leq l$ for all $z \in B_{\abs{\nu}}$, then
by Lemma \ref{lemma_lognorm_exp_bound} for $z_1, z_2 \in B_{\abs{\nu}}$ we have
$\norm{\varphi\expr{t, z_1} - \varphi\expr{t, z_2}} < e^{lt} \norm{\varphi\expr{0, z_1} - \varphi\expr{0, z_2}}$
for all $t > 0$, since  $B_{\abs{\nu}}$ is forward-invariant. It gives us another proof that
the origin is attracting for $\nu < 0$ (but it would not work for $\nu = 0$; as expected,
because this point is not hyperbolic). Moreover this  \textcolor{black}{clearly shows}  that the origin is hyperbolic in the sense of
Definition \ref{def_hyperbolicity}, with $\mathcal{Q}\expr{z} = -\norm{z}^2$.

\begin{lemma}[\textbf{Step II -- hyperbolicity of the origin before the bifurcation}]
    There exist $\nu_- < 0$ such that for $\nu \in [\nu_-, 0)$ the origin is a hyperbolic attracting fixed point.
\end{lemma}
\begin{proof}
    Fix $\nu_- < 0$ and $\nu \in [\nu_-, 0)$. We will use Theorem \ref{thm_lognorm_bounds} to bound
    $\mu_2$ on $B_{\abs{\nu}}$. Let $z = (x, y) \in B_{\abs{\nu}}$ and denote $A = DF_\nu(z)$. We have
    $$
    A =  \begin{pmatrix}
      \nu - 3x^2 + 3x^2y + y^2 & x^3 + 2xy \\
      2xy + 4x^3 &-1 + x^2
    \end{pmatrix},
    $$
    thus
    $$
    \frac{A + A^T}{2} =  \begin{pmatrix}
        \nu - 3x^2 + 3x^2y + y^2 & \frac{2xy + 4x^3 + x^3 + 2xy}{2} \\
        \frac{2xy + 4x^3 + x^3 + 2xy}{2} &-1 + x^2
      \end{pmatrix}.
    $$
    This is a symmetric matrix so it has two eigenvalues $\lambda_1, \lambda_2$
    and by the Gershgorin theorem we have for some $K, L > 0$
    that  \textcolor{black}{when}  $\abs{\nu}$ \textcolor{black}{is}  small enough there holds
    \begin{align*}
        \lambda_1 &< \nu - 3x^2 + 3x^2y + y^2 + \abs{\frac{2xy + 4x^3 + x^3 + 2xy}{2}} \\
                  &< \nu - 3x^2 + \abs{3x^2y + y^2} + \abs{\frac{2xy + 4x^3 + x^3 + 2xy}{2}} \\
                  &< \nu - 3x^2 + K \abs{\nu}^{\frac{3}{2}} \\
                  &< \nu + K \abs{\nu}^{\frac{3}{2}} < 0,
    \end{align*}
    and
    \begin{align*}
        \lambda_2 &< -1 + x^2 + \abs{\frac{2xy + 4x^3 + x^3 + 2xy}{2}} < -1 + L\abs{\nu} < 0.
    \end{align*}
    It is also easy to see that the bounds above do not depend on $z$,  thus
    the claim follows by Theorem \ref{thm_lognorm_bounds} and Lemma \ref{lemma_lognorm_exp_bound}.
\end{proof}

\begin{figure}[h!]
\begin{center}
    \includegraphics{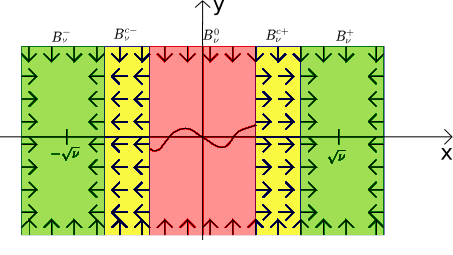}
    \caption{Illustration of the behavior after bifurcation.}
    \label{fig_proof_p4}
\end{center}
\end{figure}

Figure \ref{fig_proof_p4} illustrates the behavior after bifurcation.
We will use the following strategy of the proof (to prove that the connection from zero to the
stable equilibrium to the right
exists; proof of the existence of the other one is analogous).
It is split into three steps.
% Proof of (P3) is the most complicated part. Its proof is sketched in
% Figure \ref{fig_proof_p4}. The idea is to show that
\begin{itemize}
    \item \textbf{Step IIIa -- hyperbolicity of the origin after the bifurcation.}
    The proof is realized in Lemma \ref{lemma_zero_hyperbolic}, where we show that on the red set $B_\nu^0$
    the strong cone conditions are satisfied, hence $0$ is hyperbolic and there exists
    a point $z_* \in B_\nu^{c+} \cap W^{u}_{B_\nu^0}(0)$, where $B_\nu^{c+}$ is the
    blue set to the right in the picture.
    \item \textbf{Step IIIb -- connection sets.}
    The proof is realized in Lemma \ref{lemma_connection_sets_model}, where we show that on the blue set
    to the right $B_\nu^{c+}$
    we have $x' > 0$, thus there exists a time $t_1 > 0$
    such that $\varphi(t_1, z_*) \in B_\nu^+$ (the green set to the right in the picture).
    \item \textbf{Step IIIc -- hyperbolicity of the born fixed points.}
    The proof is realized in Lemma \ref{lemma_neg_lognorm_after_bif_model}, where we show that the green set to the right $B_\nu^+$
    is forward-invariant and for all $z \in B_\nu^+$ we have $\mu_2(DF(z)) < l < 0$,
    thus as we have discussed above there exists a unique fixed point
    $z^+ \in B_\nu^+$ and moreover this fixed point is hyperbolic and we have $\lim_{t \to \infty} \varphi(t, \varphi(t_1, z_*)) = z^+$.
\end{itemize}

\begin{lemma}[\textbf{Step IIIa -- hyperbolicity of the origin after the bifurcation}] \label{lemma_zero_hyperbolic}
    Let $$B_\nu^0 := \interval{-\frac{\sqrt{\nu}}{2}, \frac{\sqrt{\nu}}{2}} \times \interval{-\nu, \nu}.$$
    There exists $\nu_+ > 0$ such that for (\ref{eq_model})
    for $\nu \in (0, \nu_+]$ the origin is hyperbolic and
    the unstable manifold $W^u_{B_\nu^0}(0)$
    is the image of a horizontal disk in $B_\nu^0$.
\end{lemma}
\begin{proof}
    Fix $\nu_+ > 0$ and $\nu \in (0, \nu_+]$.
    In Lemma \ref{lemma_connection_sets_model} we will show that $B_\nu^0$ is an $\varepsilon$-isolating
    cuboid for sufficiently small $\nu$. Thus to use Theorems \ref{thm_unstable_manifold} and \ref{thm_hyperbolicity_characterization}
    it is enough to show that the strong cone conditions are satisfied with the quadratic form given by the matrix
    $$
    Q =  \begin{pmatrix}
      1 & 0 \\
      0 & -1
    \end{pmatrix}.
    $$
    Denote $Q_{\pm\delta} = Q \pm \delta \Id$.
    To prove that the strong cone conditions hold, we
    will show that or $z_1, z_2 \in B_\nu^0$ there exists $\lambda > 0$ such that
    for $\delta > 0$ small enough we have
    \begin{align*}
        \frac{d}{dt}\expr{\expr{\varphi\expr{t, z_1} - \varphi\expr{t, z_2}}^T Q_{\pm\delta} \expr{\varphi\expr{t, z_1} - \varphi\expr{t, z_2}}} \Big\vert_{t = 0} > \lambda \expr{z_1 - z_2}^2
    \end{align*}
    We have
    \begin{align*}
        \frac{d}{dt} &\expr{\expr{\varphi\expr{t, z_1} - \varphi\expr{t, z_2}}^T Q_{\pm\delta}\expr{\varphi\expr{t, z_1} - \varphi\expr{t, z_2}}} \Big\vert_{t = 0}  \\
        &= \expr{F_\nu\expr{z_1} - F_\nu\expr{z_2}}^T Q_{\pm\delta} \expr{z_1 - z_2} + \expr{z_1 - z_2}^T Q_{\pm\delta} \expr{F_\nu\expr{z_1} - F_\nu\expr{z_2}},
    \end{align*}
    We now want to factorize $z_1 - z_2$.
    Given that denoting $\overline{DF}_\nu := \int_{0}^1 DF_{\nu}\expr{tz_1 + (1 - t)z_2} \mathrm{d}t$
    we have
    $$F_\nu\expr{z_1} - F_\nu\expr{z_2} = \overline{DF}_\nu \expr{z_1 - z_2},$$
    thus
    \begin{align*}
        \frac{d}{dt} &\expr{\expr{\varphi\expr{t, z_1} - \varphi\expr{t, z_2}}^T Q_{\pm\delta}\expr{\varphi\expr{t, z_1} - \varphi\expr{t, z_2}}} \Big\vert_{t = 0} \\
        &=
            (z_1 - z_2)^T \expr{\overline{DF}_\nu^T Q_{\pm\delta} + Q_{\pm\delta} \overline{DF}_\nu} (z_1 - z_2)
    \end{align*}

    Thus it is sufficient to show that $A = \overline{DF}_\nu^T Q_{\pm\delta} + Q_{\pm\delta} \overline{DF}_\nu$
    is positive definite. We see that the eigenvalues of $A$ are bounded from below by
    \begin{align*}
        \lambda_1 &> 2 (1 - \delta)\inf_{z \in B_\nu^0} \frac{\partial f_\nu}{\partial x}\expr{z}
            - (1 + \delta) \sup_{z \in B_\nu^0} \expr{ \abs{\frac{\partial f_\nu}{\partial y}\expr{z} + \frac{\partial g_\nu}{\partial x}\expr{z}} } \\
                  &> 2 (1 - \delta) \nu - (1 + \delta) \sup_{z \in B_\nu^0} \expr{3x^2 + \abs{x^3 + 2xy + \frac{5}{2}x^3}} > \frac{\nu}{4} - O(\nu^{\frac{3}{2}}) > 0, \\
            \lambda_2 &> 2 (1 - \delta)\inf_{z \in B_\nu^0} \expr{-\frac{\partial g_\nu}{\partial y}\expr{z}}
            - (1 + \delta) \sup_{z \in B_\nu^0} \abs{\frac{\partial g_\nu}{\partial x}\expr{z} + \frac{\partial f_\nu}{\partial y}\expr{z}} \\
                &> 2 (1 - \delta) - \expr{1 + \delta} \sup_{z \in B_\nu^0} \expr{ x^2 + \abs{x^3 + 2xy + \frac{5}{2}x^3}} > 2(1 - \delta) - O(\nu) > 0,
    \end{align*}
    if $\nu, \delta$ are small enough.
    %  we have $\nu - 3x^2 \geq \frac{\nu}{4}$ on $B_{\nu}^0$,
    % thus for some $\beta > 0$ it holds that
    % $$\nu - 3x^2 + 3x^2y + y^2 - \abs{x^3 + 2xy} > \frac{\nu}{4} - \beta\nu^\frac{3}{2} > 0,$$
    % since $\nu > 0$. We also obviously have $-(-1 + x^2) - \abs{2xy + 4x^3} > 0.$
    In Lemma \ref{lemma_connection_sets_model} we will also show that $B_\nu^0$ is an $\varepsilon$-isolating
    cuboid. This completes the proof.
\end{proof}

We now prove that on Figure \ref{fig_proof_p4} everything from the blue region
eventually flows into the green region. It will also show the missing part of
the proof above, namely that $B_\nu^0$ is $\varepsilon$-isolating cuboid for $\nu > 0$.

\begin{lemma}[\textbf{Step IIIb -- connection sets}]  \label{lemma_connection_sets_model}
    There exists $\nu_+ > 0$ such that for each $0 < \nu \leq \nu_+$
    we have
    \begin{align}\label{eq_xprime_in_model}
        x' &> 0, \\
        x' &< 0,
    \end{align}
    respectively on the sets
    \begin{align*}
        B_\nu^{c+} &:= \interval{\frac{\sqrt{\nu}}{2}, \sqrt{\frac{\nu}{2}}} \times \interval{-\nu, \nu}, \\
        B_\nu^{c-} &:= \interval{-\sqrt{\frac{\nu}{2}}, -\frac{\sqrt{\nu}}{2}} \times \interval{-\nu, \nu}. \\
    \end{align*}

\end{lemma}
\begin{proof}
    Let $(x, y) \in B_\nu^{c+}$ (proof for $B_\nu^{c-}$ is analogous). If $\nu > 0$ is small enough,
    then we have
    $$
    x' = x(\nu - x^2) + x^3y + xy^2 > \frac{3}{8} \nu^{\frac{3}{2}} - \expr{\frac{1}{2\sqrt{2}} + \frac{1}{\sqrt{2}}}\nu^{\frac{5}{3}} > 0,
    $$
    since $\nu > 0$.
\end{proof}

We now proceed to establishing the basins of attraction in the green ($B_\nu^\pm$,
see Figure \ref{fig_proof_p4})
neighborhoods of approximate fixed points  $(\pm \sqrt{\nu}, 0)$ of (\ref{eq_model}).

\begin{lemma} [\textbf{Step IIIc -- hyperbolicity of the born fixed points}]  \label{lemma_neg_lognorm_after_bif_model}
    There exists $\nu_+ > 0$ such that for $0 < \nu \leq \nu_+$ there exist
    unique attracting fixed points $u_\pm \in B^\pm_\nu$, where
    \begin{align*}
        B^+_\nu &:= \interval{\sqrt{\frac{\nu}{2}}, \sqrt{2\nu}} \times \interval{-\nu, \nu},\\
        B^-_\nu &:= \interval{-\sqrt{2\nu}, -\sqrt{\frac{\nu}{2}}} \times \interval{-\nu, \nu}.\\
    \end{align*}
\end{lemma}
\begin{proof}
The sets
$B^\pm_\nu$ are also easily seen to be forward-invariant.
\par
We proceed to  \textcolor{black}{prove}  that the $l_2$ logarithmic norm is negative.
We have
$$
\frac{DF_\nu\expr{x, y}^T + DF_\nu\expr{x, y}}{2} =  \begin{pmatrix}
  \nu - 3x^2 + 3x^2y + y^2 & \frac{x^3 + 2xy + 2xy + 4x^3}{2} \\
  \frac{x^3 + 2xy + 2xy + 4x^3}{2} &-1 + x^2
\end{pmatrix}.
$$
We show that the eigenvalues $\lambda_{1, 2}$ of this matrix are negative.
If $\nu_+$ is small enough, then for $\nu \in [0, \nu_+]$, for $(x, y) \in B_\nu^\pm$
we have $\nu - 3x^2 \leq - \frac{\nu}{2}$
$$
\lambda_1 \leq \nu - 3x^2 + 3x^2y + y^2 + \abs{\frac{x^3 + 2xy + 2xy + 4x^3}{2}} \leq -\frac{\nu}{2} + O(\nu^\frac{3}{2}) < 0,
$$
as $\nu > 0$. We also obviously have $\lambda_2 \leq -1 + x^2 + \abs{\frac{x^3 + 2xy + 2xy + 4x^3}{2}} < 0$.
Since $B_\nu^\pm$ are forward-invariant, the conclusion follows by Lemma \ref{lemma_lognorm_exp_bound}.
\end{proof}

We are finally ready to prove the bifurcation theorem for our model.

\begin{theorem} \label{thm_bif_model}
    There exist $\nu_- < 0 < \nu_+$ such that \textcolor{black}{a} pitchfork bifurcation occurs
    in (\ref{eq_model}) on the interval $\interval{\nu_-, \nu_+}$.
\end{theorem}
\begin{proof}
    By Lemmas \ref{lemma_big_isolation_model}, \ref{lemma_zero_hyperbolic}
    it remains to show that (P3) is satisfied. Let $\nu_+ > 0$ be small enough and
    let $0 < \nu \leq \nu_+$. By Lemma \ref{lemma_zero_hyperbolic} there
    exists such a point $z_* = (x_*, y_*) \in B_\nu^{c+} \cap B_\nu^0$ which has a
    backward orbit to $0$. By Lemma \ref{lemma_connection_sets_model} there exists
    $T > 0$ such that $\varphi_\nu(T, z_*) \in B_\nu^+$ and by Lemma \ref{lemma_neg_lognorm_after_bif_model}
    everything in $B_\nu^+$ is attracted to the fixed point $u_+^\nu$. This
    gives us the required heteroclinic connection.
    \par
    Proof of the connection to $u_-^\nu$ is analogous.
    \par
    It remains to show that the fixed points
    together with the heteroclinic connections constitute the maximal invariant set
    in $R$. By Lemma \ref{lemma_big_isolation}, every point in $R$ flows into $B_\nu$.
    Let us now discuss the backward trajectories of points in $B_\nu$.
    By Lemma \ref{lemma_zero_hyperbolic} and Theorem \ref{thm_unstable_manifold},
    only points on the unstable manifold have  full backward trajectories in $B_\nu^0$. Moreover,
    since $B_\nu^0$ is an isolating cuboid, for every other point in $B_\nu^0$ its
    backward trajectory leaves $B_\nu$. Now, in $B_\nu^{c+} \cup B_\nu^+$ every point's backward
    trajectory goes into $B_\nu^0$ or leaves $B_\nu$, so only points with   
    full backward trajectories
    are the ones whose backward trajectories pass through $z_*$. The conclusion about
    the maximal invariance easily follows.

\end{proof}

Observe that the terms of higher order in $F_\nu$
(i.e. $f_\nu$ terms of form $x^k y^l$ for $k > 3$ or $l > 3$
and terms of form
$x^k y^l$ for $k > 2$ and $l \geq 1$ or  for $k \geq 4$  in $g_\nu$) would not
change our proof of the bifurcation in any meaningful way.  
\par
However, terms like $xy$  in $f_\nu$
or $y^2$  in equation for $g_\nu$ would change the proof, as in $DF_\nu$ they give rise respectively to $x$
and $2y$. Because of such terms, $\nu - 3x^2$ would no longer be dominating in $DF_\nu$ and the
presented reasoning would not work. Such terms do appear in the Kuramoto--Sivashinsky equation written
in the Fourier basis.
\par
We could deal with those terms by taking different shapes of sets, but it is a very unwieldy
approach. Instead, in Section \ref{sec_normal_forms}, using the normal form
theory we present coordinate changes
which remove those terms, introducing instead other terms of higher order, which
in turn pose no problem to our approach.
\par
We did not consider the term $y^2$ (in the equation for $x'$) in this model, \textcolor{black}{as it is not symmetrical. However, let us mention that it would not be obstacle for our method}. 
It could be removed
by the change of coordinates, but it would also suffice to take different shape
of sets in the proof -- we elaborate on this approach in the proof in the model given below,
which  \textcolor{black}{is useful there} because of the term $yz$ (\textcolor{black}{observe that this term is symmetrical}), 
where $z$ will be another direction
different from the bifurcation direction. Although those terms could also be removed
in a finite-dimensional case, in the Kuramoto--Sivashinsky there are infinitely many
terms of such form, which is why we will present an argument with the changed shape of sets.

% ===========================================================================
\subsection{Bifurcation model with an unstable direction}

\subsubsection{Statement of the model and outline of the proof}

Now we will briefly discuss the changes which occur when we add an unstable direction.
Consider the following equation

\begin{align} \label{eqn_model_unstable}
    \begin{split}
    x' &= x\expr{\nu - x^2} + x^3\expr{y + z} + x\expr{y^2 + z^2} + yz, \\
    y' &= y + x^2\expr{y + z} + x^4 + yz, \\
    z' &= -z + x^2\expr{y + z} + x^4 + yz.
    \end{split}
\end{align}

The proof in this case proceeds similarly to the proof in the previous model, although
there are some differences which we ought to stress before we proceed to the details of the
proof.
\par
In the previous model when proving \textbf{Step I} we descended
through the family of the forward-invariant sets closer and closer to the origin.
Of course, because \textcolor{black}{of the presence } of the unstable direction it is impossible now.
This is why now we flow through the descending family of the isolating cuboids.
The details of this flow
are only a bit more complicated than in the previous model.
Before the bifurcation we get every point in $R$ either
leaves $R$ or flows to the origin in the limit. Similarly,
after the bifurcation every point either flows out of $R$
or flows into $B_\nu$. We know even more -- if a point leaves
$B_\nu$, then it necessarily leaves $R$. We thus see that
as before to describe the dynamics after the bifurcation completely,
we just need to describe the dynamics on $B_\nu$.
\par
When it comes to \textbf{Step II}, i.e. proving the hyperbolicity before the bifurcation,
we can no longer use the logarithmic norms because  \textcolor{black}{ of the presence }  of the unstable direction.
But it does not change much, as verifying the strong cone conditions relies on
very similar calculations (unsurprising, because negativity of the $l_2$ logarithmic norm
is in a special case of the cone conditions).
\par
\begin{figure}[h!]
\begin{center}
    \includegraphics[scale=0.5]{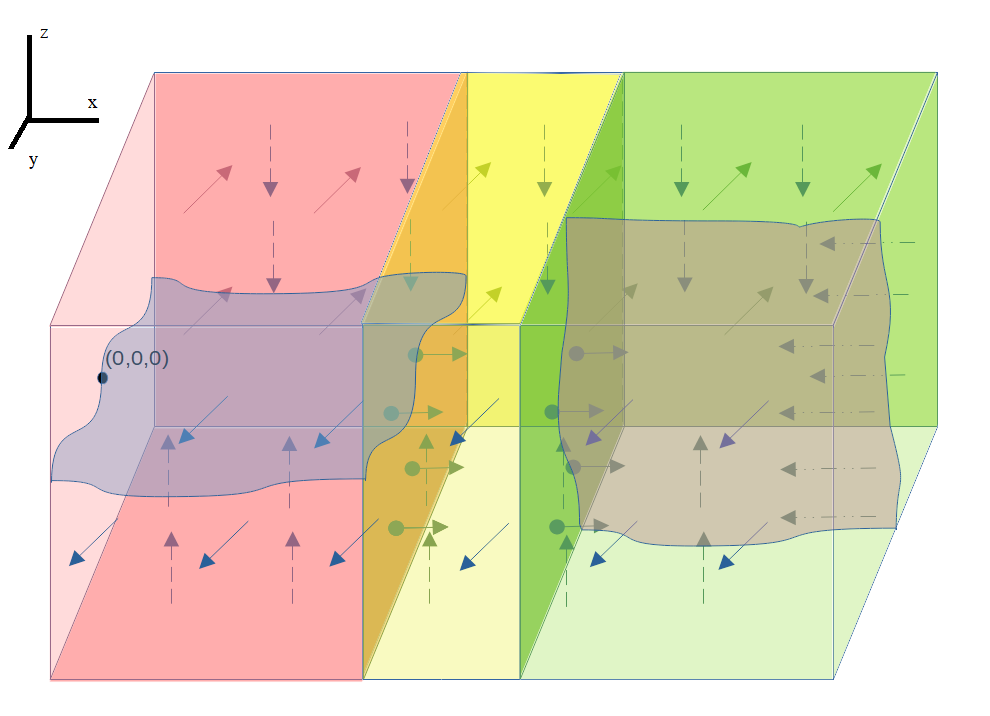}
    \caption{Illustration of the proof of the heteroclinic connection when there is an unstable direction.
    The red set is the set $B_\nu^0$, in which we verify the cone conditions to get
    the unstable manifold (the blue surface). In the yellow set $B_{\nu}^{c+}$ we have growth in the
    bifurcation direction. The green set is $B_\nu^0$, in which we brown surface is the stable manifold
    which we again get by the cone conditions.}
    \label{fig_proof_unstable}
\end{center}
\end{figure}

The illustration of the proof   \textcolor{black}{ for parameter values past the bifurcation locus, i.e. $\nu >0$}  is in the
Figure \ref{fig_proof_unstable} (one side, as previously the other one
is analogous).
$B_\nu^0, B_\nu^{c+}, B_\nu^{+}$ refer respectively to
the red, yellow and green set in the picture. The set
$B_\nu^{ct+}$ is the sum $B_\nu^{c+} \cup B_\nu^{+}$.
\par
As before, in \textbf{Step IIIa} we use the cone conditions to establish that the unstable
manifold of the origin is a graph.
However, after the bifurcation, this manifold close to the origin
(blue surface in the picture)
is a graph over not only the bifurcation direction $x$, but also the over the
additional unstable direction $y$, hence it is two-dimensional. So while in the previous
model its intersection with the exit set of the isolating cuboid close to the
origin was a point (for each of the 'branches' of the manifold),
now it is a graph over the unstable direction $y$ (see the intersection
$W^u_{B_\nu^0}(0) \cap B_\nu^{c+}$, i.e. of the blue surface with the yellow set).
\par
Previously in \textbf{Step IIIb} we could propagate the point from the unstable manifold
on the exit set close to the fixed point that was born simply by propagating
the entire exit set. Now it does not work, because not every
point from this exit set flows close to the target fixed point.
However, we prove that on $B_\nu^{ct+}$ the cones given by
the matrix
\begin{align*}
    \begin{pmatrix}
        &-1 &0 &0 \\
        &0 &1 &0 \\
        &0 &0 &-1
      \end{pmatrix}.
\end{align*}
are forward-invariant. Together with the fact that on
$B_\nu^{c+}$ we have $x' > 0$ it gives us that
there is a subset $W' \subset W^u_{B_\nu^0}(0) \cap B_\nu^{c+}$
such that for some time $T > 0$ we have that $H := \varphi_\mu(T, W') \subset B_{\nu}^{+}$
is a graph over the $y$ direction. As we will see in the moment, $H$ intersects with the
stable manifold.
% (blue subset of $W^u_{B_\nu^0}(0) \cap B_\nu^{c+}$; its images in the forward propagation are also blue)
% is a graph over $y$ in $B_\nu^+$
% % (blue subset of $S_\nu^+$)
% , let us for now call this
% graph $H$.
\par
Next, we proceed to the \textbf{Step IIIb} in which we verify the cone conditions on the set $B_\nu^+$ and conclude that
\begin{itemize}
    \item there exists a hyperbolic fixed point $z_\nu^+ \in B_\nu^+$,
    \item there exists exactly one point on $H$ which converges to $z_\nu^+$ forward in time;
          it gives us the required heteroclinic connection,
    \item the unstable (stable) manifold $W^u_{B_\nu^+}(z_\nu^+)$($W^s_{B_\nu^+}(z_\nu^+)$)
          is a graph over the $y$ direction (over the $x, z$ directions);
          we need those facts to prove that the heteroclinic connections constitute the
          maximal invariant set.
\end{itemize}

\subsubsection{Details of the proof}
We proceed to the details of the proof. Denote
$$B_\nu := \interval{-\sqrt{2\nu}, \sqrt{2\nu}} \times \interval{-\nu^{\frac{3}{2}}, \nu^{\frac{3}{2}}}^2.$$
We also denote the form given by the matrix
\begin{align*}Q_1 =
    \begin{pmatrix}
      &1 &0 &0 \\
      &0 &1 &0 \\
      &0 &0 &-1
    \end{pmatrix}
\end{align*}
by $\mathcal{Q}_1$ and the form given by the matrix
\begin{align*}
    Q_2 =
    \begin{pmatrix}
      &-1 &0 &0 \\
      &0 &1 &0 \\
      &0 &0 &-1
    \end{pmatrix}.
\end{align*}
by $\mathcal{Q}_2$.

We now prove \textcolor{black}{that there exists a set isolating the whole dynamics near bifurcation (see the end of Definition \ref{def_pitchfork})}.
We start with the case after the bifurcation.

\begin{lemma}[\textbf{Step I -- existence of the set isolating the bifurcation, situation after the bifurcation}] \label{lemma_big_isolation_model_unstable}
    There exists $\nu_+ > 0$ such that for $\nu \in (0, \nu_+]$ we have
    that for all $u \in R := B_{\nu_+}$ either $u$ leaves $R$
    or there exists $t > 0$ such that  $\varphi\expr{t, u} \in B_\nu$. Moreover,
    if a point $u \in B_\nu$ leaves $B_\nu$, then it also leaves $R$.
\end{lemma}
\begin{proof}
Fix $\nu > 0$ and $\theta \in \interval{\nu, \nu_+}$,
where $\nu_+$ is small enough. Defining $\theta' = \max\set{\frac{\theta}{2}, \nu}$
we will show that every point $u \in S_\theta \setminus S_{\theta'}$
flows into $S_{\theta'}$ or leaves $R$. First let us fix $(x, y, z) = u \in S_\theta \setminus S_{\theta'}$
and look at the vector field. We can prove analogously as in Lemma \ref{lemma_big_isolation_model}
that $\sgn x' = -\sgn x, \sgn z' = -\sgn z$ and $\sgn y' = \sgn y$. Now observe
that it easily implies if there exists sufficiently big $T > 0$ such that for all $t \in \interval{0, T}$
we have $\abs{\pi_y \varphi(t, u)} < \theta'^{\frac{3}{2}}$, then we have $\varphi\expr{T, u} \in S_{\theta'}$.
Thus to prove the lemma it is enough to show that if for some $t > 0$ we have $\abs{\pi_y \varphi(t, u)} \geq \theta'^{\frac{3}{2}}$,
then we leave $R$. In this case we leave $S_\theta$ because
$\abs{y}$ increases, so we can never go into $S_{\theta'}$. But once we leave
then we have a point in $S_{\min\set{\nu_+, 2\theta}} \setminus S_\theta$ such that
for it $\abs{y} = \theta^{\frac{3}{2}}$, so for analogous reasons it leaves $S_{\min\set{\nu_+, 2\theta}}$
and so on until we leave $R$.
\end{proof}
% \begin{proof}
% Fix $\nu > 0$. We first prove that if for $\theta > \nu$ we have
% that $(x, z) \in \pi_{\expr{x, z}} S_{theta}$
% and $\abs{y} = \theta^{\frac{3}{2}}$,
% then for some $t > 0$ we have $\varphi\expr{t, (x, y, z)} \in S_{\min{2\theta, \nu_+}} \setminus S_{\theta}$.
% Observe that it is simply a reversal of a technique by which we have shown in the proof
% of Lemma \ref{lemma_big_isolation_model} that everything flows close to the origin. We calculate
% $$y' = $$
% We can prove by the same reasoning as in the proof
% of Lemma $\ref{lemma_big_isolation_model}$ that if $u \in R$ stays for a sufficiently long
% time $T > 0$ in $R$, then for $t > T$ we have $\pi_{\expr{x, z}} \varphi\expr{t, u} \in \pi_{\expr{x, z}} B_\nu$,
% where $\pi_{\cdot}$ denotes a projection
% on the given directions. It is also easy
% to verify that for $u \in R$ such that $\pi_{\expr{x, z}} u \in \pi_{\expr{x, z}} B_\nu$
% and $u \not\in B_\nu$ we have $\sgn y' = \sgn y,$ so those points leave $R$.
% \end{proof}

Slightly modifying the proof above, we get the following lemma about the
situation before the bifurcation.

\begin{lemma}[\textbf{Step I -- existence of the set isolating the bifurcation, situation before the bifurcation}] \label{lemma_zero_attracting_unstable}
    There exist $\nu_- < 0 < \nu_+$ such that for each $\nu \in \interval{\nu_-, 0}$
    we have that for all $u \in R := B_{\nu_+}$ either $u$ leaves $R$
    or we have $\varphi\expr{t, u} \to 0$.
\end{lemma}

Now we state the hyperbolicity result. Compared to the proof of the Lemma \ref{lemma_zero_hyperbolic},
we need to use the cone conditions instead of the logarithmic norms in the proof of the
hyperbolicity before the bifurcation (because of the additional unstable direction).
Nevertheless, the computations are almost the same, so we skip them.

\begin{lemma}[\textbf{Step II and IIIa -- hyperbolicity of the origin before and after the bifurcation}]  \label{lemma_zero_hyperbolic_unstable}
    There exists $\nu_- < 0 < \nu_+$ such that for (\ref{eqn_model_unstable}) the origin
    is hyperbolic for $\nu \in \interval{\nu_-, \nu_+} \setminus \set{0}$.
    Moreover, for $\nu \in (0, \nu_+]$ the unstable manifold $W^u_{B_\nu^0}(0)$, where
    $B_\nu^0 := \interval{-\frac{\sqrt{\nu}}{2}, \frac{\sqrt{\nu}}{2}} \times \interval{-\nu^\frac{3}{2}, \nu^\frac{3}{2}}^2$, is
    the image of a horizontal disk in $B_\nu^0$ with respect to the cones given by
    $\mathcal{Q}_1$.
\end{lemma}

Now we proceed to  \textcolor{black}{prove}  that although not all points from the set $B_\nu^0 \cap B_\nu^{c\pm+}$
flow into $B_\nu^\pm$, we can nevertheless say that the image of some subset
of $W^u_{B_\nu^0} \cap B_\nu^{c\pm}$ is a graph over the unstable direction $y$ in the set $B_\nu^+$.

\begin{lemma}[\textbf{Step IIIb -- connection sets, cone invariance}] \label{lemma_model_cone_field}
    If $\nu > 0$ is small enough, then the sets
    \begin{align*}
        B_\nu^{ct+} &:= \interval{\frac{\sqrt{\nu}}{2}, \sqrt{2\nu}} \times \interval{-\nu^{\frac{3}{2}}, \nu^{\frac{3}{2}}}^2, \\
        B_\nu^{ct-} &:= \interval{-\sqrt{2\nu}, \frac{\sqrt{\nu}}{2}} \times \interval{-\nu^{\frac{3}{2}}, \nu^{\frac{3}{2}}}^2.
    \end{align*}
    are $\varepsilon-$isolating cuboids. Moreover, if $h$ is a horizontal disk with respect to the cones
    $\mathcal{Q}_2$ in $B_\nu^{ct\pm}$, then $\varphi_\mu\expr{t, h\expr{\pi_{y}B_\nu^{ct\pm}}}$ contains the image of another horizontal
    disk with respect to the cones $\mathcal{Q}_2$ for all $t > 0$.
\end{lemma}
\begin{proof}
    We will verify the assumptions of Lemma \ref{disk_image_lemma}.
    First we check that those sets are indeed $\varepsilon$-isolating cuboids. Let $(x, y, z) \in B_\nu^{ct+}.$
    Then if $\nu$ is small enough we can find $d > 0$ such that for all $(x, y, z)$ with
    $y = \nu^{\frac{3}{2}}$ (case $y = -\nu^{\frac{3}{2}}$ is analogous) we have
    \begin{align*}
        y' &= y + x^2\expr{y + z} + x^4 + yz \geq \\
            &\geq \nu^{\frac{3}{2}} - O\expr{\nu^\frac{5}{2} + \nu^2 + \nu^3} > 0.
    \end{align*}
    Thus we see that we can pick $\varepsilon$ such that the condition the condition (I2) holds.
    % We also see that $y'$ is easily seen to be uniformly separated from zero for such points
    \par
    Proper behavior on $x, z$ directions can be checked by the computations as in the previous model.
    \par
    Now we verify the remaining assumption of the Lemma \ref{disk_image_lemma},
    i.e. that the positive cones are invariant. We do this by showing
    that if $\mathcal{Q}_2\expr{\varphi\expr{t, u_1} - \varphi\expr{t, u_2}} > 0$,
    then $\frac{d}{dt} \mathcal{Q}_2\expr{\varphi\expr{t, u_1} - \varphi\expr{t, u_2}}\Big\vert_{t = 0} > 0$.
    \par
    As in the proof of Lemma \ref{lemma_zero_hyperbolic}, we get that
    for $u_1, u_2 \in B_\nu^{ct+}$ and defining $A := \overline{DF}_\nu^T Q_{2} + Q_{2} \overline{DF}_\nu$

    $$ \frac{d}{dt} \expr{\varphi\expr{t, u_1} - \varphi\expr{t, u_2}}^T Q_2 \expr{\varphi\expr{t, u_1} - \varphi\expr{t, u_2}}
                = (u_1 - u_2)^T A (u_1 - u_2). $$

    Denote $u_1 - u_2 =: \expr{x, y, z}$.
    Then we can symbolically write
    $$
    (u_1 - u_2)^T A (u_1 - u_2) = \sum_{a, b \in \set{x, y, z}} A_{ab} ab.
    $$

    Now what we will show is that $A_{yy}, A_{zz} > 0$ and that $A_{yy}$ dominates
    $A_{xx}$ (which on some part of the set is negative)
    and all off-diagonal $A$ terms. Then, considering that $\mathcal{Q}_2\expr{u_1 - u_2} > 0$
    is equivalent to $y^2 > x^2 + z^2$, we have for sufficiently small $\nu > 0$
    (ignoring a non-negative term $A_{zz}z^2$).

    \begin{align*}
        \sum_{a, b \in \set{x, y, z}} A_{ab} ab >
            \expr{A_{yy} - \abs{A_{xx}} - \sum_{a, b \in \set{x, y, z}, a \neq b} \abs{A_{ab}}}y^2 > 0.
    \end{align*}

    Thus we calculate
    \begin{align*}
        A_{xx} &\geq \nu - 3x^2 - \sup_{u \in B_\nu^{ct+}} \expr{3x^2\expr{y + z} + y^2 + z^2}
                \geq -5 \nu - O(\nu^{\frac{5}{2}}), \\
        A_{yy} &\geq 1 - \sup_{u \in B_\nu^{ct+}} \expr{x^2 + z} = 1 - O(\nu) > 0, \\
        A_{zz} &\geq 1 - \sup_{u \in B_\nu^{ct+}} \expr{x^2 + y} = 1 - O(\nu) > 0.
    \end{align*}

    We see that we can similarly show that the off-diagonal terms are all at worst $O(\nu)$
    and thus we have shown what we wanted to. We also see that the bounds do not depend on
    the choice of $u_1, u_2$.
\end{proof}

The proof of result about the growth in the bifurcation direction almost does not change, so we skip it.

\begin{lemma}[\textbf{Step IIIb -- connection sets, growth in the bifurcation direction}] \label{lemma_connection_sets_model_unstable}
    We have
    \begin{align*}
        x' &> 0, \\
        x' &< 0,
    \end{align*}
    respectively on the sets
    \begin{align*}
        B_\nu^{c+} &:= \interval{\frac{\sqrt{\nu}}{2}, \frac{\sqrt{\nu}}{\sqrt{2}}} \times \interval{-\nu^{\frac{3}{2}}, \nu^{\frac{3}{2}}}^2, \\
        B_\nu^{c-} &:= \interval{-\frac{\sqrt{\nu}}{\sqrt{2}}, \frac{\sqrt{\nu}}{2}} \times \interval{-\nu^{\frac{3}{2}}, \nu^{\frac{3}{2}}}^2.
    \end{align*}
\end{lemma}

The following lemma states that the new fixed points are indeed born and
that they are hyperbolic. Moreover, their stable manifolds are graphs over
$x, z$ directions. To prove this it is enough to verify the cone conditions
and use the Theorem \ref{thm_stable_manifold}, so the computations are again skipped.

\begin{lemma}[\textbf{Step IIIc -- hyperbolicity of the born fixed points}] \label{lemma_neg_lognorm_after_bif_model_unstable}
    There exists $\nu_+ > 0$ such that for $0 < \nu \leq \nu_+$ there exist
    unique fixed points $u^\pm_\nu \in B^\pm_\nu$, where
    \begin{align*}
        B^+_\nu &:= \interval{\sqrt{\frac{\nu}{2}}, \sqrt{2\nu}} \times \interval{-\nu^{\frac{3}{2}}, \nu^{\frac{3}{2}}}^2,\\
        B^-_\nu &:= \interval{-\sqrt{2\nu}, -\sqrt{\frac{\nu}{2}}} \times \interval{-\nu^{\frac{3}{2}}, \nu^{\frac{3}{2}}}^2.
    \end{align*}
    Moreover, those fixed points are hyperbolic with respect to the cones given by
    $\mathcal{Q}_2$ and $W^u_{B^\pm_\nu}\expr{u^\pm_\nu}$ ($W^s_{B^\pm_\nu}\expr{u^\pm_\nu}$)
    are the images of horizontal (vertical) disks.
\end{lemma}

We are ready to prove that the bifurcation occurs in the model with an unstable direction.

\begin{theorem} \label{thm_bif_model_unstable}
    There exist $\nu_- < 0 < \nu_+$ such that \textcolor{black}{a} pitchfork bifurcation occurs
    in (\ref{eqn_model_unstable}) on the interval $\interval{\nu_-, \nu_+}$.
\end{theorem}
\begin{proof}
    Observe that by Lemma \ref{lemma_zero_hyperbolic_unstable} the unstable manifold $W^u_{B_\nu^0}$ is the image of a horizontal disk
    with respect to cones $\mathcal{Q}_1$, thus the set $W^* := W^u_{B_\nu^0} \cap B_\nu^{c+}$ is clearly the
    image of a horizontal disk in $B_\nu^{ct+}$ with respect to the cones $\mathcal{Q}_2$.
    Thus due to the Lemmas \ref{lemma_model_cone_field} and \ref{lemma_connection_sets_model_unstable}  we get that there exists
    a time $T > 0$ such that $\varphi\expr{T, W^*}$ is the image of some horizontal disk in
    $B_{\nu}^+$, which we denote by $h$. Thus by Lemmas \ref{lemma_neg_lognorm_after_bif_model_unstable}
    and \ref{lemma_fixed_point}  we get that for some
    $u \in h(\pi_y B_\nu^+)$ we have $\lim_{t \to \infty} \varphi(t, u) = z_+$.
    \par
    It remains to show that we have obtained the maximal invariant set.
    By Lemma \ref{lemma_big_isolation_model_unstable}, we only need to focus on $B_\nu$ for $\nu > 0$
    (case $\nu \leq 0$ being obvious).
    \par
    Observe that
    by Lemma \ref{thm_unstable_manifold} and by the fact that $B_\nu^0$ is an $\varepsilon$-isolating
    cuboid, every point in $B_\nu^0$ leaves $B_\nu$ backward in time or is on the unstable manifold,
    so only points in
    $B_\nu^0$ which can be in the maximal invariant set are the points from $W^u_{S_\nu^0}\expr{0}$.
    \par
    Points in the sets $B_\nu^\pm$ backward in time either leave $B_\nu$ or go into $B_\nu^0$.
    Thus only those points which have some point from $W^u_{S_\nu^0}\expr{0}$ on their backward
    trajectory can be in the maximal invariant set.
    \par
    Now consider $u \in B_\nu^\pm$. Then this point's backward trajectory may
    \begin{itemize}
        \item leave $B_\nu$, so it is not in the maximal invariant set.
        \item stay in $B_\nu$. Then $u \in W^u_{B_\nu^\pm}(u_\nu^\pm)$.
              Now, forward in time $u$ either leaves $B_\nu$, so it
              is not in the maximal invariant set, or stays in $B_\nu^\pm$.
              So for $u$ to be in the maximal invariant set, we would
              also need to have $u \in W^s_{B_\nu^\pm}(u_\nu^\pm)$,
              thus by Lemma \ref{lemma_neg_lognorm_after_bif_model_unstable}
              we have $u = u_\nu^\pm$, because horizontal and vertical
              disk can intersect only in one point.
        \item flow into $B_\nu^{c+}$, so all we have said before applies.
    \end{itemize}
\end{proof}
%
% \begin{itemize}
%     \item In $B_\nu^0 := \interval{-\frac{\sqrt{\nu}}{\sqrt{2}}, \frac{\sqrt{\nu}}{\sqrt{2}}} \times \interval{-\nu^{\frac{3}{2}}, \nu^{\frac{3}{2}}}^2$
%     the cone conditions for the flow are satisfied with quadratic form $\mathcal{Q}_1$
%     given by the matrix
%     $$
%     \begin{pmatrix}
%       &1 &0 &0 \\
%       &0 &1 &0 \\
%       &0 &0 &-1
%     \end{pmatrix}.
%     $$
%     This set is also an isolating cuboid, thus unstable manifold of the source is
%     the image of a horizontal disk over the $x, y$ directions.
%
%     \item In $B_\nu^+ := \interval{\frac{\sqrt{\nu}}{2}, \sqrt{2\nu}} \times \interval{-\nu^{\frac{3}{2}}, \nu^{\frac{3}{2}}}^2$
%     and in $B_\nu^- := \interval{-\sqrt{2\nu}, -\frac{\sqrt{\nu}}{2}, } \times -\interval{-\nu^{\frac{3}{2}}, \nu^{\frac{3}{2}}}^2$ the cone conditions for the flow are satisfied with quadratic form $\mathcal{Q}_2$
%     given by the matrix
%     $$
%     \begin{pmatrix}
%       &-1 &0 &0 \\
%       &0 &1 &0 \\
%       &0 &0 &-1
%     \end{pmatrix}.
%     $$
%     Those sets are also isolating cuboid, thus in those sets there exist unique
%     fixed points $z_\pm$ whose stabale manifolds are vertical disks over the $x, z$ directions.
% \end{itemize}
%

\section{The method of self-consistent bounds}
\label{sec:method}

In this section we recall the method of self-consistent bounds developed in \cite{ZM, ZNS, ZA}.
Let $H $ be a real Hilbert space. We study the following
equation
\begin{equation}
\label{eq_pde} \frac{du}{dt} = F(u)
\end{equation}
where the domain of $F$ is  dense in $H$. By a solution of
(\ref{eq_pde}) we understand a function $u:[0,t_{max}) \to
\dom(F)$ such that $u$ is differentiable and (\ref{eq_pde}) is
satisfied for all $t \in [0,t_{max})$.
The scalar product in $H$  will be denoted by $(u | v)$ and norm by $\norm{u}$,
for $u, v \in H$.
Throughout this section we assume that there \textcolor{black}{exists} a set $I \subset
\mathbb{Z}^d$  \textcolor{black}{(in our work we only need $d=1$, but in this section we present the theory for the more general case, following the cited works \cite{ZM, ZNS, ZA})}
and a sequence of subspaces $H_k \subset H$ for $k
\in I$ such that $\dim H_k=d_k < \infty$ and $H_k$ and $H_{k'}$
are mutually orthogonal for $k \neq k'$. Let $A_k: H \to H_k$ be
the orthogonal projection onto $H_k$. We assume that for each $u
\in H$ holds
\begin{equation}
   u=\sum_{k \in I} u_k, \label{eq_H-decmp}
\end{equation}
where $u_k := A_k u$.
Analogously, for a function $B$ with its range in $H$
we set $B_k(u) := A_k B(u)$ (in particular, $F_k(u) = A_k F(u)$). Equation (\ref{eq_H-decmp}) implies that
$H = \overline{\bigoplus_{k \in I} H_k}$.

For $k \in \mathbb{Z}^d$ we define
\begin{displaymath}
  \abs{k}=\sqrt{\sum_{i=1}^d k_i^2}
\end{displaymath}

 For $n > 0$ we set
\begin{align*}
  X_n &= \bigoplus_{\abs{k} \leq n, k \in I } H_k \\
  Y_n &= X_n^\bot.
\end{align*}
By $P_n:H \to X_n$ and $P_n^\perp :H\to Y_n$ we will denote the
orthogonal projections onto $X_n$ and onto $Y_n$, respectively.

% \cite[Def. 2.1, 2.11]{ZM}
\begin{definition}
We say that $F:H \supset \dom (F) \to H$
is admissible  if the following conditions are satisfied for any
$i, j \in I$, such that $\dim X_i, \dim X_j >0 $
\begin{itemize}
\item $ X_i \subset \dom(F)$,
\item $P_i F : X_i \to X_i$ is a $C^1$ function,
\item  $\frac{\partial F_i}{\partial x_j}$ exists and is continuous.
\end{itemize}
\end{definition}

\begin{definition}
Assume $F$ is admissible. For a given number $n>0$ the ordinary
differential equation
\begin{equation}
  x' = P_n F(x), \textnormal{ where } x \in X_n,  \label{eq_galproj}
\end{equation}
will be called \emph{the $n$-th Galerkin projection} of
(\ref{eq_pde}).

By $\varphi^n(t,x)$ we denote the local flow on $X_n$ induced by
(\ref{eq_galproj}).
\end{definition}

\begin{definition}
\label{defn:selfconsistent} Assume $F$ is an admissible function.
Let $m,M \in \mathbb{N}$ with $m \leq M$. Consider a compact set
$W\subset X_m$ and a sequence of compact sets $B_k \subset H_k$ for $\abs{k}
> m$, $k \in I$.
We say that $W \oplus \Pi_{k \in I, \abs{k} > m} B_k$ forms \emph{self-consistent  bounds for $F$}
if

\begin{itemize}
\item[(C1)] For $k \in I$ such that $\abs{k} > M$ we have $0 \in B_k$ .
\item[(C2)] Let $\hat{a}_k := \max_{a \in B_k} \norm{a}$ for $\abs{k} > m$, $k \in I$. Then
$\sum_{\abs{k} > m, k \in I} \hat{a}_k^2 < \infty$. In particular
\begin{equation*}
   W \oplus \Pi_{\abs{k} > m} B_k \subset H
\end{equation*}
and for every $u \in  W \oplus \Pi_{k \in I,\abs{k} > m} B_k $ holds,
$\norm{P_n^\perp  u} \leq \sum_{\abs{k}>n, k \in I}\hat{a}_k^2$.
\item[(C3)] The function $u\mapsto F(u)$ is continuous on
$W\oplus\prod_{k \in I, \abs{k}>m}B_k \subset H$.
Moreover, if for $k \in I$ we define $F_k=\max_{u \in W \oplus \prod_{k \in I, \abs{k} >m} B_k}
\norm{F_k(u)} $, then $\sum F_k^2 < \infty$.
\item[(C4a)] For $|k|>m$, $k \in I$  $B_k$ is given by
\begin{align}
  B_k&=[a_{k}^-,a_k^+], \qquad a_s^- < a_s^+, \:
  s=1,\dots,d_1 \label{eq:C4ain}
\end{align}
 Let $u \in W \oplus \Pi_{|k| > m}B_k$.
 Then for $|k| > m$
  \begin{align}
     u_{k}=a_{k}^-  & \Rightarrow   F_{k}(u) >   0,  \label{eq:C4ain-} \\
     u_{k}=a_{k}^+  & \Rightarrow   F_{k}(u) <   0.  \label{eq:C4ain+}
 \end{align}
\end{itemize}
\end{definition}
The equations
(\ref{eq:C4ain-}, \ref{eq:C4ain+}) are called  \emph{isolation equations}.
In our case those equations will be verified alongside with the proof that
the set is forward-invariant or an isolating cuboid.
% Typically the definition of the self-consistent bounds includes also isolation
% conditions (hence justifying the \emph{self-consistent} in name), but in our case
% it would be redundant, as in the sequel we will always work with either forward-invariant
% sets or with isolating cuboids.
\par
If the choice of $F$ is clear from the context, then we often drop $F$
 and we will speak simply  about \emph{self-consistent bounds}.
\par
 Given self-consistent  bounds $W$ and $\{B_k\}_{k \in I, \abs{k}>m}$,
 by $T$ (the tail) we will denote
\[
T:=\prod_{\abs{k} > m} B_k \subset Y_m.
\]

\subsection{Convergence of the Galerkin projections}
\label{sec:cov-proj}

% \subsection{Convergence of the Galerkin projections and the Lipschitz constants}

In this section we get the existence of the (local) semiflows for
the self-consistent bounds which are forward-invariant or which
are isolating cuboids. We also get that if every Galerkin projection
of the self-consistent bounds is forward invariant or an isolating
cuboid, then the same is true for the set itself.
\par
We now list two theorems about the convergence of the Galerkin
projections, a global version in a forward-invariant set and a local one
in an isolating cuboid. We cite here slightly modified result \cite[Theorem 13]{ZNS}.
Proof of this theorem in \cite{ZNS} is based on logarithmic norms (defined and discussed
in Appendix \ref{app_lognorms}), but it is done
only for the $l_2$ norm. In the computer-assisted proof of the
heteroclinic connections we use the max norm, so we rephrase this
theorem to work for every logarithmic norm. We discuss it in
a bit more detailed way in Appendix \ref{app_lognorms}.
\par
The following condition tells us that the logarithmic norms
of the projections are uniformly bounded. It is the most
important condition in proving the convergence of the Galerkin
projections.
\begin{definition}
  \label{def_condition_D}
  Consider self-consistent bounds $V$ for $F$. We define the
  following condition
  \begin{itemize}
    \item[(D)]  there exists $l \in \mathbb{R}$ such that for any
    any $n \in \mathbb{N}$ we have
    \begin{align*}
      \mu\expr{D P_n F}  \textcolor{black}{\leq}  l.
    \end{align*}
    \end{itemize}
\end{definition}

By Theorem \ref{thm_lognorm_bounds}, we easily get the following way of
verifying the condition (D).
\begin{lemma}
Let $R \subset H$ be convex self-consistent bounds for $F$.
If there exist $l \in \mathbb{R}$ such that and for all $i \in I$and $x \in R$ we have
\begin{itemize}
  \item for norm $\norm{\cdot} = \norm{\cdot}_2$ is the $l_2$ norm
  \begin{align*}
    \frac{\partial F_{i}}{\partial x_{i}}\expr{x}  +
    \sum_{k,\: k \neq i}
        \frac{1}{2} \expr{
            \abs{\frac{\partial F_{i}}{\partial x_{k}}(x)} +
            \abs{\frac{\partial F_{k}}{\partial x_{i}}(x)}}  \leq l.
  \end{align*}
  then the condition (D) holds with the $l_2$ norm.
  \item for the maximum norm $\norm{\cdot} = \norm{\cdot}_\infty$
  $$
   \frac{\partial F_{i}}{\partial x_{i}}(x) +
  \sum_{k,\: k \neq i} \abs{\frac{\partial F_{i}}{\partial
  x_{k}}(x)}  \leq l,
  $$
  then the condition (D) holds with the maximum norm.
\end{itemize}
\end{lemma}

We now proceed to the proof of the convergence on the
forward-invariant self-consistent bounds. We start by defining
the forward-invariance.

\begin{definition}
  We say than a closed set $R \subset X$ is \emph{forward-invariant} with respect
  to the semiflow $\varphi$
  if for all $x \in R$ and $t \in \mathbb{R}_+$ we have $\varphi(t, x) \in R$.
\end{definition}

\begin{theorem}
\label{thm:conv-on-Trap}
  Assume that $\norm{\cdot}$ is the $l_2$ or the maximum norm on $H$
  and let $\mu$ be the associated logarithmic norm.
  Let $R \subset H$ be a convex self-consistent bounds for $F$ such that
  the condition (D) holds on $R$.
  \par
  Assume that $P_n(R)$ is forward-invariant for the $n$-dimensional
  Galerkin projection of (\ref{eq_pde}) for all $n > M$. Then
  \begin{itemize}
    \item[1.]{\bf Uniform convergence and existence}. For a fixed $x_0 \in R$,
    let  $x_n:[0,\infty] \to P_n(R)$ be a solution of $x'=P_n(F(x))$, $x(0)=P_nx_0$.
    Then $x_n$ converges uniformly on compact  intervals to a function
     $x^*:[0,\infty]\to R$, which
    is a solution of (\ref{eq_pde}) and $x^*(0)=x_0$.  The
    convergence of $x_n$ on compact time intervals is uniform with
    respect to $x_0 \in R$.
    \item[2.] {\bf Uniqueness within $R$}. There exists only one solution of the initial
     value problem (\ref{eq_pde}), $x(0)=x_0$ for any $x_0 \in R$, such that $x(t) \in R$
     for $t>0$.
    \item[3.]  {\bf Lipschitz constant}. Let $x:[0,\infty] \to R$ and
      $y:[0,\infty] \to R$ be solutions of (\ref{eq_pde}), then
      \begin{align}\label{eq_lip_constant}
         \norm{y(t) - x(t)} \leq e^{lt}\norm{x(0) - y(0)}
     \end{align}
    \item[4.] {\bf Semiflow.} The map $\varphi:\mathbb{R}_+ \times R \to R$, where
       $\varphi(\cdot,x_0)$ is a unique solution of equation (\ref{eq_pde}), such that
       $\varphi(0,x_0)=x_0$ defines a semiflow on $R$.
  \end{itemize}
\end{theorem}

The following theorem is our main tool in proving the existence of the
attracting  fixed points.
\begin{theorem} \cite[Theorem 3.8]{ZA}
\label{thm:attractfp} Let same assumptions on $R$ and $F$ as in Theorem \ref{thm:conv-on-Trap} hold.
Assume that $l <0$.
Then there exists a unique fixed point $x^* \in R$ for (\ref{eq_pde}).
Moreover, for every $y \in R$
\begin{align*}
  \norm{\varphi(t,y) - x^*} &\leq e^{lt}\norm{y-x^*}, \textnormal{for } t \geq
  0,\\
  \lim_{t\to \infty} \varphi(t,y) &= x^*.
\end{align*}
\end{theorem}

\begin{remark}
    \label{rem_stable_hyperbolicity}
    Observe that if we equip $H$ with $l_2$ norm,
    then the attracting fixed point $x^*$ in Theorem \ref{thm:attractfp} is
    hyperbolic in the sense of Definition \ref{def_hyperbolicity}. To prove this
    it is enough to apply Theorem \ref{thm_hyperbolicity_characterization} with
    the $l_2$ norm and $\mathcal{Q} = -\id$.
\end{remark}

The fact that $R$ is forward-invariant was essential in proof of Theorem \ref{thm:conv-on-Trap}
to ensure that  the solutions stay in $R$ for all times $t > 0$. Nevertheless, if we
only could assume that $\varphi(t, x)$ stay in $R$ for times $0 \leq t \leq t(x)$,
the inequality (\ref{eq_lip_constant}) would still be valid for those times.
Thus, a slight modification of the proof of \cite[Theorem 13]{ZNS}
yields the proof of the following result.

% Let $V$ be convex self-consistent bounds for  $\dot{x} = F(x)$. Assume that
% for some $m \leq M$ and $\varepsilon > 0$ there exists a homeomorphism $c_u:X_m \to \mathbb{R}^m$ such that
% for the ball $B(0, 1) \subset X_{m}$ we have
% that $P_n\expr{V} = c_u^{-1} (\overline{B}(0,1)) \oplus \bigoplus_{i = m' + 1} ^ n B_k$ is an $\varepsilon$-isolating cuboid for $\varphi_n$ for all $n > M$.

% Assume moreover that the following condition holds.
% \begin{itemize}
% \item[(D)] There exists $l \in \mathbb{R}$, such that for any $i$
% and $x \in V^\varepsilon := c_u^{-1}(\overline{B}(0,1+\varepsilon)) \oplus \bigoplus_{i = m + 1} ^ \infty B_k$
% we have
% \begin{align*}
% \frac{\partial F_{i}}{\partial x_{i}}\expr{x}  +
% \sum_{k,\: k \neq i}
%     \frac{1}{2} \expr{
%         \abs{\frac{\partial F_{i}}{\partial x_{k}}(x)} +
%         \abs{\frac{\partial F_{k}}{\partial x_{i}}(x)}}  \leq l.
% \end{align*}
% \end{itemize}

\begin{theorem}
\label{thm:local-conv}
Assume that $\norm{\cdot}$ is the $l_2$ or the maximum norm on $H$
and let $\mu$ be the associated logarithmic norm.
Let $m \leq M$ be natural numbers. Let $c_u:X_m \to \mathbb{R}^m$ be a homeomorphism
and consider the self-consistent bounds for $F$ of the form
$V := c_u^{-1}(\overline{B}(0,1)) \oplus \bigoplus_{k = m + 1} ^ \infty B_k$.
Assume that there exists $\varepsilon > 0$ such that for each
$n > M$ the set $P_n\expr{V} = c_u^{-1} (\overline{B}(0,1)) \oplus \bigoplus_{i = m + 1} ^ n B_k$
is an $\varepsilon$-isolating cuboid for $\varphi_n$.
Moreover, assume that the condition (D) holds for
$V^\varepsilon := c_u^{-1}(\overline{B}(0,1+\varepsilon)) \oplus \bigoplus_{i = m + 1} ^ \infty B_k$.
\par
Then
  \begin{itemize}
    \item[1.]{\bf Local uniform convergence and existence} For a fixed $x_0 \in V^\varepsilon$,
    let  $x_n:[0,t_{max}(n,x_0)] \to P_n(V)^\varepsilon$ be a solution of $x'=P_n(F(x))$, $x(0)=P_nx_0$
    defined on the maximum interval of existence.
    Then  $t_{\max}(n, x_0) \to t_{\max}(x_0)$ and  $x_n$ converges uniformly on
     compact intervals contained in $\interval{0,t_{\max}(x_0)}$ to a function
     $x^*:\interval{0,t_{max}(x_0)} \to V^\varepsilon$, which
    is a solution of (\ref{eq_pde}) and for which $x^*(0) = x_0$.  The
    convergence of $x_n$ on compact time intervals is uniform with
    respect to $x_0 \in V^\varepsilon$.
    \item[2.] {\bf Uniqueness within $V^\varepsilon$}. For any $x_0 \in V^\varepsilon$ there exists only
     one solution of the initial value problem (\ref{eq_pde}), $x(0)=x_0$ such
     that $x(t) \in V^\varepsilon$ for $0 < t < t_{\max}(x_0)$.
    \item[3.]  {\bf Lipschitz constant}. Let $x:[0,t_1] \to V^\varepsilon$ and
      $y:[0,t_1] \to V^\varepsilon$ be solutions of (\ref{eq_pde}), then
      \begin{align*}
         \norm{y(t) - x(t)} \leq e^{lt}\norm{x(0) - y(0)}.
      \end{align*}
    \item[4.] {\bf Local semiflow.} The partial map
       $\varphi:\mathbb{R}_+ \times H \rightharpoonup H$, where
       $\varphi(\cdot, x_0)$ is the unique solution of equation (\ref{eq_pde}) such that
       $\varphi(0, x_0) = x_0$, defines a local semiflow on $V$.
    \item[5.] {\bf Isolation}. The set $V$ is an $\varepsilon$-isolating cuboid for $\varphi$.
  \end{itemize}
\end{theorem}

For set $V$ as in theorem above we call directions $\set{1,\ldots,m}$ the \emph{exit}
and the remaining directions the \emph{entry} directions.

\subsection{Verification of the cone conditions} \label{sec_cc_ver}

We will need a slight modification of Lemma 3.1 proved in \cite{ZBif}.

\begin{lemma} \label{lemma_averaging}
    Let $V \subset H$ be self-consistent bounds for an admissible function $F$.
    Moreover, assume that the following condition holds.
    \begin{itemize}
        \item[(F)] if we set $d_{ij} := \max_{x \in V} \abs{\frac{\partial F_i}{\partial x_j}}$, then for
        each $i \in \mathbb{N}$ the sum
            $$
                \sum_{j = 1}^\infty d_{ij} \sup_{x, y \in V} \abs{x_j - y_j}
            $$
            converges.
    \end{itemize}
    Then for all $x, y \in V$ we have
    $$
    F(x) - F(y) = \overline{DF} \expr{x - y},
    $$
    where $$\overline{DF}_{ij} := \int_{0}^{1}\frac{\partial F_i}{\partial x_j}\expr{tx + \expr{1-t}y} \,\mathrm{d}t.$$
\end{lemma}
\begin{proof}
By Lemma 3.1 from \cite{ZBif}, we have
$$
    F_{i}(x) - F_i(y) =
        \sum_{j = 1}^\infty
            \int_{0}^{1} \frac{\partial F_i}{\partial x_j}\expr{tx + \expr{1-t}y} \,\mathrm{d}t \expr{x_j - y_j}.
$$
Thus, it is enough to prove that $P_n F(x) - P_n F(y) \to F(x) - F(y)$, which
immediately follows by (C3).
\end{proof}

\begin{theorem} \label{thm_unstable_manifold_selfcb}
    Let $V$ be self-consistent bounds for $F$
    and assume that $F$ defines a continuous local semiflow $\varphi$ on $V^\varepsilon$.
    Moreover, assume that $V$ is an $\varepsilon$-isolating cuboid with unstable dimension $m$.
    \par
 Let $Q$ be an infinite diagonal matrix such
that $Q_{ii}=1$ for $i = 1, \dots, m$ and
$Q_{ii}=-1$ for $i > m$.
We denote $[DF(V)]_{ij} := \set{t \in \mathbb{R} \mid \exists z \in V \; t = \frac{\partial F_i}{\partial x_j}(z)}$.
Assume that there exists $\lambda > 0$ such that for all $\delta > 0$ small enough and
for each symmetric matrix
$$A \in \left( [DF(V)]^T \expr{Q \pm \delta \Id} + \expr{Q \pm \delta \Id} [DF(V)] \right)$$
we have
\begin{align} \label{eq_wAw}
    \expr{w | A w} \geq \lambda \norm{w}^2,
\end{align}
for all $w \in H$ such that $w = z_1 - z_2$ for some $z_1, z_2 \in V$.
\par
Finally, let $\mathcal{Q}$ be a quadratic form defined by the matrix $Q$ on the vector space
generated by the set $\set{w \in H \mid \exists z_1, z_2 \in V \; w = z_1 - z_2}.$
Then there exists a unique fixed point  $z_0 \in V$ which is hyperbolic and for which
 $W^u_{V}(z_0)$ ($W^s_{V}(z_0)$) is the image of a horizontal (vertical)
 disk in $V$ with respect to $\mathcal{Q}$.
\end{theorem}
\begin{proof}
    We will show that the strong cone conditions are satisfied on $V$.
    % with the quadratic form $\mathcal{Q}$ defined by the matrix $Q$ on the vector space
    % generated by the set $$\set{w \in H \mid \exists z_1, z_2 \in V \; w = z_1 - z_2}.$$
    Consider $z_1, z_2 \in V$.
   %  By the condition \ref{eq_cone_cond_convergence1} we have
   %  \begin{align*}
   %    \expr{w | A w} &= \expr{P_n w + P_n^\perp w | A \expr{P_n w + P_n^\perp  w}} \\
   %    &= \expr{P_n w| A P_n w} + 2\expr{ P_n^\perp w | A P_n w} + \expr{ P_n^\perp w | A P_n^\perp  w} \\
   %    &\geq \lambda \norm{P_nw}^2 +
   %     2( P_n^\perp w | A P_n w) + ( P_n^\perp w | A P_n^\perp  w).
   % \end{align*}
   % Due to the condition (\ref{eq_cone_cond_convergence2}), $2( P_n^\perp w | A P_n w) + ( P_n^\perp w | A P_n^\perp  w) \xrightarrow[n \to \infty]{} 0$, thus we get
   % \begin{align}\label{eq_gershgorin_like_bound}
   %     \expr{w | A w} \geq \lambda \norm{w}^2.
   % \end{align}
   % \par
   Due to Lemma \ref{lemma_averaging}, there exists
   a matrix $\overline{DF} \in \expr{[ DF(V)]}$ such that
   \begin{align*}
       \frac{d}{dt}\expr{\varphi(t, z_1) - \varphi(t, z_2)} = F\expr{\varphi{t, z_1}} - F\expr{\varphi{t, z_2}} = \overline{DF}\expr{\varphi(t, z_1) - \varphi(t, z_2)},
   \end{align*}
   thus for small enough $\delta > 0$ we have
   \begin{align*}
       \frac{d}{dt} &\expr{\varphi(t, z_1) - \varphi(t, z_2) | \expr{Q \pm \delta \Id}  \expr{\varphi(t, z_1) - \varphi(t, z_2)}} \\
       &=  \expr{\frac{d}{dt}\expr{\varphi(t, z_1) - \varphi(t, z_2)} | \expr{Q \pm \delta \Id}  \expr{\varphi(t, z_1) - \varphi(t, z_2)}} \\
          &+ \expr{\varphi(t, z_1) - \varphi(t, z_2) | \expr{Q \pm \delta \Id}  \frac{d}{dt} \expr{\varphi(t, z_1) - \varphi(t, z_2)}} \\
       &=  \expr{\overline{DF}\expr{\varphi(t, z_1) - \varphi(t, z_2)} | \expr{Q \pm \delta \Id}  \expr{\varphi(t, z_1) - \varphi(t, z_2)}} \\
          &+ \expr{\varphi(t, z_1) - \varphi(t, z_2) | \expr{Q \pm \delta \Id}  \overline{DF} \expr{\varphi(t, z_1) - \varphi(t, z_2)}} \\
       &=  \expr{\varphi(t, z_1) - \varphi(t, z_2)}^T \expr{\overline{DF}^T \expr{Q \pm \delta \Id}
       + \expr{Q \pm \delta \Id} \overline{DF} }\expr{\varphi(t, z_1) - \varphi(t, z_2)} \\
       &\geq \lambda \norm{\varphi(t, z_1) - \varphi(t, z_2)}^2 > 0,
   \end{align*}
   by (\ref{eq_wAw}). For $t = 0$ we get
   \begin{align*}
       \frac{d}{dt} &\expr{\varphi(t, z_1) - \varphi(t, z_2) | \expr{Q \pm \delta \Id}  \expr{\varphi(t, z_1) - \varphi(t, z_2)}} |_{t = 0} > \\
       &> \lambda \norm{z_1 - z_2}^2 =  \lambda \expr{\norm{P_n\expr{z_1 - z_2}}^2 + \norm{P_n^\bot\expr{z_1 - z_2}}^2}.
   \end{align*}
   The inequality above easily implies (\ref{eqn_strong_cc_1}, \ref{eqn_strong_cc_2}). Thus the claim
   follows by Theorems \ref{thm_stable_manifold}, \ref{thm_unstable_manifold} and
   \ref{thm_hyperbolicity_characterization}.
\end{proof}

Now we give sufficient conditions for (\ref{eq_wAw}).

\begin{lemma} \label{lemma_cc_verification}
    Let $V$ be self-consistent bounds for $F$ and let $Q$ be as in Theorem \ref{thm_unstable_manifold_selfcb}.
    \par
 Assume that for each symmetric matrix $A \in \left( [DF(V)]^T Q + Q [DF(V)] \right)$
 we have

 \begin{align} \label{eq_cone_cond_convergence1}
   A P_n w \xrightarrow[n \to \infty]{} Aw
 \end{align}
 for all $m \in \mathbb{N}$ and $w \neq 0$ such that $w=z_1 - z_2$, where $z_1, z_2 \in V$.
 \par
 If for some $\lambda > 0$ the following inequality
\begin{equation}
  2 \inf_{x \in V}\left|\frac{\partial F_i}{\partial x_i}(x) \right| - \sum_{j, j \neq i}
     \sup_{x \in V} \left| Q_{jj}\frac{\partial F_j}{\partial x_i}(x) +
     Q_{ii} \frac{\partial F_i}{\partial x_j}(x)  \right| \geq
     \lambda  \label{eq_cc-diag-dom}
\end{equation}
holds for all $i \in \mathbb{N}$, then
also condition %s
% (D) and
(\ref{eq_wAw}) holds.
\end{lemma}
\begin{proof}
    Let $\delta > 0$ and consider the infinite symmetric matrix
    $$A \in \expr{[DF(V)]^T \expr{Q \pm \delta \Id} + \expr{Q \pm \delta \Id} [DF(V)]^T}.$$
    Let $n > 0$ and let $A_n: P_n V \to P_n V$ be such that for $w \in V$ we have
    $A P_n w = A_n P_n w$ (it is easy to check that it is well-defined). Observe
    that for $i \leq n$ $$2\inf_{x \in V} \abs{\frac{\partial F_i}{\partial x_i}(x)} - 2\delta$$
    is a lower bound for the Gershgorin center of the $i$-th row of $A_n$ and $$\sup_{x \in V} \left| Q_{jj}\frac{\partial F_j}{\partial x_i}(x) +
    Q_{ii} \frac{\partial F_i}{\partial x_j}(x)  \right|$$ is the upper bound for the
    Gershgorin radius of the $i$-th row of $A_n$.
    Thus if $\delta > 0$ is small enough, then by (\ref{eq_cc-diag-dom}) and the Gershgorin theorem we have
    \begin{align} \label{eq_cone_cond_convergence2}
        \expr{P_n w | A \expr{P_n w}} \geq \lambda_\delta \norm{P_n w}^2,
    \end{align}
    where $\lambda_\delta = \lambda - 2\delta$.
    It follows that
     \begin{align*}
       \expr{w | A w} &= \expr{P_n w + P_n^\perp w | A \expr{P_n w + P_n^\perp  w}} \\
       &= \expr{P_n w| A P_n w} + 2\expr{ P_n^\perp w | A P_n w} + \expr{ P_n^\perp w | A P_n^\perp  w} \\
       &\geq \lambda_\delta \norm{P_n w}^2 +
        2( P_n^\perp w | A P_n w) + ( P_n^\perp w | A P_n^\perp  w).
    \end{align*}
    Due to the condition (\ref{eq_cone_cond_convergence1}),
    $2( P_n^\perp w | A P_n w) + ( P_n^\perp w | A P_n^\perp  w) \xrightarrow[n \to \infty]{} 0$
    and thus (\ref{eq_wAw}) follows.
    % \par
    % Condition (D) holds obviously by (\ref{eq_cc-diag-dom}).
\end{proof}

Above we dealt with the case when all directions are either strongly unstable or strongly stable.
In that case we could get results about the invariant manifolds.
However, in the proof of the bifurcation when unstable directions
are present we need a result about the cone invariance when
not every direction is strongly stable or unstable (see proof of the Lemma \ref{lemma_model_cone_field}
for a motivation).
This is why  {\textcolor{red} we deal with } {\textcolor{black} we consider} the setting when we have some central directions, i.e. directions
whose stability depends on the part of the set we consider (in our case
we will have one central direction -- the bifurcation direction).

\begin{theorem} \label{thm_disk_image_selfcb}
    Let $V$ be self-consistent bounds for $F$
    and assume that $F$ defines a continuous local semiflow $\varphi$ on $V$.
    Moreover, assume that $V$ is an $\varepsilon$-isolating cuboid with the unstable directions
    $x_1, \dots, x_m$.
    For $z \in V$ we write $z = (x, y)$, where $x$ are unstable directions.
    % Let $m \leq M$ be natural numbers. Let $c_u:X_m \to \mathbb{R}^m$ be a homeomorphism
    % and consider the self-consistent bounds of the form
    % $V := c_u^{-1}(\overline{B}(0,1)) \oplus \bigoplus_{k = m + 1} ^ \infty B_k$.
    % Assume that there exists $\varepsilon > 0$ such that for each
    % $n > M$ the set $P_n\expr{V} = c_u^{-1} (\overline{B}(0,1)) \oplus \bigoplus_{i = m + 1} ^ n B_k$
    % is an $\varepsilon$-isolating cuboid for $\varphi_n$.
    % Assume moreover that conditions (D) and (F) are satisfied.
    \par
 Let $Q$ be an infinite diagonal matrix such that $Q_{ii}=1$ for $i = 1, \dots, m$ and
$Q_{ii}=-1$ for $i = m + 1, \dots$
\par
Assume that
for some $\lambda > 0$ we have for all
$A \in \left( [DF(V)]^T Q + Q [DF(V)] \right)$

\begin{align}
    \expr{w_x, A_{xx} w_x} &> \lambda \norm{w_x}^2, \label{eq_waw21} \\
    \expr{w_y, A_{yy} w_y} &> -\frac{\lambda}{4} \norm{w_y}^2, \label{eq_waw22} \\
    \norm{A_{xy}} + \norm{A_{yx}} &< \frac{\lambda}{4}, \label{eq_waw23}
\end{align}
for all $w \in H$ such that $w = z_1 - z_2$ for some $z_1, z_2 \in V$.
Finally, let $\mathcal{Q}$ be a quadratic form defined by the matrix $Q$ on the vector space
generated by the set $\set{w \in H \mid \exists z_1, z_2 \in V \; w = z_1 - z_2}.$
Then there exists $T > 0$ such that
for every horizontal disk $h$ in $S$ with respect to $\mathcal{Q}$
the set $\varphi\expr{\im h, t}$ contains
the image of another horizontal disk for all $t \in \interval{0, T}$.
\end{theorem}
\begin{proof}
    Reasoning analogously as in the proof of Theorem \ref{thm_unstable_manifold_selfcb},
    by Theorem \ref{disk_image_lemma} it is enough to prove that the positive cones are invariant.
    We show it by proving that $$\frac{d}{dt}\mathcal{Q}\expr{\varphi(t, z_1) - \varphi(t, z_2)}\mid_{t = 0} > 0$$
    whenever $\mathcal{Q}\expr{z_1 - z_2} > 0$.
    \par
    Let
    $z_1 = (x_1, y_1), z_2 = (x_2, y_2) \in V$ be
    such that $\expr{\expr{z_1 - z_2}, Q \expr{z_1 - z_2}} > 0$ and denote $w = z_1 - z_2$.
    For such $w$ we have $\norm{w_x}^2 > \norm{w_y}^2$.
    For $A \in \left( [DF(V)]^T Q + Q [DF(V)] \right)$, by similar computations as in the proof of
    Theorem \ref{thm_unstable_manifold_selfcb}  we get for $t = 0$
    \begin{align*}
        \frac{d}{dt} &\expr{\varphi(t, z_1) - \varphi(t, z_2) | Q  \expr{\varphi(t, z_1) - \varphi(t, z_2)}}  \geq \\
         &\geq \lambda \norm{w_x}^2 - \frac{\lambda}{4} \norm{w_y}^2 - \expr{w_x, A_{xy} w_y} - \expr{w_y, A_{yx} w_x} \\
         &\geq \expr{\lambda - \frac{\lambda}{4} - \frac{\lambda}{4}} \norm{w_x}^2 > 0.
    \end{align*}
\end{proof}

Analogously to the Lemma \ref{lemma_cc_verification} we can prove that
when we the first $m$ directions are central, next $n$ are strongly unstable
and remaining are stable, the conditions (\ref{eq_waw21}) -- (\ref{eq_waw23}) hold.

\begin{lemma} \label{lemma_cc_verification2}
    % Let $m \leq M$ be natural numbers. Let $c_u:X_m \to \mathbb{R}^m$ be a homeomorphism
    % and consider the self-consistent bounds of the form
    % $V := c_u^{-1}(\overline{B}(0,1)) \oplus \bigoplus_{k = m + 1} ^ \infty B_k$.
    % Assume that there exists $\varepsilon > 0$ such that for each
    % $n > M$ the set $P_n\expr{V} = c_u^{-1} (\overline{B}(0,1)) \oplus \bigoplus_{i = m + 1} ^ n B_k$
    % is an $\varepsilon$-isolating cuboid for $\varphi_n$.
    Let $V$ be self-consistent bounds for $F$ and assume that for some $m, n \geq 0$ the directions
    $i = m + 1, \dots, m + n$ are unstable and the directions $i \geq m + n + 1$ are stable.
    Let $Q$ be a diagonal
    matrix such that $Q_{ii} = -1$ for $i = 1, \ldots, m$ or $i > m + n$ and such that
    $Q_ii = 1$ otherwise.
    % and assume that $F$ defines an associated continuous local semiflow $\varphi$ on $V$.
    \par
 Assume that for each $A \in \left( [DF(V)]^T Q + Q [DF(V)] \right)$
 we have (\ref{eq_cone_cond_convergence1}).
 If for some $\lambda > 0$ we have
 \begin{equation} \label{eq_cc-diag-dom-3}
    -2 \inf_{x \in V}\left|\frac{\partial F_1}{\partial x_i}(x) \right| - \sum_{j, j \neq 1}
       \sup_{x \in V} \left| Q_{jj}\frac{\partial F_j}{\partial x_1}(x) +
       Q_{ii} \frac{\partial F_1}{\partial x_j}(x)  \right| \geq
       -\frac{\lambda}{4},
  \end{equation}
  for $i \leq m$ and
\begin{equation} \label{eq_cc-diag-dom-2}
  2 \inf_{x \in V}\left|\frac{\partial F_i}{\partial x_i}(x) \right| - \sum_{j, j \neq i}
     \sup_{x \in V} \left| Q_{jj}\frac{\partial F_j}{\partial x_i}(x) +
     Q_{ii} \frac{\partial F_i}{\partial x_j}(x)  \right| \geq
     \lambda,
\end{equation}
for all $i > m$,
then the conditions (\ref{eq_waw21}) -- (\ref{eq_waw23}) hold.
\end{lemma}

\section{Analytical proof of  \textcolor{black}{a}  pitch\-fork bi\-fur\-ca\-tion in the Ku\-ra\-moto--Si\-va\-shin\-sky equation} \label{sec_proof_pitchfork_ks}

The goal of this section is to prove that the KS equation has the pitchfork bifurcation at $\mu=1$. This is formulated as Theorem~\ref{thm_bif} at the end of this section.

As in \cite{ZM} we write (\ref{eq_ks_pde}) in Fourier basis

\begin{align} \label{eqn_ks_fourier_with_mu}
\begin{split}
  \dot{a_k}&=k^2(1- \textcolor{black}{\mu}  k^2)a_k - k \sum_{i=1}^{k-1} a_i a_{k-i} + 2 k \sum_{i=1}^{\infty} a_i a_{k+i}=:F_k(\mu,a), \quad k=1,2,\dots
\end{split}
\end{align}

We also denote
\begin{align*}
    \lambda_k(\mu) &= k^2(1 - \mu k^2). \\
\end{align*}
The bifurcation occurs with respect to the parameter $\lambda_1(\mu)$ (when $\mu = 1$). We could
of course write other $\lambda_k$ in terms of $\lambda_1$, thus making it a new parameter,
but it is more convenient to leave the equation in form where $\mu$ is the parameter.

\textcolor{black}{In the remainder of this section, for functions whose only argument is $\mu$ and the value of $\mu$ can be clearly deduced from the context, we do not write this argument explicitly.}
\par
First, we transform the KS equation to the normal form. When we are done with that,
we proceed to the realization of \textbf{Steps I--III} which were described in Section \ref{subsection_pitchfork_model}.

\subsection{Normal forms} \label{sec_normal_forms}

Let us discuss the transformation of the KS equation (\ref{eqn_ks_fourier_with_mu}) to the normal form.
The idea behind the normal forms is to introduce the transformation of variables to simplify the equation.
{\textcolor{red} Such approach was the basis of, for example, KAM theory.}
\par
Our goal is to transform (\ref{eqn_ks_fourier_with_mu}) to
\begin{align*}
    b_1' &= \lambda_1(\mu) b_1 + c(\mu)b_1^3 + \textnormal{higher order terms}, \\
    b_k' &= \lambda_k(\mu) b_k + \textnormal{higher order terms}, k = 2, 3, \dots,
\end{align*}
where $c(\mu) < 0$. To see what higher order terms are acceptable for our
purposes see the discussion at the end of Section \ref{subsection_pitchfork_model}.
\par
In general,  transformations we use are given by the inverse, i.e.
$$a_k = b_k + p(a_1, \dots, a_{k - 1}, b_k, a_{k + 1}, \dots),$$ where
$a_i$s are the old variables and $b_k$ is the transformed variable.
When the equation is of the form $a_i' = F_i(a_1, a_2, \dots)$, then
writing the equation in the new coordinates for $i \neq k$ is easy
-- it is simply
$a_i' = F_i\expr{a_1, a_2, \dots, a_{k - 1}, b_k + p\expr{a_1, \dots, a_{k - 1}, b_k, a_{k + 1}, \dots}, a_{k + 1}}$.
We sketch here how to derive the vector field in the new coordinates for the KS equation.
\par

To discuss  \textcolor{black}{orders of magnitude} of the terms that arise, we need to work on some sets.
Those sets are of the form
\begin{equation}
S_C := \interval{-\sqrt{2} C, \sqrt{2} C} \times \prod_{k = 2}^\infty \interval{\frac{-C^3}{k^s}, \frac{C^3}{k^s}},  \label{eq_SCdef}
\end{equation}
where $C > 0$.

\textcolor{black}
{ 
As in our analysis of the ODE models performed in Section~\ref{sec_pitchfork_models}, we expect $C \approx |b_1(\mu)|$, where $(b_1(\mu),b_2(\mu),\dots)$ is a
fixed point born during bifurcation. Moreover, we want the set $S_C$ to be an isolating cuboid for each Galerkin projection with each face being either entry or exit, so that the linear part dominates the nonlinear parts for $k\geq 2$, while the behaviour in the bifurcation direction will be decided by $b_1'=\lambda_1 b_1 + c(\mu)b_1^3$ which should dominate the remaining terms.
}

We now look at estimates for the sums appearing in $F$ on $S_C$. \textcolor{black}{For $k\geq 2$ we want them to be dominated by the linear terms $\lambda_k a_k$, when evaluated at the boundary, which gives $|\lambda_k a_k| = O(C^3/k^{s-4})$.} It is not hard to see  (see Appendix \ref{sec_app_normal_form})  that on $S_C$  for $k \neq 2$
\textcolor{black}{we have}   $- k \sum_{i=1}^{k-1} a_i a_{k-i} + 2 k \sum_{i=1}^{\infty} a_i a_{k+i} = \frac{O(C^4)}{k^{s - 2}}$ (we will discuss the case $k = 2$ later), hence it is dominated by the linear term for all $k>2$ for $C$ small enough  (we have $k^{s - 2}$ because our bound for the series 'lose' two powers, \textcolor{black}{which is not optimal}).
In the last equation the expression $O(C^4)$ appears because of the terms of the form $a_1 a_i, i > 1$.
We see that if the terms containing $a_1$ are omitted, then those sums are of order $\frac{O(C^6)}{k^{s - 2}}$.
As for $k = 2$, there is a term $-2a_1^2$ which makes the sum to be $O(C^2)$ (which is why
we want to remove it too), but aside
from this term the same remarks hold as for the other $k$'s.

\par
We want to transform
$F$ to the form similar to the discussed models in the Section \ref{sec_pitchfork_models}.
One term we need to remove is $2a_1 a_2$ from $F_1$. 
The transformation which allows us to do this is given by
$a_1 = b_1 + c b_1 a_2$ for some $c \in \mathbb{R}$ to be picked later.
We want to find the equation for $b_1'$. To this end we calculate 
\begin{align*}
    \frac{d}{dt} \expr{b_1 + c b_1 a_2} = b_1' (1 + ca_2) + cb_1 a_2'.
\end{align*}
Given that we have $a_2' = \lambda_2 a_2 - 2b_1^2 + O(C^4)$ ($C^4$ comes from the term $b_1 a_3$)
and that $b_1 O(C^4) = O(C^5)$, we get
\begin{align} \label{eqn_nf_intro1}
    \frac{d}{dt} \expr{b_1 + c b_1 a_2} = b_1' (1 + ca_2) + c \lambda_2 b_1 a_2 - 2cb_1^3 + O(C^5).
\end{align}
On the other hand, we have
\begin{align} \label{eqn_nf_intro2}
    a_1' = \lambda_1 b_1 + c \lambda_1 b_1 a_2 + 2 b_1 a_2 + O(C^5).
\end{align}

When $\abs{ca_2} < 1$, then $\frac{1}{1 + ca_2} = 1 - ca_2 + \sum_{i = 2}^\infty (-ca_2)^i = 1 - ca_2 + O(C^5)$.
Thus comparing (\ref{eqn_nf_intro1}) to (\ref{eqn_nf_intro2}) and rearranging the terms we get

\begin{align*}
    b_1' = \lambda_1 b_1 + 2cb_1^3 + (2 - c \lambda_2)b_1 a_2 + (c\lambda_1 + 2 - c \lambda_2) b_1^2 a_2 + O(C^5).
\end{align*}

Taking $c = \frac{2}{\lambda_2}$, we see that the term $b_1 a_2$ is removed.
\textcolor{black}{Observe that this works only as long as $\lambda_2 \neq 0$, i.e. $\mu \neq \frac{1}{4}$.
This is why we will limit ourselves to $\mu > \frac{1}{2}$, which makes sure we are far from the resonances
for any change of coordinates we will use.}
Another good thing that  happened \textcolor{black}{after the change of coordinates} is that we get the normal form of  \textcolor{black}{a}  pitchfork
bifurcation $\lambda_1 b_1 + \frac{4}{\lambda_2}b_1^3$. Looking closely at the calculations
we see that it was produced by the term $-2b_1^2$ of $F_2$, so the normal
form of  \textcolor{black}{a}  pitchfork bifurcation was in some way 'entangled' in the original equation.
\par
Now we need to remove the term $(c\lambda_1 + 2 - c \lambda_2)b_1^2 a_2$ from the new $F_1$
and $-2b_1^2$ from $F_2$ plus some other terms which would present some additional technical
difficulties (see discussion of terms at the end of Section \ref{subsection_pitchfork_model}).
We also need to calculate bounds for the derivatives of the $F$ in the normal form,
in order to bound the logarithmic norms and verify the cone conditions.
All of this is done with care in Appendix \ref{sec_app_normal_form} and summarized in the following lemma.

\begin{lemma} \label{lemma_snd_normal_form}
  
  There exists $1 > \bar{C} > 0$ such that 
  equation (\ref{eqn_ks_fourier_with_mu}) can be transformed by an \textcolor{black}{analytical} 
  change of variables into
  \begin{align} \label{eqn_ks_normal_form}
      \begin{split}
      \dot{b}_1 &= \lambda_1 b_1 + \frac{4}{\lambda_2}b_1^3 + h_1\expr{\mu, b_1, b_2, a_3, \ldots}, \\
      \dot{b}_k &= \lambda_k b_k + h_k\expr{\mu, b_1, b_2, b_3, \ldots}, k = 2, 3, 4, \ldots \\
  \end{split}
  \end{align}
  where $h_1, h_2, h_3, \ldots$ are such that \textcolor{black}{for every $n \in \mathbb{Z}_{\geq 1}$ 
  there exist constants $\alpha > 0, \alpha_n > 0, \beta > 0, \beta_n > 0, M > 0$ for which the following
  conditions are satisfied: 
  \begin{itemize}
    \item For every $\mu \geq \frac{1}{2}$ and $C \in \interval{0, \bar{C}}$ we have
    \begin{align} \label{eqn_bound_norma_form_nonlinear_parts}
      \begin{split}
      \abs{h_1} &< \alpha_1 C^5, \\
      \abs{h_k} &< \frac{\alpha_k C^4}{k^{s-2}}, \quad k= 2, 3 \ldots
      \end{split}
    \end{align}
    \item For every $\mu \geq \frac{1}{2}$ and $C \in \interval{0, \bar{C}}$ we have
    \begin{align}\label{eqn_bound_norma_form_nonlinear_derivatives}
      \begin{split}
      2\abs{\frac{\partial h_1}{\partial b_1}} + \sum_{i \neq k} \abs{\frac{\partial h_1}{\partial b_i}} + \abs{\frac{\partial h_i}{\partial b_1}} &< \beta_1 C^3, \\
      2\abs{\frac{\partial h_k}{\partial b_k}} + \sum_{i \neq k} \abs{\frac{\partial h_k}{\partial b_i}} + \abs{\frac{\partial h_i}{\partial b_k}} &< \beta_k kC, \quad k=2, 3, \ldots
      \end{split}
    \end{align}
    \item For every $n > M$ we have $\alpha_n = \alpha$.
    \item For every $n > M$ we have $\beta_n = \beta$.
  \end{itemize}
} % endcolor green
\end{lemma}
%\begin{proof}[Proof is provided in Section \ref{sec_app_normal_form}]
%\end{proof}

For $z = (b_1, b_2, b_3, \ldots)$ we denote $\pi_{b_k} z := b_k$.  \textcolor{black}{We also} denote by $\hat{F}_\mu$ and $\varphi_\mu$ respectively the right hand side of (\ref{eqn_ks_normal_form}) and the semiflow associated with it  for a fixed $\mu > 0.$ We may drop $\mu$ when its value is clear from the context.
\par
\textcolor{black}
{
We will need to verify that conditions (F) and (\ref{eq_cone_cond_convergence1}) are satisfied for (\ref{eqn_ks_normal_form}) on $S_{\bar C}$. Hence in order for the assertion of Theorem~\ref{thm_unstable_manifold_selfcb} to hold true we just need to verify condition (\ref{eq_cc-diag-dom}) from Lemma~\ref{lemma_cc_verification}.
It follows from the following lemma, which we obtain
}
by a slight modification of the proof of Theorem 7.1 from \cite{ZBif}.
\begin{lemma} \label{lemma_ks_norm_form_F}
    Fix $\mu > 0$ and let $\bar C \in \expr{0, 1}$ be such that Lemma \ref{lemma_snd_normal_form} holds.
    Denote $$d_{ij} := \max_{x \in S_{\bar C}} \abs{\frac{\partial \hat{F}_{\mu, i}}{\partial x_j}}$$
    Then
    $$
        \sum_{i = 1}^\infty \sum_{j = 1}^\infty d_{ij} \sup_{x, y \in S_{\bar C}} \abs{x_j - y_j} < \infty.
    $$
\end{lemma}

\subsection{Proof of the bifurcation}
\label{sec_proof_of_bifurcation}

We proceed to the realization of \textbf{Steps I--III}
which were described in Section \ref{subsection_pitchfork_model}. Our point of departure is system
(\ref{eqn_ks_normal_form}).

Let $\frac{1}{2} < \mu < 1$. Throughout this section we denote
\begin{align*}
  C(\mu) &:= \sqrt{\frac{-\lambda_1(\mu)\lambda_2(\mu)}{4}},\\
  B_\mu &:= S_{C(\mu)},
\end{align*}
\textcolor{black}{where $S_C$ is defined in (\ref{eq_SCdef}).}

Also, when writing $\alpha_i, \alpha, \beta_i, \beta$ we always mean respectively
$\alpha_i(\mu_+), \alpha(\mu_+), \beta_i(\mu_+), \beta(\mu_+)$, where $\mu_+$ will
be always known from the context.

We will repeatedly use the following remark.
\begin{remark} \label{rem_lambda1}
  For $1 > \mu > \frac{1}{2}$ we have $\frac{C(\mu)^2}{3} < \lambda_1(\mu) < C(\mu)^2$.
  Consequently, $\lambda_1(\mu) = O(C(\mu)^2)$.
\end{remark}

\textcolor{black}{
Now we proceed to \textbf{Step I -- existence of the set isolating the bifurcation}. Since in our setting
existence of the semiflow is connected with the isolation, we prove the mentioned existence alongside with
realizing this step.
}

\begin{lemma}[\textbf{Step I -- existence of the set isolating the bifurcation}] \label{lemma_big_isolation}

  There exists $\frac{1}{2} \leq \mu_+$ such that the semiflow
  associated with (\ref{eqn_ks_normal_form})
  is defined on $R := B_{\mu_+}$ \textcolor{black}{for every $\mu \geq \mu_+$} and
   \textcolor{black}
   {
 \begin{itemize}
     \item[(i)] for each $\mu \in [\mu_+, 1)$ there exists time $T > 0$ such that for all $t >T$ we have $\varphi_\mu(t, B_{\mu_+}) \subset B_\mu$,
    \item[(ii)] for each $\mu > 1$ and any positive $D \leq C\expr{\mu_+}$ there exists a time $T > 0$ such that for all $t >T$ we have $\varphi_\mu(t, B_{\mu_+}) \subset S_D$.
 \end{itemize}
   }
%   \begin{itemize}
%     \item[(i)] $B_{\mu_+}$ satisfies assumptions of Theorem \ref{thm:conv-on-Trap} for all $\mu$ close to $1$ and for each $\mu_+ < \mu < 1$ there exists time $T > 0$ such that for all $t >T$ we have $\varphi_\mu(t, B_{\mu_+}) \subset B_\mu$,
%    \item[(ii)] for each $\mu \geq 1$ and any positive $D \leq C\expr{\mu_+}$ there exists a time $T > 0$ such that for all $t >T$ we have $\varphi_\mu(t, B_{\mu_+}) \subset S_D$.
%   \end{itemize}

\end{lemma}
\begin{proof}
% \textbf{PZ: nie rozumiem co to ma oznaczac:  jaki "proceeding part"? \textcolor{red}{JK:Chodzi o to ze tam dalej w dowodzie wykazujemy strzalki do srodka, robienie tego tutaj wydaje mi sie mniej naturalne wiec pisze tylko ze za chwile to wykazemy}} 
 It will be evident from the \textcolor{black}{following argument} that the set $P_n R$ is forward-invariant for every $n > 1$ and condition (D) 
 \textcolor{black}{(see Definition \ref{def_condition_D})}  is manifestly satisfied because of (\ref{eqn_bound_norma_form_nonlinear_derivatives}). Thus the existence
  of the semiflow follows by Theorem \ref{thm:conv-on-Trap}.  
  \par
  First we fix $\mu$ such that $\mu_+ \leq \mu < 1$. Observe that since $\mu_+ > \frac{1}{2},$ the function
  $C(\mu)$ is decreasing to $0$ on the interval $\interval{\mu_+, 1}$.
  \par
  \textcolor{black}{
    We want to prove that there exists $T > 0$ such that $\varphi\expr{T, S_{\mu_+}} \subset S_\mu$. We cannot do this completely straightforwardly, because $C(\mu_+)$ may be much bigger than $C(\mu)$. This is why instead
    we will prove a stronger condition, namely that when we take any $\mu_+ \leq \theta < 1$ and $\theta' =\min\set{\frac{1 + \theta}{2}, \mu}$, there exists some time $T > 0$ such that for every $t > T$
  we have $\varphi(S_\theta, t) \subset S_{\theta'}$. }
  \textcolor{black}{Thanks to this choice $C(\theta)$ is of the order of $C(\theta')$.} \textcolor{black}{Indeed,} since $\frac{1}{2} \leq \mu_+ < 1$, we have
  \begin{align}
    \label{eq_same_order}
    \begin{split}
      1 \leq \frac{C(\theta)}{C(\theta')} &\leq \sqrt{\frac{\expr{1 - \theta}\expr{1 - 4\theta}}{\expr{1 - \frac{1 + \theta}{2}}\expr{1 - 4\frac{1 + \theta}{2}}}}\\
      &= \sqrt{\frac{2\expr{1 - 4\theta}}{\expr{1 - 4\frac{1 + \theta}{2}}}} < \sqrt{2} \\
    \end{split}
\end{align}
\par
Now consider $z \in B_{\mu_+}$ such that $\pi_{b_1} z \geq \sqrt{2}C(\theta')$ (case $\leq -\sqrt{2}C(\theta')$ is analogous).
If $\mu \geq \frac{1 + \theta}{2}$, then $\theta' = \mu$ and  $\frac{\lambda_2\expr{\theta'}}{\lambda_2\expr{\mu}} = 1$.
Otherwise, since $\theta \geq \mu_+ \geq \frac{1}{2}$, we have $\frac{1 + \theta}{2} \geq \frac{3}{4}$.
Thus we have that
$\frac{\lambda_2(\theta')}{\lambda_2(\mu)}
\geq \frac{\lambda_2\expr{\frac{1 + \theta}{2}}}{\lambda_2\expr{\mu}}
> \frac{\lambda_2(\frac{3}{4})}{\lambda_2(1)}
= \frac{2}{3}$, which by (\ref{eqn_bound_norma_form_nonlinear_parts}) and Remark \ref{rem_lambda1} gives us that
\begin{align*}
  \dot{b}_1 &=    b_1\expr{\lambda_1(\mu) + \frac{4}{\lambda_2(\mu)}b_1^2} + h_1(b_1, b_2, \ldots)\\
            &\leq \sqrt{2}C(\theta')\expr{\lambda_1(\mu) + \frac{4}{\lambda_2(\mu)} \expr{\frac{-2\lambda_1(\theta')\lambda_2(\theta')}{4}}} + \abs{h_1(b_1, b_2, \ldots)} \\
            &<    \sqrt{2}C(\theta')\expr{\lambda_1(\mu) - \frac{4}{3}\lambda_1(\theta')} + \alpha_1 C(\theta)^6 \\
            &< -\frac{\sqrt{2}}{3} C(\theta') \lambda_1(\theta') + 8 \alpha_1 C(\theta')^6 < 0, \numberthis \label{eq_condition_on_isolation_b1}
\end{align*}
if $\mu_+$ is sufficiently close to $1$.
\par

% If $\pi_{b_2} z \geq C(\theta')^3$ (case $\pi_{b_2} z \leq -C(\theta')^3$ is analogous), then
% \begin{align*}
%     \dot{b}_2 &= \lambda_2\expr{\mu}b_2 + h_2(b_1, b_2, \ldots) \leq \\
%               &\leq \lambda_2\expr{\mu}C(\theta')^3 + \abs{h_2(b_1, b_2, \ldots)} < \\
%               &< \lambda_2\expr{\mu}C(\theta')^3 + \alpha_2 C(\theta)^4 < \\
%               &< \lambda_2\expr{\mu}C(\theta')^3 + 4 \alpha_2 C(\theta')^4 < 0, \numberthis \label{eq_condition_on_isolation_b2}
% \end{align*}
% if $\mu_+$ is sufficiently close to $1.$
% \par
% If $\pi_{b_3} z \geq C(\theta')^3$ (case $\leq -C(\theta')^4$ is analogous), then
% \begin{align} \label{eq_condition_on_isolation_b3}
%   \dot{b}_3 < \lambda_3\expr{\mu} C(\theta')^3 + 8 \abs{d(\mu_+)} C(\theta')^3 + 4 \alpha_3 C(\theta')^4 < 0,
% \end{align}
% for $\mu_+$ close enough to $1$, because for all $x, y \in \interval{\frac{1}{2}, 1}$ we have $8 \abs{d(x)} + \lambda_3(y) < 0$.
% \par

  Let $z \in B_{\mu_+}$ be such that $\pi_{b_k} z \geq \frac{C(\theta')^4}{k^s}$ for some $k \geq 2$ (case $\leq -\frac{C(\theta')^4}{k^s}$ is analogous). Using the fact that $\abs{\frac{k^2(1 - \mu k^2)}{k^s}} \geq \frac{1}{k^{s-2}}$,
  for $k \geq 2$ and $\mu > \frac{1}{2}$, we have by (\ref{eqn_bound_gen_form_nonlinear_parts})
  \begin{align*}
      \dot{b}_k &= \lambda_k(\mu) b_k + h_k(b_1, b_2, \ldots) \\
                &\leq \lambda_k\expr{\mu} \frac{C(\theta')^3}{k^s} + \abs{h_k(b_1, b_2, \ldots)} \\
                &< \lambda_k\expr{\mu} \frac{C(\theta')^3}{k^s} + \frac{\alpha_k C(\theta)^4}{k^{s-2}} \leq \\
                &\leq \lambda_k\expr{\mu} \frac{C(\theta')^3}{k^s} + \frac{4\alpha_k C(\theta')^4}{k^{s-2}} < 0, \numberthis \label{eq_condition_on_isolation}
  \end{align*}
  for each $k \geq 2$  if $\mu_+$ is sufficiently close to $1$.
  \par
  It follows that there exists $M > 0$ such that for the sequence $\theta_0 = \mu_+, \theta_1, \theta_2, \ldots, \theta_M = \mu$ given by $\theta_i := \min \set{\frac{1 + \theta_{i - 1}}{2}, \mu}, i=1,\dots,M$ there exist
  times $t_1, t_2, \ldots, t_M$ such that for all $t > t_i$ we have $\varphi_\mu \expr{t_i, B_{\theta_{i-1}}} \subset B_{\theta_{i}}, i = 1,\ldots,M,$ when $\mu_+$ is sufficiently close to $1$.
  This concludes the proof when $\mu_+ \leq \mu < 1$.
  \par
  Now we fix $\mu \geq 1$ and a positive $D \leq C\expr{\mu_+}$. Observe that for some $\xi \in \left[\mu_+, 1\right)$
  we have $S_D = B_\xi$. Case of variables $k \geq 2$ is
  analogous as when $\mu_+ \leq \mu < 1$. For $k = 1$, first observe that $\lambda_1(\mu) \leq 0$ (in consequence, we will simply ignore this term in the inequality below). Now if for any $\theta \in \interval{\mu_+, \xi}$ we set $\theta' = \min \set{\xi, \frac{1 + \theta}{2}}$,
  then for $z \in B_{\mu_+}$ such that $\pi_{b_1} z \geq \sqrt{2}C(\theta_1)$ we have
  \begin{align*}
    \dot{b}_1 &=    b_1\expr{\lambda_1(\mu) + \frac{4}{\lambda_2(\mu)}b_1^2} + h_1(b_1, b_2, \ldots)\\
              &< -\frac{4\sqrt{2}}{3} C(\theta') \lambda_1(\theta') + 8 \alpha_1 C(\theta')^6 < 0, \numberthis \label{eq_condition_on_isolation_stable}
  \end{align*}
  if $\mu_+$ is close enough to $1$. The remaining part of the proof is analogous to the proof when $\mu_+ \leq \mu < 1$.
\end{proof}

Now we prove that before bifurcation the origin is not only attracting,
but also hyperbolic.

\begin{lemma} [\textbf{Step II-- hyperbolicity of the origin before the bifurcation}] \label{lemma_neg_lognorm_before_bif}
    If $\mu > 1$ is small enough, then
    $0$ is a hyperbolic attracting fixed point.
\end{lemma}
\begin{proof}

    Let $\mu > 1$. The set $S_{-\lambda_1\expr{\mu}}$ is forward-invariant
    if $\mu > 1$ is small enough by Lemma \ref{lemma_big_isolation}.
    Thus if we show that on this set the $l_2$
    logarithmic norm is negative, by Remark \ref{rem_stable_hyperbolicity} we will get the result.
    \par
    Since $\frac{24}{\lambda_2} < 0$,
    by Lemma (\ref{eqn_bound_norma_form_nonlinear_derivatives}) we have for some $l < 0$
    \begin{align*}
        2\frac{\partial \hat{F}_1}{\partial b_1} + \sum_{i=2}^{\infty} \abs{\frac{\partial \hat{F}_1}{\partial b_i}} + \abs{\frac{\partial \hat{F}_i}{\partial b_1}} &=
        2\lambda_1 + \frac{24}{\lambda_2}b_1^2 + 2\frac{\partial h_1}{\partial b_1} +  \sum_{i=2}^{\infty} \abs{\frac{\partial h_1}{\partial b_i}} + \abs{\frac{\partial h_i}{\partial b_1}} < \\
        &< 2\lambda_1 + \beta \abs{\lambda_1^3} < l,
    \end{align*}
    and
    \begin{align*}
        \frac{\partial \hat{F}_k}{\partial b_k} + \sum_{i \neq k} \abs{\frac{\partial \hat{F}_k}{\partial b_i}} + \abs{\frac{\partial \hat{F}_i}{\partial b_k}} &=
        2\lambda_k + 2\frac{\partial h_k}{\partial b_k} + \sum_{i \neq k} \abs{\frac{\partial h_k}{\partial b_i}} + \abs{\frac{\partial h_i}{\partial b_k}} = \\
        &< 2k^2\expr{1 - \mu k^2} + \beta_k k \expr{\mu - 1} < l,
    \end{align*}
    for all $k \geq 2$ if $\mu$ is small enough. Thus condition $(Dd)$ with $l < 0$ is verified.
    \par
\end{proof}

Let us remark that the dynamics for $\mu > 1$ is so simple that we could prove that the origin
is an attracting hyperbolic fixed point by taking a bit different shape of sets for
larger $\mu$, but we are not interested in this result here.
\par
Lemmas \ref{lemma_big_isolation}, \ref{lemma_neg_lognorm_before_bif} completely describe the behavior before the bifurcation.
Now, we proceed to describing what happens when we pass through $\mu = 1$. We first establish,
by checking the cone conditions, that the unstable manifold of the source point
is a graph over the $b_1$ direction.

% \begin{lemma} \label{lemma_local_isolation}
%   There exists $\mu_+$ such that for all $\mu_+ \leq \mu < 1$
%   for all $t > 0$ we have $\varphi_\mu(t, B_\mu) \subset B_\mu$.
% \end{lemma}
% \begin{proof}
%   Analogous to the proof of Lemma \ref{lemma_big_isolation}.
% \end{proof}

\begin{lemma}[\textbf{Step IIIa -- hyperbolicity of the origin after the bifurcation}] \label{lemma_unstable_manifold_in_ks}
  Let
  \begin{align*}
    B_\mu^0 &:= \set{z \in B_\mu \mid \abs{\pi_{b_1}z} \leq \frac{C(\mu)}{2}}
  \end{align*}
  Then there exists $\mu_+ > 0$ such that for all $\mu > \mu_+$
  there exists a unique fixed point   $z_0 \in B_\mu^0$. Moreover,
   $W^u_{B_\mu^0}(z_0)$ is the image of a horizontal disk in $B_\mu^0$.
  % There exists $\mu_+ > 0$ such that for each $1> \mu > \mu_+$ such that cone conditions are satisfied on $B_\mu^0$ with
  % \begin{align*}
  %   Q = \begin{bmatrix}
  %   1 & 0 & 0 &\dots \\
  %   0 & -1 & 0 &\dots \\
  %   0 & 0  & -1 &\dots \\
  %   \vdots & \vdots & \vdots & \ddots\\
  %   \end{bmatrix}
  % \end{align*}
\end{lemma}
\begin{proof}
    We will verify assumptions of Theorem \ref{thm_unstable_manifold_selfcb}.
    \par
    Fix $\mu \in [\mu_+, 1)$. $B_\mu^0$ is an $\varepsilon$-isolating cuboid by
    Lemmas \ref{lemma_big_isolation} and \ref{lemma_connection_sets}.
    Thus it is enough to verify assumptions
    of Lemma \ref{lemma_cc_verification}.

    By (\ref{eqn_bound_norma_form_nonlinear_derivatives}) we have
    \begin{align*}
        2\frac{\partial \hat{F}_1}{\partial b_1} - \sum_{i = 2}^\infty \expr{\abs{\frac{\partial \hat{F}_1}{\partial b_i}} + \abs{\frac{\partial \hat{F}_i}{\partial b_1}}} &\geq
        2\lambda_1 + \frac{24}{\lambda_2}b_1^2 + 2\frac{\partial h_1}{\partial b_1} -  \sum_{i = 2}^\infty \expr{ \abs{\frac{\partial h_1}{\partial b_i}} + \abs{\frac{\partial h_i}{\partial b_1}}} \\
        &\geq 2\lambda_1 + \frac{24}{\lambda_2}\frac{C^2}{4} - \beta_1 C^3.
    \end{align*}
    Observe that $2\lambda_1 + \frac{24}{\lambda_2}\frac{C^2}{4} = \frac{1}{2}\lambda_1 > 0$,
    thus by Remark \ref{rem_lambda1}
    \begin{align*}
        2\frac{\partial \hat{F}_1}{\partial b_1} - \sum_{i = 2}^\infty \expr{\abs{\frac{\partial \hat{F}_1}{\partial b_i}} + \abs{\frac{\partial \hat{F}_i}{\partial b_1}}}
        \geq \frac{1}{2}\lambda_1 -  \beta_1 C^3 > 0, \numberthis \label{eqn_cone_conditions_in_proof_of_bif1}
    \end{align*}
    if $\mu$ is close enough to $1$.
    \par
    Now for some $l > 0$ we have analogously as in the proof of Lemma \ref{lemma_neg_lognorm_before_bif}
    \begin{align*}
        -2\frac{\partial \hat{F}_k}{\partial b_k} - \sum_{i \neq k} \expr{ \abs{\frac{\partial \hat{F}_1}{\partial b_i}} + \abs{\frac{\partial \hat{F}_i}{\partial b_i}} } &= -2\lambda_k - \beta_k kC > l, \numberthis \label{eqn_cone_conditions_in_proof_of_bif2}
    \end{align*}
    for all $k \geq 2$ if $\mu_+$ is close enough to $1$.
\end{proof}

We now show that points on the unstable manifold of the source point (excluding the source itself) eventually go into $B_\mu^-$ or $B_\mu^+$.

\begin{lemma}[\textbf{Step IIIb -- connection sets}]\label{lemma_connection_sets}
  Define for $\mu > 0$
  \begin{align*}
    B_\mu^{c+} := \set{z \in B_\mu \mid \frac{C(\mu)}{2} \leq \pi_{b_1}z \leq \frac{C(\mu)}{\sqrt{2}}},
    B_\mu^{c-} := \set{z \in B_\mu \mid -\frac{C(\mu)}{2} \geq \pi_{b_1}z \geq -\frac{C(\mu)}{\sqrt{2}}}.
  \end{align*}
  Then for $\mu < 1$ close enough to $1$ we have $\dot{b}_1 > 0$ on $B_\mu^{c+}$ and
  $\dot{b}_1 < 0$ on $B_\mu^{c-}$.
\end{lemma}
\begin{proof}
  Fix $\mu > 0$. Assume $z \in B_\mu^{c+}$ (case of $z \in B_\mu^{c-}$ is analogous). Since we have
  $$
  \lambda_1  + \frac{4}{\lambda_2}b_1^2 \geq \lambda_1  + \frac{4}{\lambda_2}\frac{C^2}{2} = \frac{\lambda_1}{2},
  $$
   by (\ref{eqn_bound_norma_form_nonlinear_parts}) and Remark \ref{rem_lambda1} it follows that
  \begin{align*}
      \dot{b}_1 &= b_1(\lambda_1  + \frac{4}{\lambda_2}b_1^2) + h_1(b_1, b_2, \ldots) \\
                &\geq \frac{1}{2}C\frac{\lambda_1}{2} - \abs{ h_1(b_1, b_2, \ldots)} \\
                &\geq \frac{1}{4}C \lambda_1 - \alpha_1 C^5 > 0. \numberthis \label{eqn_connection_sets}
  \end{align*}
  % \textbf{NOTE here and in other places we can improve the inequalities, because we aren't on
  % the whole set $B_{\mu}$! But we would have to restate the normal form lemma.
  % Looking at the numbers it doesn't look like it is a bottleneck for the proof of
  % the bifurcation anyway.}
\end{proof}

Now we will prove that two attracting fixed points are born by establishing basins
of attraction near to the approximate fixed points of (\ref{eqn_ks_normal_form}).
\begin{lemma}[\textbf{Step IIIc -- hyperbolicity of the born fixed points}  ] \label{lemma_neg_lognorm_after_bif}
  Let
  \begin{align*}
    B_\mu^+ &:= \set{z \in B_\mu \mid \pi_{b_1}z \geq \frac{C(\mu)}{\sqrt{2}}} \\
    B_\mu^- &:= \set{z \in B_\mu \mid \pi_{b_1}z \leq -\frac{C(\mu)}{\sqrt{2}}}
  \end{align*}
  Then for $\mu < 1$ close enough to $1$ there exist unique fixed points $z_\mu^\pm \in B_\mu^\pm$
  of $\varphi_{\mu}$.
  Moreover, those fixed points are hyperbolic and attracting and those are \textcolor{black}{the}  only points
  which have  \textcolor{black}{a}  full backward trajectory respectively in $B_\mu^\pm$.
  \begin{proof}
      Let $\mu < 1$. $B_\mu^\pm$ are easily seen to be forward-invariant by Lemma \ref{lemma_big_isolation}.
      and Lemma \ref{lemma_connection_sets}.
      Using (\ref{eqn_bound_norma_form_nonlinear_derivatives}) on $B_\mu^+ \cup B_\mu^-$ we have by Remark \ref{rem_lambda1}
      \begin{align*}
          2\frac{\partial \hat{F}_1}{\partial b_1} + \sum_{i=2}^{\infty} \abs{\frac{\partial \hat{F}_1}{\partial b_i}} + \abs{\frac{\partial \hat{F}_1}{\partial b_1}} &=
          2\lambda_1 + \frac{24}{\lambda_2}b_1^2 + 2\frac{\partial h_1}{\partial b_1} + \sum_{i=2}^{\infty} \abs{\frac{\partial h_1}{\partial b_i}} + \abs{\frac{\partial h_i}{\partial b_1}} \leq \\
          &\leq 2 \lambda_1 + \frac{24}{\lambda_2}\frac{C^2}{2} + \beta_1 C^3 = \\
          &= -\lambda_1 + \beta_1 C^3 < \\
          &< 0, \numberthis \label{eqn_neg_lognorm_condition_1}
      \end{align*}
      when $\mu$ is close to $1$.
     For some $l < 0$ we also have
      \begin{align*}
          2\frac{\partial \hat{F}_k}{\partial b_k} + \sum_{i \neq k} \abs{\frac{\partial \hat{F}_k}{\partial b_i}} + \abs{\frac{\partial \hat{F}_i}{\partial b_k}} &=
          2\lambda_k + 2\frac{\partial h_k}{\partial b_k} + \sum_{i \neq k} \abs{\frac{\partial h_k}{\partial b_i}} + \abs{\frac{\partial h_i}{\partial b_k}} = \\
          &< 2\lambda_k + \beta_k k C < l, \numberthis \label{eqn_neg_lognorm_condition_2}
      \end{align*}
      for all $k \geq 2$ if $\mu$ is close to $1$.
      \par
      The thesis follows by Theorem \ref{thm:attractfp} and Remark \ref{rem_stable_hyperbolicity}.
  \end{proof}
\end{lemma}

Now we are ready to prove the main result in this section.

\begin{theorem} \label{thm_bif}
    There exist $\mu_+ < 1 < \mu_-$ such that the
    pitchfork bifurcation occurs in (\ref{eqn_ks_normal_form}) on $\interval{\lambda_1\expr{\mu_-}, \lambda_1\expr{\mu_+}}$.
\end{theorem}
\begin{proof}
    (P1) is true by Lemmas \ref{lemma_big_isolation}, \ref{lemma_neg_lognorm_before_bif}, \ref{lemma_unstable_manifold_in_ks}.
    (P2) is true by Lemma \ref{lemma_big_isolation}.
    \par
    We now prove (P3) is true. By Lemma \ref{lemma_unstable_manifold_in_ks}
    there exists $z \in B_\mu$ such that $\pi_{b_1} z \in \pi_{b_1} B_\mu^0 \cap B_\mu^{c+}$ which has
    a backward trajectory to $0$. Assume that $\pi_{b_1} z > 0$ (other case is analogous).
    By Lemma \ref{lemma_connection_sets} there exists a time $t_1 > 0$ such that $\varphi(t_1, z) \in B_\mu^+$
    and by Lemma \ref{lemma_neg_lognorm_after_bif} we have $\lim_{t \to \infty} \varphi\expr{t, \varphi\expr{t_1, z}} = z_\mu^+$. Proof for $z_\mu^-$ is analogous. By Lemma \ref{lemma_neg_lognorm_after_bif},
    $z_\mu^\pm$ are also hyperbolic.
    \par
    It remains to show that we have described the maximal
    invariant set. By Theorem \ref{thm_unstable_manifold}, only points in $B_\mu^0$ which can
    be in the maximal invariant set are those on the unstable manifold (because only those
    have \textcolor{black}{a}  full backward trajectory). By Lemma \ref{lemma_connection_sets}, any full backward
    trajectory of a point in $B_\mu^{c\pm}$ must have a point in $B_\mu^0$, so it must
    have a point from $W^u_{B_{\mu}^0}(0)$ by what we have just proved. Now, by Lemma
    \ref{lemma_neg_lognorm_after_bif}, all points in $B_\mu^{\pm}$ with full backward trajectory
    either leave $B_\mu$ or have a point from $B_\mu^{c\pm}$ on a backward trajectory.
    Conclusion easily follow by considering in turn which points can have a full forward trajectory.
\end{proof}

 \section{Inequalities for the verification of a given bifurcation range}
\label{sec_gen_bif_thm}

The analytical proof of  \textcolor{black}{a}  pitchfork bifurcation in the Kuramoto--Sivashinsky equation we have
just provided does not tell us in what range of the parameters the established dynamics is valid.
We now want to extract from the proof of Theorem \ref{thm_bif} the conditions general enough
to use them in the computer assisted proofs. We also add a few conditions to account for
the unstable directions. We use those conditions in Section \ref{sec_appendix_data_from_the_proof_of_bifurcation}
to give an explicit range when $\mu = 1$ and to prove the bifurcation (with an unstable direction) when $\mu = 0.25$.
\par
We again work on the sets of the form
\begin{align}
    S_C := \interval{-\zeta C, \zeta C} \times \prod_{i = 2}^\infty \interval{-\frac{C^\omega}{i^s}, \frac{C^\omega}{i^s}},   \label{eq:SCdef2}
\end{align}
where $\omega \geq 3$ and $\zeta > 1$.
\par
Equation we  consider is
\begin{align} \label{eqn_gen_bif_form}
    \begin{split}
        a_1' &= \lambda_1\expr{\mu}a_1 + c\expr{\mu}a_1^3 + h_1\expr{a}, \\
        a_k' &= \lambda_k\expr{\mu}a_k + h_k\expr{a}, \quad k > 1.
    \end{split}
\end{align}

\textcolor{black}{We introduce some notation first. Let $\mu_1 < \mu_2$  For any function of parameter $f(\mu)$ we will
use in this subsection we define $ \overline{f}\expr{\interval{\mu_1, \mu_2}} := \max_{\interval{\mu_1, \mu_2}}f(\mu)$ 
and $\underline{f}\expr{\interval{\mu_1, \mu_2}} := \min{\interval{\mu_1, \mu_2}}f(\mu)$.
We also shorten the notation for the bound on the entire parameter range on which $f$ is defined.
Throughout this subsection we also assume that those bounds exist and are finite.
In the case we use for our assisted proof all of maxima and minima will be realized at the interval's
endpoints or at the bifurcation parameter $\mu_b$ due to monotonicity, but it is not necessary to assume it in general.}

Assume that $\lambda_1$ is a continuous increasing function, $\lambda_1\expr{\mu_-} < 0$
and $\textcolor{black}{ \overline{\lambda_1}} > 0$ for some $\mu_\pm \in \mathbb{R}$. As before, $\lambda_1$
will be our bifurcation parameter. Note that here we assume that
$\mu_+ > \mu_-$, so the direction in which $\lambda_1$ grows is reversed compared to the KS equation,
but the modification is straightforward. Let $\mu_b \in \interval{\mu_-, \mu_+}$ be
the solution of $\lambda_1\expr{\mu_b} = 0$.
Assume that $c\expr{\mu}$ is a negative  function \textcolor{black}{bounded away from zero on $\interval{\mu_-, \mu_+}$}  
 that there exists $m \geq 1$ such that for $\mu \in \interval{\mu_-, \mu_+}$
\begin{align*}
    \lambda_i\expr{\mu} &> 0, \quad \textnormal{for } 1 < i \leq m, \\
    \lambda_i\expr{\mu} &< 0, \quad \textnormal{for } i > m.
\end{align*}
Case $m = 1$ corresponds to the case of no unstable directions.
We also denote for $\mu \geq \mu_b$
\begin{align*}
    C\expr{\mu} := \sqrt{\frac{\lambda_{1}\expr{\mu}}{-c\expr{\mu}}},
\end{align*}
so that when $\abs{a_1} = C\expr{\mu}$ we have $\lambda_1\expr{\mu}a_1 + c\expr{\mu}a_1^3 = 0$.
\par

We assume that there exist $\alpha_1, \alpha_2, \dots > 0$ such that
\begin{align} \label{eqn_bound_gen_form_nonlinear_parts}
    \begin{split}
    \abs{h_k} &< \frac{\alpha_k C\expr{\mu}^{n_k}}{k^{p}}, \quad k = 1, 2 \ldots
\end{split}
\end{align}
and $\beta_1, \beta_2, \dots > 0$
such that
\begin{align}\label{eqn_bound_gen_form_nonlinear_derivatives}
    \begin{split}
    2\abs{\frac{\partial h_k}{\partial b_k}} + \sum_{i \neq k} \abs{\frac{\partial h_k}{\partial b_i}} + \abs{\frac{\partial h_i}{\partial b_k}} &< \beta_k k C^{d_k}\expr{\mu}, \quad k = 1, 2, \ldots
    \end{split}
\end{align}
on $S_{C\expr{ \overline{\mu}}}$ for all $\mu \in \interval{\mu_-, \mu_+},  \overline{\mu} \in (\mu_b, \mu_+]$. 
\par
With the assumptions above, we have similarly to Lemma \ref{lemma_ks_norm_form_F} the lemma below.
\begin{lemma}
        Denote $d_{ij} := \max_{x \in S_{ \textcolor{black}{ \overline{C}}}} \abs{\frac{\partial \hat{F}_{\mu, i}}{\partial x_j}}$.
        Then the conditions (F) and (\ref{eq_cone_cond_convergence1}) are satisfied for (\ref{eqn_gen_bif_form}) on $S_{ \textcolor{black}{ \overline{C}}}$.
\end{lemma}

\par
Let $K > 1, l > 0, 1 > \gamma_+ > \gamma_- > 0$. Assume also that on the set

$$\interval{-\gamma_- C, \gamma_- C} \times \prod_{i = 2}^\infty \interval{-\frac{C^\omega}{i^s}, \frac{C^\omega}{i^s}}$$

we have additional derivatives bound

\begin{align*}
  2\abs{\frac{\partial h_1}{\partial b_1}} + \sum_{i \neq 1} \abs{\frac{\partial h_1}{\partial b_i}} + \abs{\frac{\partial h_i}{\partial b_1}} &<  \overline{\beta_1} C^{d_1\expr{\mu}},
\end{align*}
for some $0 <  \overline{\beta_1} < \beta_1$. Consider the following conditions.

\textcolor{black}{Now we list a seemingly intimidating list of inequalities which we assume in 
our theorem. Nevertheless those inequalities arise quite naturally in the proof and after stating all of them
we summarize the dynamical character of each inequality.}
\begin{align}
    \frac{\textcolor{black}{\underline{\lambda_k}\expr{\interval{\mu_b, \mu_+}}}}{k^s} - \frac{\alpha_k K^{n_k} \textcolor{black}{ \overline{C}}^{n_k - \omega}}{k^{p}}  &> 0, \quad k = 2, \dots, m, \label{eqn_cond1} \\
    1 \leq \frac{C(\mu)}{C\expr{\frac{\mu_b + \mu}{2}}} &\leq K, \quad \mu \in (\mu_b, \mu_+] \label{eqn_cond2} \\
    \expr{1 - \frac{\textcolor{black}{ \overline{c}\expr{\interval{\mu_b, \mu_+}}}}{ \textcolor{black}{ \overline{C}}} \zeta^2} \zeta  + K^{n_1} \alpha_1 \frac{ \textcolor{black}{ \overline{\lambda_1}}^{\frac{n_1 - 3}{2}}}{\expr{-\textcolor{black}{ \overline{c}\expr{\interval{\mu_b, \mu_+}}}}^{\frac{n_1 - 1}{2}}} &< 0, \label{eqn_cond4}  \\
    \frac{\textcolor{black}{ \overline{\lambda_k}}}{k^s} + \frac{\alpha_k K^{n_k} \textcolor{black}{ \overline{C}}^{n_k - \omega}}{k^p} &< 0, \quad k > m, \label{eqn_cond5}  \\
    -\textcolor{black}{2} + \beta_1 \textcolor{black}{\abs{\underline{\lambda_1}}}^{\omega - 1} &< 0, \label{eqn_cond6} \\
    2\textcolor{black}{\underline{\lambda_k}\interval{\mu_b, \mu_+}} - \beta_k k  \textcolor{black}{\abs{\underline{\lambda_1}}} &> l, \quad k = 2, \dots, m, \label{eqn_cond8}  \\
    2\textcolor{black}{ \overline{\lambda_k}\expr{\interval{\mu_b, \mu_+}}} + \beta_k k  \textcolor{black}{\abs{\underline{\lambda_1}}} &< -l, \quad k > m, \label{eqn_cond9}  \\
    (2 - 6\gamma_-^2) -  \overline{\beta_1} \frac{ \textcolor{black}{ \overline{\lambda_1}}^{\frac{d_1 - 2}{2}}}{\expr{-\textcolor{black}{ \overline{c}\expr{\interval{\mu_b, \mu_+}}}}^{\frac{d_1}{2}}} &> 0, \label{eqn_cond10} \\
    (2 - 6\gamma_+^2) + \beta_1 \frac{ \textcolor{black}{ \overline{\lambda_1}}^{\frac{d_1 - 2}{2}}}{\expr{-\textcolor{black}{ \overline{c}\expr{\interval{\mu_b, \mu_+}}}}^{\frac{d_1}{2}}} &< 0, \label{eqn_cond10'} \\
    2\textcolor{black}{\underline{\lambda_k}\expr{\interval{\mu_b, \mu_+}}} - \beta_k k  \textcolor{black}{ \overline{C}} &> l, \quad k = 2, \dots m \label{eqn_cond11}\\
    2\textcolor{black}{ \overline{\lambda_k}} + \beta_k k  \textcolor{black}{ \overline{C}} &< -l, \quad k > m, \label{eqn_cond12}\\
    \gamma (1 - \gamma^2) - \alpha_1 \frac{ \textcolor{black}{ \overline{\lambda_1}}^\frac{n_1 - 3}{2}}{\expr{-\textcolor{black}{ \overline{c}\expr{\interval{\mu_b, \mu_+}}}}^\frac{n_1 - 1}{2}} &> 0 \; \textnormal{for all } \gamma \in \interval{\gamma_-, \gamma_+}, \label{eqn_cond13} \\
    \max\set{\abs{2 - 6\gamma_-^2}, \abs{2 - 6\gamma_+^2}}\textcolor{black}{ \overline{\lambda_1}} + \beta_1 \textcolor{black}{ \overline{C}}^{d_1} &< \frac{l}{4}. \label{eqn_cond14}
\end{align}
% 10 + \beta_1 \frac{\textcolor{red}{\lambda_1\expr{\mu_+}} \textcolor{black}{ \overline{\lambda_1}}^{\frac{d_1 - 3}{2}}}{\expr{-\textcolor{red}{c\expr{\mu_b}}\textcolor{black}{ \overline{c}\expr{\interval{\mu_b, \mu_+}}}}^{\frac{d_1 - 1}{2}}} &< \frac{l}{4\textcolor{red}{\lambda_1\expr{\mu_+}} \textcolor{black}{ \overline{\lambda_1}}}, \label{eqn_cond14}
% \frac{1}{2} - \frac{1}{2 \zeta } - \alpha_1 \abs{\frac{\textcolor{red}{\lambda_1\expr{\mu_+}} \textcolor{black}{ \overline{\lambda_1}}^{\frac{k_1 - 3}{2}}}{\expr{-\textcolor{red}{c\expr{\mu_b}}\textcolor{black}{ \overline{c}\expr{\interval{\mu_b, \mu_+}}}}^{\frac{k_1 - 1}{2}}}} &> 0, \label{eqn_cond13}\\
% \frac{1}{2} - \beta_1 \abs{\frac{\textcolor{red}{\lambda_1\expr{\mu_+}} \textcolor{black}{ \overline{\lambda_1}}^{\frac{1}{2}}}{\expr{-\textcolor{red}{c\expr{\mu_b}}\textcolor{black}{ \overline{c}\expr{\interval{\mu_b, \mu_+}}}}^\frac{3}{2}}} &> 0, \label{eqn_cond10} \\
Observe that if there exist $\alpha, \beta, N > 0$ such that $\alpha_k = \alpha, \beta_k = \beta$
for $k > N$, then it is easy to see that to verify the conditions (\ref{eqn_cond5}, \ref{eqn_cond9}, \ref{eqn_cond12}),
it is enough to verify them for $k = m + 1, \dots, N + 1$.
\par
Throughout the remainder of this section we assume that all of the conditions above hold.
We also denote as before
$$B_\mu = S_{C\expr{\mu}}$$
for $\mu > \mu_b$.
\textcolor{black}{Let us list a summary of dynamical properties for whose verification those conditions are used
\begin{itemize}
\item condition (\ref{eqn_cond1}): used twice in Lemma~\ref{lemma_big_isolation_gen} to check the isolation of macroscopic cuboid in the unstable dierections,
\item condition  (\ref{eqn_cond2}):  used once in Lemma~\ref{lemma_big_isolation_gen} - to ensure that 
  $C(\theta)$ is of the same order as $C(\theta')$ when $\theta$ and $\theta'$ are close to each other. 
  It allows us to prove that when $\theta > \theta'$, points from $S_\theta$ flow into $S_{\theta'}$, 
  which we need in order to check the isolation of macroscopic cuboid. See also remark about (\ref{eq_same_order}) in proof of Lemma \ref{lemma_big_isolation}.
\item condition  (\ref{eqn_cond4}):  used once in Lemma~\ref{lemma_big_isolation_gen} - to check the isolation of macroscopic cuboid in the bifurcation direction.
\item (\ref{eqn_cond5} ):  used once in Lemma~\ref{lemma_big_isolation_gen} - to check the isolation of macroscopic cuboid for the tail
\item (\ref{eqn_cond6}):   used once in Lemma~\ref{lemma_neg_lognorm_before_bif_gen} -
 needed to check the cone conditions near the origin before the bifurcation,
More precisely to check that the condtion (\ref{eq_cc-diag-dom}) holds when $i$ is the bifurcation 
coordinate and when the bifurcation direction is treated as stable,
\item (\ref{eqn_cond8}): used once in Lemma~\ref{lemma_neg_lognorm_before_bif_gen} 
needed to check the cone conditions near the origin before the bifurcation. 
More precisely to check that the condtion (\ref{eq_cc-diag-dom}) holds when $i$ is an unstable  
coordinate and when the bifurcation direction is treated as stable,
\item (\ref{eqn_cond9}):  used once in Lemma~\ref{lemma_neg_lognorm_before_bif_gen} - 
needed to check the cone conditions near the origin before the bifurcation. 
More precisely to check that the condtion (\ref{eq_cc-diag-dom}) holds when $i$ is a stable  
coordinate and when the bifurcation direction is treated as stable 
\item (\ref{eqn_cond10}): used once in Lemma~\ref{lemma_unstable_manifold_in_gen} - same as (\ref{eqn_cond6}) but after the bifurcation, with the bifurcation direction treated as unstable,
\item (\ref{eqn_cond10'}):  used once in Lemma~\ref{lemma_neg_lognorm_after_bif_gen} - same as (\ref{eqn_cond10}) but for the born fixed points with the bifurcation direction treated as stable,
\item (\ref{eqn_cond11}):   used once in Lemma~\ref{lemma_unstable_manifold_in_gen} and implicitly used in Lemma~\ref{lemma_neg_lognorm_after_bif_gen} - same as (\ref{eqn_cond8}) but after the bifurcation, with the bifurcation direction treated as unstable,
\item (\ref{eqn_cond12}):  used once in Lemma~\ref{lemma_unstable_manifold_in_gen} and implicitly used in Lemma~\ref{lemma_neg_lognorm_after_bif_gen} - same as (\ref{eqn_cond9}) but after the bifurcation, with the bifurcation direction treated as unstable,
\item (\ref{eqn_cond13}):  used once in Lemma~\ref{lemma_connection_sets_gen} - to obtain estimate $a_1'>0$ in the connecting region
\item (\ref{eqn_cond14}): used once in Lemma~\ref{lemma_bif_gen_cone_field} - cone invariance in the connecting sets
 \end{itemize}
}

\begin{lemma}[\textbf{Step I -- existence of the set isolating the bifurcation}] \label{lemma_big_isolation_gen}
    Let $R := B_{\mu_+}$. For $\mu \in \interval{\mu_-, \mu_+}$ there exists a local continuous
    semiflow $\varphi_\mu$ associated with (\ref{eqn_gen_bif_form}).
   \textcolor{black}
   {
    Moreover
 \begin{itemize}
     \item[(i)] for each $\mu \in (\mu_b, \mu_+]$ and for all $x \in R$ and for all $\xi \geq \mu$ we have that $x$ leaves $R$ or there exists $t > 0$ such that $\varphi(t, x) \in B_\xi$.
    \item[(ii)] for each $\mu \in \interval{\mu_-, \mu_b}$ and for all $x \in R$ and for all $\xi > 0$ we have that $x$ leaves $R$ or there exists $t > 0$ such that $\varphi(t, x) \in B_\xi$.
 \end{itemize}
   }
\end{lemma}
\begin{proof}
  Fix $\mu \in \interval{\mu_-, \mu_+}$. To verify the assumptions of Theorem \ref{thm:local-conv},
  we first prove that for $n > m$
  the set is $P_n B_{\mu_+}$ is an $\varepsilon$-isolating cuboid. Let $k = 2, \dots, m$
  and assume that $z \in B_{\mu_+}$ is such that $\pi_{a_k} z = \frac{C^\omega}{k^s}$ (case $= -\frac{C^\omega}{k^s}$ is
  analogous). Then by (\ref{eqn_cond1})
  \begin{align} \label{eqn_help_isolation}
      a_k' = \lambda_k a_k + h_k\expr{a} > \lambda_k\expr{\mu} \frac{C\expr{\mu_b}^\omega}{k^s} - \frac{\alpha_k}{k^{p}} C\expr{\mu_b}^{n_k} > 0.
  \end{align}
  We see that $a_k'$ is uniformly separated from $0$ on the whole exit set of $P_n B_{\mu_+}$,
  so this proves the (I2) condition.
  The (I1) condition will be discussed in the further part of the proof.
  The condition (D) is manifestly satisfied because of (\ref{eqn_bound_gen_form_nonlinear_derivatives}).
  \par
  Now fix $\mu \in (\mu_b, \mu_+]$.
  Consider some $\theta \in \interval{\mu, \mu_+}$ and denote $\theta' = \max \set{\frac{\mu_b + \theta}{2}, \mu}$.
  We will show that there exists some time $T > 0$ such that for every $t > T$
  we have $\varphi(t, B_\theta) \subset B_{\theta'}$.
  We have by (\ref{eqn_cond2}) that
  \begin{align*}
       1 \leq \frac{C(\theta)}{C(\theta')} &\leq K. %\sqrt{\frac{\lambda_{1}\expr{\theta}c\expr{\theta'}}{c\expr{\theta}\lambda_1\expr{\theta'}}}
\end{align*}
\par
First we will prove, by checking the signs of the vector field
in the bifurcation and stable directions, that for $z \in B_\theta \setminus B_{\theta'}$
if for all $k = 2, \dots, m$ and for sufficiently long time $T > 0$ we have for
all $t \leq T$
$\pi_{a_k} \varphi\expr{t, z} < \frac{C(\theta')^{3}}{k^s}$,
then we have $\varphi\expr{T, B_\theta} \subset B_{\theta'}$.
Then we will prove that if in turn we have $\pi_{a_k} z \geq \frac{C(\theta')^{3}}{k^s}$ for some $k \in \set{2, \dots, m}$,
then $z$ leaves $R$.
If we prove those two things, then the claim follows very similarly as in the proof of Lemma \ref{lemma_big_isolation}.
\par
Consider $z \in B_{\theta}$ such that $\pi_{a_1} z \geq \zeta C(\theta')$ (case $\leq -\zeta C(\theta')$ is analogous).
Since $\lambda_1\expr{\mu} \leq \lambda_1\expr{\theta'},$ we have by (\ref{eqn_cond4})
\begin{align*}
  \dot{a}_1 &= a_1\expr{\lambda_1(\mu) + c\expr{\mu}a_1^2} + h_1(a_1, a_2, \ldots)\\
            &\leq \zeta  C(\theta')\expr{\lambda_1(\mu) + \zeta^2 c\expr{\mu} C\expr{\theta'}^2} + \alpha_1 K^{n_1} C\expr{\theta'}^{n_1} \\
            &< \left(1 - \frac{\textcolor{black}{ \overline{c}\expr{\interval{\mu_b, \mu_+}}}}{ \textcolor{black}{ \overline{C}}} \zeta^2\right)\zeta  C(\theta') \lambda_1(\theta') +  \alpha_1 K^{n_1} C\expr{\theta'}^{n_1} < 0.
\end{align*}
\par
  Let $z \in B_{\theta}$ be such that $\pi_{a_k} z \geq \frac{C(\theta')^{3}}{k^s}$ for some $k > m$ (case $\leq -\frac{C(\theta')^{n_k}}{k^s}$ is analogous). By (\ref{eqn_cond5}) we have
  \begin{align*}
      \dot{a}_k &= \lambda_k(\mu) a_k + h_k(a_1, a_2, \ldots) \\
                &\leq \lambda_k\expr{\mu} \frac{C(\theta')^\omega}{k^s} + \abs{h_k(a_1, a_2, \ldots)} \\
                &< \lambda_k\expr{\mu} \frac{C(\theta')^\omega}{k^s} + \frac{\alpha_k K^{n_k} C(\theta')^{n_k}}{k^p} < 0.
  \end{align*}
\par
  Those inequalities easily imply first of the two claims mentioned above. Now
  we proceed to the second one, namely that points with sufficiently large unstable direction (relatively to
  the other directions) leave.
  \par
  Let $z \in B_\theta$ be such that $\pi_{a_k} z \geq \frac{C(\theta')^{3}}{k^s}$
  for some $k \in \set{2, \dots, m}$ (case $\leq -\frac{C(\theta')^{n_k}}{k^s}$ is analogous).
  By (\ref{eqn_cond1}) we have
  \begin{align*}
      \dot{a}_k &= \lambda_k(\mu) a_k + h_k(a_1, a_2, \ldots) \\
                &\geq \lambda_k\expr{\mu} \frac{C(\theta')^\omega}{k^s} - \abs{h_k(a_1, a_2, \ldots)} \\
                &\geq \lambda_k\expr{\mu} \frac{C(\theta')^\omega}{k^s} - \frac{\alpha_k K^{n_k} C(\theta')^{n_k}}{k^p} > 0.
  \end{align*}
  This proves that for any $\theta \geq \mu$ if $\pi_{a_k} z \geq \frac{C(\theta')^{3}}{k^s}$,
  then for some $t > 0$ we have $\varphi\expr{t, z} \in B_{\min\set{2\theta, \mu_+}}$ and
  $\abs{\pi_{a_k} \varphi\expr{t, z}} \geq \frac{C(\theta)^{3}}{k^s}$.
  It is obvious by the same reasoning we leave $B_{2\theta}$ and so on, until
  we leave $R$.
\end{proof}

\begin{lemma}[\textbf{Step II -- hyperbolicity of the origin before the bifurcation}] \label{lemma_neg_lognorm_before_bif_gen}
    The origin is a hyperbolic
    fixed point for $\mu \in [\mu_-, \mu_b)$.
\end{lemma}
\begin{proof}
    Let $\mu \in [\mu_-, \mu_b)$. We will show that on $S_{-\lambda_1\expr{\mu}}$ assumptions of
    Theorem \ref{thm_unstable_manifold_selfcb} are satisfied. We have by (\ref{eqn_cond6})
    \begin{align*}
        \textcolor{black}{2}\frac{\partial F_1}{\partial a_1} + \sum_{i=2}^{\infty} \abs{\frac{\partial F_1}{\partial a_i}} + \abs{\frac{\partial F_i}{\partial a_1}} &=
        2\lambda_1 + 6c a_1^2 + 2\frac{\partial h_1}{\partial a_1} +  \sum_{i=2}^{\infty} \abs{\frac{\partial h_1}{\partial a_i}} < \\
        &< \textcolor{black}{2}\lambda_1 + \beta_1 \abs{\lambda_1}^\omega < 0,
    \end{align*}
    (we ignored the negative term $6ca_1^2$).
    \par
    For $k = 2, \dots, m$ we have by (\ref{eqn_cond8})
    \begin{align*}
        2\frac{\partial F_k}{\partial a_k} - \sum_{i \neq k} \abs{\frac{\partial F_k}{\partial a_i}} + \abs{\frac{\partial F_i}{\partial a_k}} &=
        2\lambda_k + 2\frac{\partial h_k}{\partial a_k} - \sum_{i \neq k} \abs{\frac{\partial h_k}{\partial a_i}} > \\
        &> 2\lambda_k - \beta_k k \abs{\lambda_1} > l.
    \end{align*}
    For $k > m$ we have by (\ref{eqn_cond9})
    \begin{align*}
        \textcolor{black}{2}\frac{\partial F_k}{\partial a_k} + \sum_{i \neq k} \abs{\frac{\partial F_k}{\partial a_i}} + \abs{\frac{\partial F_i}{\partial a_k}}
            < 2\lambda_k + \beta_k k \abs{\lambda_1} < -l.
    \end{align*}
\end{proof}

\begin{lemma} [\textbf{Step IIIa -- hyperbolicity of the origin after the bifurcation}] \label{lemma_unstable_manifold_in_gen}
  Let
  \begin{align*}
    B_\mu^0 &:= \set{z \in B_\mu \mid \abs{\pi_{a_1}z} \leq \gamma_- C(\mu)}
  \end{align*}
  Then all $\mu \in (\mu_b, \mu_+]$
  there exists a unique fixed point   $z_0 \in B_\mu^0$. Moreover,
   $W^u_{B_\mu^0}(z_0)$ is the image of a horizontal disk in $B_\mu^0$.
\end{lemma}
\begin{proof}
    We will verify assumptions of Theorem \ref{thm_unstable_manifold_selfcb} \textcolor{black}{on $B_\mu^0$}. We have by (\ref{eqn_cond10})
    \begin{align*}
        2\frac{\partial F_1}{\partial a_1} - \sum_{i = 2}^\infty \expr{\abs{\frac{\partial F_1}{\partial a_i}} + \abs{\frac{\partial F_i}{\partial a_1}}} &\geq
        2\lambda_1 + 6c a_1^2 + 2\frac{\partial h_1}{\partial a_1} -  \sum_{i = 2}^\infty \expr{ \abs{\frac{\partial h_1}{\partial a_i}} + \abs{\frac{\partial h_i}{\partial a_1}}} \\
        &\geq (2 - 6\gamma_-^2)\lambda_1 -  \overline{\beta}_1 C^{d_1} > 0.
    \end{align*}
    % &\geq 2\lambda_1 - \frac{3}{2}\lambda_1 - \beta_1 C^3 > 0

    For $k = 2, \dots, m$ we have by (\ref{eqn_cond11})
    \begin{align*}
        2\frac{\partial F_k}{\partial a_k} - \sum_{i \neq k} \expr{ \abs{\frac{\partial F_1}{\partial a_i}} + \abs{\frac{\partial F_i}{\partial a_i}} } &> 2\lambda_k - \beta_k kC > l.
    \end{align*}
    For $k > m$ we have by (\ref{eqn_cond12})
    \begin{align*}
        2\frac{\partial F_k}{\partial a_k} + \sum_{i \neq k} \expr{ \abs{\frac{\partial F_1}{\partial a_i}} + \abs{\frac{\partial F_i}{\partial a_i}} } &\textcolor{black}{<} 2\lambda_k + \beta_k kC < -l.
    \end{align*}
\end{proof}

\begin{lemma} [\textbf{Step IIIb -- connection sets, growth in the bifurcation direction}] \label{lemma_connection_sets_gen}
  Define for $\mu > \mu_b$
  \begin{align*}
    B_\mu^{c+} &:= \set{z \in B_\mu \mid \gamma_- C(\mu) \leq \pi_{a_1}z \leq \gamma_+ C(\mu)},\\
    B_\mu^{c-} &:= \set{z \in B_\mu \mid -\gamma_+ C(\mu) \geq \pi_{a_1}z \geq -\gamma_+ C(\mu)}.
  \end{align*}
  Then for $\mu \in (\mu_b, \mu_+]$ on those sets we have respectively
  \begin{align*}
    a_1' &> 0, \\
    a_1' &< 0. \\
  \end{align*}
\end{lemma}
\begin{proof}
  Fix $\mu > \mu_b$ and $\gamma \in \interval{\gamma_-, \gamma_+}$.
  By (\ref{eqn_cond13}) we have on $B_{\mu, \gamma}^{c+}$
  we have for $b_1 = \gamma C$
  \begin{align*}
    a_1' &= \lambda_1 a_1 + ca_1^3 + h_1\expr{a} \\
         &= a_1 \expr{\lambda_1 + c \gamma^2 C^2} + h_1\expr{a} \\
         &= a_1 \expr{\lambda_1 - c \gamma^2 \frac{\lambda_1}{c}} + h_1\expr{a}  \\
         &= \gamma (1 - \gamma^2) C \lambda_1 - \alpha_1 C^{n_1} > 0.
  \end{align*}

\end{proof}

\begin{lemma} [\textbf{Step IIIb -- connection sets, cone invariance}] \label{lemma_bif_gen_cone_field}
    Define $B^{ct\pm}_\mu := B^{c\pm}_\mu \cup B^{+}_\mu$.
    Then for $\mu \in (\mu_b, \mu_+]$ the sets $B_\mu^{ct\pm}$ are $\varepsilon-$isolating cuboids. Moreover,
  there exists a time $T > 0$ such that for all horizontal disks $h$ with respect
  to the cones given by the matrix
  \begin{align*}
      Q_{11} &= -1, \\
      Q_{ii} &= 1, \quad i = 2, \dots, m, \\
      Q_{ii} &= -1, \quad i > m,
  \end{align*}
  set $\varphi\expr{t, \im h}$ contains the image of another horizontal disk.
\end{lemma}
\begin{proof}
    Given what we have computed in \ref{lemma_unstable_manifold_in_gen} it is enough to
  observe that on $S^{ct\pm}$ by (\ref{eqn_cond14}) we have
  \begin{align*}
    2\abs{\frac{\partial F_1}{\partial a_1}} + \sum_{i=2}^{\infty} \abs{\frac{\partial F_1}{\partial a_i}} + \abs{\frac{\partial F_i}{\partial a_1}}
    &\leq \abs{2\lambda_1 + 6c a_1^2} + \beta_1 C^{d_1} \\ &\leq \max\set{\abs{2 - 6\gamma_-^2}, \abs{2 - 6\gamma_+^2}}\lambda_1 + \beta_1 C^{d_1} < \frac{l}{4},
\end{align*}
  % \begin{align*}
  %     2\abs{\frac{\partial F_1}{\partial a_1}} + \sum_{i=2}^{\infty} \abs{\frac{\partial F_1}{\partial a_i}} + \abs{\frac{\partial F_i}{\partial a_1}}
  %     \leq \abs{2\lambda_1 + 6c a_1^2} + \beta_1 C^{d_1} < 10\lambda_1 + \beta_1 C^{d_1} < \frac{l}{4},
  % \end{align*}
  % because $-12 \lambda_1 \leq 6ca_1^2 \leq -\frac{3}{2} \lambda_1 $
\end{proof}

\begin{lemma} [\textbf{Step IIIc -- hyperbolicity of the born fixed points}] \label{lemma_neg_lognorm_after_bif_gen}
  Define for $\mu > \mu_b$
  \begin{align*}
    B_\mu^+ &:= \set{z \in B_\mu \mid \pi_{a_1}z \geq \gamma_+ C(\mu)} \\
    B_\mu^- &:= \set{z \in B_\mu \mid \pi_{a_1}z \leq -\gamma_+ C(\mu)}.
  \end{align*}
  Then for $\mu \in (\mu_b, \mu_+]$ there exist unique fixed points $z_\mu^\pm \in B_\mu^\pm$ of $\varphi_{\mu}$.
  Moreover, those fixed points are hyperbolic with respect to the cones given by
  the matrix $Q$ such that
  \begin{align*}
      Q_{11} &= -1, \\
      Q_{ii} &= 1, \quad i = 2, \dots, m, \\
      Q_{ii} &= -1, \quad i > m.
  \end{align*}
  and $W^u_{B^\pm_\nu}\expr{u^\pm_\nu}$ ($W^s_{B^\pm_\nu}\expr{u^\pm_\nu}$)
  are the images of horizontal (vertical) disks.
 \end{lemma}
  \begin{proof}
      We have by (\ref{eqn_cond10'})
      \begin{align*}
          2\frac{\partial F_1}{\partial a_1} + \sum_{i=2}^{\infty} \abs{\frac{\partial F_1}{\partial a_i}} + \abs{\frac{\partial F_i}{\partial a_1}} &=
          2\lambda_1 + 6ca_1^2 + 2\frac{\partial h_1}{\partial a_1} + \sum_{i=2}^{\infty} \abs{\frac{\partial h_1}{\partial a_i}} + \abs{\frac{\partial h_i}{\partial a_1}} \\
          &\leq  (2 - 6\gamma_+^2)\lambda_1 + \beta_1 C^{d_1} < 0.
      \end{align*}
     Proof in \textcolor{black}{the}  remaining directions is analogous as in the proof of Lemma \ref{lemma_unstable_manifold_in_gen}.
  \end{proof}

Modifying the proof of Theorem \ref{thm_bif} as we have modified the proof
of Theorem \ref{thm_bif_model} to prove Theorem \ref{thm_bif_model_unstable},
we get the main result.

\begin{theorem} \label{thm_bif_gen}
     Pitchfork bifurcation occurs in (\ref{eqn_gen_bif_form}) on $\interval{\lambda_1\expr{\mu_-}, \lambda_1\expr{\mu_+}}$.
\end{theorem}

\section{Proof of the heteroclinic connection away from the bifurcation}

Aside from the proof of bifurcation we want to prove that the heteroclinic connections
arising from it can be continued for further parameters. We will do this by
verifying the assumptions of the following theorem.

\begin{theorem} \label{thm_proof_of_conn}
    Consider equation (\ref{eq_pde}) and assume that the associated local semiflow. 
    Assume that
    \begin{itemize}
        % \item[(i)] \textbf{assumptions of Theorem \ref{thm_unstable_manifold_selfcb}}. Let $m \leq M$ be natural numbers. Let $c_u:X_m \to \mathbb{R}^m$ be a homeomorphism
        % and consider the self-consistent bounds of the form
        % $V := c_u^{-1}(\overline{B}(0,1)) \oplus \bigoplus_{k = m + 1} ^ \infty B_k$.
        % Assume that there exists $\varepsilon > 0$ such that for each
        % $n > M$ the set $P_n\expr{V} = c_u^{-1} (\overline{B}(0,1)) \oplus \bigoplus_{i = m + 1} ^ n B_k$
        % is an $\varepsilon$-isolating cuboid for $\varphi_n$.
        % Assume moreover that conditions (D) is satisfied
        % and (F') are satisfied.
        \item[(i)] there exists a set $S \subset H$ (isolating cuboid of the source) satisfies the assumptions of Theorem \ref{thm_unstable_manifold_selfcb},
        \item[(ii)] there exists a set $R \subset H$ (basin of attraction of the target) satisfies the assumptions of Theorem \ref{thm:attractfp},
        \item[(ii)] there exists a set $S' \subset S$ such that $P_{m}^\perp S' = P_{m}^\perp S$ (where $m$ is
        as in Theorem \ref{thm_unstable_manifold_selfcb})
        and a time $t > 0$ such that $\varphi(\interval{0, t}, S')$ exists and $\varphi(t, S') \subset R$.
    \end{itemize}
    Then, there exist fixed points $u_0 \in S, u_+ \in R$ such that there exists
    a heteroclinic connection from $u_0$ to $u_+$.
\end{theorem}
\begin{proof}
    Proof is analogous to the proof of the condition (P3) in Theorem \ref{thm_bif}.
\end{proof}

We elaborate on verifying (i) in Section \ref{sec_cc_and_ln_ver}. We prove existence of
the set $A$ and we verify (ii) by using the rigorous integration algorithm presented in \cite{ZPer2}.
With this approach we prove the following theorem.

\begin{theorem} \label{thm_comp_assisted_connections}
    For system (\ref{eq_ks_pde}) with $\mu \in \set{0.99, 0.75}$ there exists  a heteroclinic
    connection between two fixed points, the unstable zero solution
    and the attracting fixed point.
\end{theorem}

Numerical data from the proof is contained in Section \ref{sec_data_heteroclinic}.

\appendix
\section{Logarithmic norms} \label{app_lognorms}

Logarithmic norms allow us to obtain \textcolor{black}{the}  one-sided (with respect to time)
Lipschitz constants for the flows induced by the ODEs.
Here we recall their definition and basic properties, following \textcolor{black}{the}  presentation  from \cite{HNW}
given in \cite{ZM}.

\begin{definition}\cite[Def. 1.10.4]{HNW}
    Let $A$ be a square matrix on the normed space $(\mathbb{R}^n, \norm{\cdot})$. We call
    \begin{align*}
        \mu(A) := \lim_{h \to 0^+} \frac{\norm{I + hA} - 1}{h}
    \end{align*}
    the \emph{logarithmic norm} of $A$.
\end{definition}

The following theorem gives us computable formulas for the logarithmic norms.

\begin{theorem}\cite[Thm. 1.10.5]{HNW}\label{thm_lognorm_bounds}
    The logarithmic norm satisfies the following formulas.
    \begin{itemize}
        \item If $\norm{\cdot} = \norm{\cdot}_2$ is the $l_2$ norm, then $\mu_2(A) = \textnormal{largest eigenvalue of } \frac{A + A^T}{2}$.
        \item If $\norm{\cdot} = \norm{\cdot}_{\infty}$ is the max norm, then $\mu_\infty(A) = \max_k \sum_{k = 1}^\infty \abs{a_{kk}} + \sum_{i \neq k} \abs{a_{ki}}$.
    \end{itemize}
\end{theorem}

Now let $f: \mathbb{R}^n \to \mathbb{R}^n$ and consider the following autonomous ODE $x' = f(x)$
and denote by $\varphi$ the associated flow.
The usefulness of the logarithmic norms in our context comes from the fact that given
the solution $\varphi(t, x_0)$ and its perturbation $y(t)$ they give us an upper bound
for $\norm{y(t) - \varphi(t, x_0)}$. More precisely the following Lemma from \cite{ZM}
based on Theorem 1.10.6 from \cite{HNW} holds.

\begin{lemma} \label{lemma_lognorm_exp_bound}
    Let $y: [0, T] \to \mathbb{R}^n$ be a piecewise $C^1$ function and assume that $\varphi(\cdot, x_0)$
    is defined on $\interval{0, T}$. Suppose that $Z$ is a convex set such that there exists $l \in \mathbb{R},
    \delta > 0$ such that
    \begin{align*}
        y\expr{\interval{0, T}}, \varphi(\interval{0, T}, x_0) &\subset Z, \\
        \mu\expr{\frac{\partial f}{\partial x}\expr{\eta}} &\leq l, \; \textnormal{for } \eta \in Z, \\
        \norm{\frac{dy}{dt}\expr{t} - f\expr{y\expr{t}}} &\leq \delta.
    \end{align*}
    Then for $t \in \interval{0, T}$ we have
    \begin{align*}
        \norm{y(t) - \varphi(t, x_0)} &\leq e^{lt} \norm{y(0) - x_0} + \delta \frac{e^{lt} - 1}{l}, \textnormal{for } l \neq 0, \\
        \norm{y(t) - \varphi(t, x_0)} &\leq e^{lt} \norm{y(0) - x_0} + \delta t, \textnormal{for } l = 0.
    \end{align*}
\end{lemma}

Let us remark that when $y(t)$ is a solution of the ODE, then of course this theorem is true with $\delta = 0$.
\par
The lemma above is used in the proof of Theorem \ref{thm:conv-on-Trap}, which is proved in \cite{ZM},
where \textcolor{black}{the} constant in the lemma is uniform for all Galerkin projections.
Let $F$ be  \textcolor{black}{an}  admissible function on self-consistent
bounds $V$ (see Section \ref{sec:method}).
 \textcolor{black}{In}  the mentioned paper \textcolor{black}{the}  lemma is used with the $l_2$ norm, but we can also use it with any other norm,
we only need to adjust the constant $l$\textcolor{black}{; there} it is taken such that for every $i \in I$
and $x \in V$
\begin{align*}
    \frac{\partial F_{i}}{\partial x_{i}}\expr{x}  +
\sum_{k,\: k \neq i}
    \frac{1}{2} \expr{
        \abs{\frac{\partial F_{i}}{\partial x_{k}}(x)} +
        \abs{\frac{\partial F_{k}}{\partial x_{i}}(x)}}  \leq l.
\end{align*}
We indeed see that by Theorem \ref{thm_lognorm_bounds} and the Gershgorin theorem
it is a bound for the $l_2$ logarithmic norm for every Galerkin projection.
For the max norm we can see directly by the same theorem that we can take $l$ such that
for every $i \in I$ and $x \in V$
\begin{align*}
    \frac{\partial F_{i}}{\partial x_{i}}\expr{x}  +
\sum_{k,\: k \neq i}
        \abs{\frac{\partial F_{i}}{\partial x_{k}}(x)} \leq l.
\end{align*} 
\section{Elementary properties of the Brouwer degree} \label{app_brouwer_degree}

In this section we recall the definition and properties of the Brouwer degree we use.
For a detailed exposition see for example \cite[Chapter 12]{Sm}.
\par
Let $D \subset S \subset \mathbb{R}^n$ and assume that $D$ is an open set. Let
$f: S \to \mathbb{R}^n$ be a continuous function and pick $c \in \mathbb{R}^n$.
Suppose that $f^{-1}(c) \cap D$ is compact. If $\overline{D}$ is compact
and $\overline{D} \subset S$,
then the last condition is satisfied when $f^{-1}(c) \cap \partial D = \emptyset$.
 \textcolor{black}{If}  $f$ is a smooth map, then $f^{-1}(c)$ is finite. In this case if for
all $x \in f^{-1}(c)$ we have
$\det Df(x) \neq 0$ (it is then said
that $c$ is  \textcolor{black}{a}  \emph{regular value}), then the Brouwer degree $d(f, D, c)$ can be defined as
$$
\sum_{x \in f^{-1}(c)} \sgn \det Df(x).
$$
Then we can extend the definition of the degree to $c$ which is not a regular value and to $f$
which is not smooth.

\begin{theorem}[Solution property] \label{thm_solution_property}
    If $d(f, D, c) \neq 0$, then there exists $x \in D$ such that $f(x) = c$.
\end{theorem}

\begin{theorem}[Homotopy property] \label{thm_homotopy_property}
    Let $H:\interval{0, 1} \times S \to \mathbb{R}^n$ be continuous. Suppose that
    $$ \bigcup_{t \in \interval{0, 1}} H_t^{-1}(c) \cap D$$
    is compact. Then for $t \in \interval{0, 1}$ we have $d(H_t, D, c) = d(H_0, D, c)$.
\end{theorem} 
\section{Normal forms}

\subsection{General considerations}

Throughout this section $H$ is a Hilbert space.
\par
Let $k > 0$ and consider an infinite-dimensional ODE

\begin{align*}
    a_k' &= \lambda_k a_k + p(a) + f_k(a), \\
    a_i' &= \lambda_i a_i + f_i(a), i \neq k,
\end{align*}
where $p(a)$ = $\sum_j m_j$ and $m_j(a) = d_j \prod_{i = 1}^\infty a_i^{\alpha_{i, j}}$ (only finitely many $\alpha_{i, j}$ are different from $0$). For now we will limit ourselves to the formal considerations and
thus we make no assumptions on $\lambda_i$ and $f_i$.
\par
Now we discuss the transformation which removes the term
$p$. New variables are denoted by $$b = \expr{a_1, \dots, b_k, a_{k + 1}, \dots}$$ and the transformation is given by its inverse
$a_k\expr{b} := b_k + \sum_{j} c_j m_j(b)$, where $c_j \in \mathbb{R}$ will be chosen later.
We also denote $a(b) := \expr{a_1, \dots, a_k\expr{b}, a_{k + 1}, \dots}$.
We have

\begin{align} \label{eqn_norm_form_misc1}
    a_k' = b_k' \expr{1 + \sum_j c_j \frac{\partial m_j}{\partial a_k} \expr{b}}
           + \sum_j c_j \sum_{i \neq k} \expr{\lambda_i a_i + f_i\expr{a\expr{b}}} \frac{\partial m_j}{\partial a_i} \expr{b},
\end{align}

and on the other hand
we have

\begin{align} \label{eqn_norm_form_misc2}
    a_k' &= \lambda_k b_k + \sum_j \expr{c_j \lambda_k m_j\expr{b} + m_j\expr{b} + \overline m_j\expr{b}}
            + f_k\expr{a\expr{b}},
\end{align}

 where $$\overline m_j \expr{b} := m_j\expr{a\expr{b}} - m_j \expr{b}.$$

 For $\abs{x} < 1$ we have $\expr{1 + x}^{-1} = 1 - x + g(x),$ where
\begin{align} \label{g_def}
    g(x) := \sum_{j = 2}^\infty \expr{-1}^j x^j.
\end{align}
Thus comparing (\ref{eqn_norm_form_misc1}) and (\ref{eqn_norm_form_misc2}) we get by a simple rearrangement
\begin{align}
    \begin{split} \label{eqn_gen_normal_form}
    b_k' &= \lambda_k b_k + \sum_j \expr{c_j\lambda_k  + 1 - c_j \sum_{i} \lambda_i \alpha_{i, j}} m_j\expr{b} \\
         &+ R_k^1\expr{p, \lambda_1, \dots}\expr{b} + R_k^2\expr{p, \lambda_1, \dots, f_1, \dots}\expr{b}, \\
    a_i' &= \lambda_i a_i + f_i\expr{a\expr{b}}, i \neq k,
    \end{split}
\end{align}
where
\begin{align*}
    R_k^1&\expr{p, \lambda_1, \dots}\expr{b} := \sum_j \overline m_j \expr{b} + \lambda_k b_k g\expr{\sum_j c_j \frac{\partial m_j}{\partial a_k} \expr{b}} + \\
    &+ \expr{\sum_j \expr{c_j\lambda_k + 1 - c_j\sum_{i \neq k} \lambda_i \alpha_{i, j}} m_j\expr{b} + \overline{m_j}\expr{b}}\expr{- \sum_j c_j \frac{\partial m_j}{\partial a_k} \expr{b} + g\expr{\sum_j c_j \frac{\partial m_j}{\partial a_k} \expr{b}}} \\
    R_k^2&\expr{p, \lambda_1, \dots, f_1, \dots}\expr{b} := \expr{f_k\expr{a\expr{b}} - \sum_j c_j \sum_{i \neq k} f_i\expr{a\expr{b}} \frac{\partial m_j}{\partial a_i}\expr{b}} \cdot \\
      & \expr{1 - \sum_j c_j \frac{\partial m_j}{\partial a_k} \expr{b} + g\expr{\sum_j c_j \frac{\partial m_j}{\partial a_k} \expr{b}}}.
\end{align*}

Now assume that all $m_j$ have the same formal order
(i.e. in each $\sum_i \alpha_{i, j}$ is the same; otherwise it may happen that removing $m_{j_1}$ produces
$m_{j_2}$ and the change of variables does not reach the desired end).
Then we see that to remove term $p$ we need to have

\begin{align} \label{eqn_removal_cond}
    c_j = \frac{1}{\lambda_k - \sum_{i} \lambda_i \alpha_{i, j}}
\end{align}

which is possible to satisfy if and only if

\begin{align} \label{eqn_res_cond}
    \lambda_k - \sum_{i} \lambda_i \alpha_{i, j} \neq 0.
\end{align}

Terms which satisfy (\ref{eqn_res_cond}) are called \emph{non-resonant}.
\par

\subsection{Bounds for the sums arising in the Kuramoto--Sivashinsky equation}
\label{sec_ks_bounds_normal_forms}

Now we use bounds for nonlinear part of (\ref{eqn_ks_fourier_with_mu}) from \cite{ZM}
for sets which are symmetric with respect to $0$ on every coordinate. 
Then to use those
bounds in the analytical proof of the bifurcation we apply
them to the sets of the form
\begin{align*}
S_C := \interval{-\sqrt{2} C, \sqrt{2} C} \times \prod_{k = 2}^\infty \interval{\frac{-C^3}{k^s}, \frac{C^3}{k^s}},
\end{align*}
We introduce the following notations

\begin{align*}
  \lambda_k(\mu) &= k^2(1 - \mu k^2), \\
  \FS F_k^{\geq l}(a_1, a_2, \ldots) &= - k \sum_{i=l}^{k-l} a_i a_{k-i}, \quad \FS F_k := \FS F_k^{\geq 1} \\
  \IS F_k^{\geq l}(a_1, a_2, \ldots) &= 2k \sum_{i=l}^{\infty} a_i a_{k+i}, \quad \IS F_k := \IS F_k^{\geq 1} \\
  N_k^{\geq l}(a_1, a_2, \dots) &= \FS^{\geq l} F_k(a_1, a_2, \ldots) + \IS^{\geq l} F_k(a_1, a_2, \ldots),  \\
  F_k(\mu, a_1, a_2, \ldots) &= \lambda_k\expr{\mu}a_k + N_k(a_1, a_2, \ldots).
\end{align*}

Let $M \geq 1$ and consider a set
$$W := \prod_{k = 1}^M \interval{-A_k, A_k} \times \prod_{k = M + 1}^\infty \interval{-\frac{D}{k^s}, \frac{D}{k^s}},$$
where $A_k \geq 0$ for $k = 1, \dots, M$.
\begin{lemma} \cite[Lemma 3.1]{ZM}
    Assume that $1 \leq k \leq M$. Then for any $a \in W$ we have
    \begin{align*}
        \abs{\IS F_k^{\geq M + 1}\expr{a}} \leq \frac{D^2}{\expr{k + M + 1}^s \expr{s - 1} M^{s - 1}}.
    \end{align*}
\end{lemma}

\begin{lemma} \cite[Lemma 3.6]{ZM}
    Assume that $k > M$. Then for any $a \in W$ we have
    \begin{align*}
        \abs{\IS F_k^{\geq M + 1}\expr{a}} \leq \frac{2D^2}{k^{s - 2}\expr{M + 1}^{s} \expr{s - 2} M^{s - 1}}.
    \end{align*}
\end{lemma}

\begin{lemma} \cite[Lemma 3.5]{ZM}
    Assume that $k > 2M$. Then for any $a \in W$ we have
    \begin{align*}
        \abs{\FS F_k^{\geq M + 1}\expr{a}} \leq \frac{D^2}{k^{s - 2}}\expr{\frac{4^s}{\expr{2M + 1}^{s - 1}} + \frac{2^s}{\expr{s - 1} M^{s}}}.
    \end{align*}
\end{lemma}

% \begin{lemma} \cite[Lemma 3.1]{ZM}
%     Assume that $1 \leq k \leq M$. Then for any $a \in W$ we have
%     \begin{align*}
%         \abs{\IS F_k^{\geq l}\expr{a}} \leq 2k\sum_{i = l}^{M - k} A_i A_{k + i} + 2kD \sum_{i = M - k + 1}^M \frac{A_i}{\expr{k + i}^s}
%         + \frac{D^2}{\expr{k + M + 1}^s \expr{s - 1} M^{s - 1}}.
%     \end{align*}
% \end{lemma}
%
% \begin{lemma} \cite[Lemma 3.6]{ZM}
%     Assume that $k > M$. Then for any $a \in W$ we have
%     \begin{align*}
%         \abs{\IS F_k^{\geq_l}\expr{a}} \leq \frac{2D}{k^{s - 2} \expr{M + 1}}\expr{\sum_{i = l}^{M} A_i + \frac{D}{\expr{M + 1}^{s - 1} \expr{s - 2} M^{s - 1}}}.
%     \end{align*}
% \end{lemma}
%
% \begin{lemma} \cite[Lemma 3.5]{ZM}
%     Assume that $k > 2M$. Then for any $a \in W$ we have
%     \begin{align*}
%         \abs{\FS F_k\expr{a}} \leq \frac{D}{k^{s - 2}}\expr{\frac{2^{s + 1}}{2M + 1}\sum_{i = 1}^{M} A_i + \frac{D 4^s}{\expr{2M + 1}^{s - 1}} + \frac{D 2^s}{\expr{s - 1} M^{s}}}.
%     \end{align*}
% \end{lemma}
Observe that $\FS F_k^{\geq M + 1} = 0$ for $k \leq 2M$.

We will use the bounds above in the computer assisted proof of the bifurcation.
In the analytical proof we only need to know the order of those sums on $S_C$.
To reach normal form, we will transform only $a_1, a_2, a_3$.
The corollary below deals with terms which do not contain any of those.
\begin{cor} \label{cor_ks_bounds}
    There exist constants $\eta, \xi > 0$ such that for ali $k > 0$ we have
    \begin{align*}
        \sup_{a \in S_C} \abs{\IS F_k(a)^{\geq 4}} &< \frac{\eta C^6}{k^{s - 2}}, \\
        \sup_{a \in S_C} \abs{\FS F_k(a)^{\geq 4}} &< \frac{\xi C^6}{k^{s - 2}}.
    \end{align*}
\end{cor}

% OLD
% \begin{lemma}\label{lemma_ks_bounds_prev}
%     \begin{align*}
%         \sup_{b \in S_C} \abs{\IS F_1(b)} &\leq \sqrt{2}C^4 + C^6 + \frac{C^6}{4^s} + \frac{C^6}{5^s 3^{s - 1} \expr{s - 1}}, \\
%         \sup_{b \in S_C} \abs{\IS F_2(b)} &\leq 4\expr{\sqrt{2}C^4 + \frac{C^6}{4^s} + \frac{C^6}{5^s 2^{s - 1} \expr{s - 1}}}, \\
%         \sup_{b \in S_C} \abs{\IS F_k(b)} &\leq C^4 \expr{\frac{\sqrt{2}}{3k^{s - 2}} + }, \quad k > 2 \\
%         \sup_{b \in S_C} \abs{\IS F_k(b)} &\leq \frac{\sqrt{2}C^3}{k^{s - 2}}\expr{C + \frac{C^3}{2^{s - 1}\expr{s - 1}}}, \quad k > 1,\\
%         \sup_{b \in S_C} \abs{\FS F_k(b)} &\leq \frac{\sqrt{2}C^3}{k^{s - 2}}
%                                         \expr{\frac{2^{s + 1}}{3} C
%                                             + \frac{\sqrt{2}C^3 4^s}{3^{s + 1}}
%                                             + \frac{\sqrt{2}C^3 2^s}{\expr{s - 1}}}, \quad k > 2.
%     \end{align*}
% \end{lemma}
% \begin{proof}
%     fdfds
% \end{proof}

Now we provide bounds for sums of the derivatives of the Kuramoto--Sivashinsky equation. 
We return to $M \geq 1$ and the set $W$ as above. 
The two simple estimates in the following lemmas will be used throughout the
entire section.

\begin{lemma}
    For $l > 0$ we have
    \begin{equation*}
        \sum_{k=l}^\infty \frac{1}{k^s} < \int_{l-1}^\infty \frac{dx}{x^s}=\frac{1}{(s-1)(l-1)^{s-1}}.
      \end{equation*}
\end{lemma}

\begin{lemma}\label{lemma_sum_b_translated}
    Let $N > 0$. For $\iota = -\lfloor \frac{N}{2} \rfloor, \dots, 1, 2, \dots$ such that $N + \iota > M$ we have
    \begin{align*}
        \sup_{a \in W} \sum_{k = N}^\infty k \abs{a_{k + \iota}} &\leq \frac{\expr{1 + \max{\set{-\sgn{\iota}, 0}}}D}{\expr{s - 1}\expr{N + \iota - 1}^{s - 2}}, \quad \iota = 0, 1, \ldots.
    \end{align*}
\end{lemma}
\begin{proof}
   For $k \geq N$ and $\iota \geq -\frac{N}{2}$ we have $k + \iota \geq \frac{k}{2}$ if $\iota < 0$ and $k + \iota \geq k$ if $\iota \geq 0$.
    % and thus $\abs{ka_{k + \iota}} \leq \frac{kD}{\expr{k + \iota}^s} \leq \frac{2D}{\expr{k + \iota}^{s - 1}}$,
    Thus it is enough to use the estimate from the previous lemma. 
    % following estimate
    % \begin{equation*}
    %   \sum_{k=l}^\infty \frac{1}{k^s} < \int_{l-1}^\infty \frac{dx}{x^s}=\frac{1}{(s-1)(l-1)^{s-1}}.
    % \end{equation*}
\end{proof}

% As a corollary we get the following lemma.
%
% \begin{lemma}\label{lemma_sum_b_translated}
%     For $\iota = -3, -2, \dots, 1, 2, \dots$ we have
%     \begin{align*}
%         \sup_{a \in S_C} \sum_{k = M + 1}^\infty k \abs{a_{k + \iota}} &\leq \frac{D}{\expr{s - 1}\expr{M + \iota}^{s - 1}}, \quad \iota = 0, 1, \ldots.
%     \end{align*}
% \end{lemma}

%
% \begin{lemma}\label{lemma_sum_b_translated}
%     We have for $i \geq 2$
%     \begin{align*}
%         \sup_{b \in S_C} \sum_{k = l}^\infty k \abs{b_{k + \iota}} &\leq \frac{C^3}{\expr{s - 1}\expr{i - 1}^{s - 1}}, \quad \iota = 0, 1, \ldots, \\
%         \sup_{b \in S_C}\sum_{k = l}^\infty k \abs{b_{k - 1}} &\leq \frac{2C^3}{\expr{s - 1}\expr{i - 1}^{s - 1}}.
%     \end{align*}
% \end{lemma}
% \begin{proof}
%     It is enough to observe that
%     \begin{equation*}
%       \int_{l}^\infty \frac{dx}{x^s}=\frac{1}{(s-1)l^{s-1}}  <  \sum_{k=l}^\infty \frac{1}{k^s} < \int_{l-1}^\infty \frac{dx}{x^s}=\frac{1}{(s-1)(l-1)^{s-1}}.
%     \end{equation*}
% \end{proof}

\begin{lemma}\label{lemma_is_fs_derivatives_sum}
    % Let $e(k) = 1$ if $k \leq M$ and $e(k) = 0$ if otherwise. Then
    % We have
    % \begin{align*}
    %     \sup_{a \in W} \sum_{i = 1}^\infty \abs{\frac{\partial N_k}{\partial a_i}}
    %         &\leq 4k \expr{ \sum_{i = 1}^M A_i - e(k) A_k + \frac{D}{M^{s - 1} \expr{s - 1}}},\\
    %     \sup_{a \in W} \sum_{i = 1}^\infty \abs{\frac{\partial N_i}{\partial a_k}}
    %         &\leq 2M\sum_{i = 1}^M A_i + \frac{2D}{\expr{s - 1}M^{s - 2}} +
    %         4M \sum_{l=1}^{2M} A_i + \frac{2D}{\expr{s - 1}\expr{2M}^{s - 2}}, \quad k \leq M, \\
    %     \sup_{a \in W} \sum_{i = 1}^\infty \abs{\frac{\partial N_i}{\partial a_k}}
    %         &\leq 2k \frac{D}{\expr{s - 1} M^{s - 1}} - 2 \frac{D}{\expr{s - 1} M^{s - 2}} \\ &+ 2\expr{k - 1} \sum_{i = 1}^M A_i + 4k \sum_{i = 1}^M A_i + \frac{2D}{\expr{s - 1} k^{s - 2}}, \quad k > M.
    % \end{align*}

    We have
    \begin{align*}
        \sup_{a \in W} \sum_{i = 1}^\infty \abs{\frac{\partial N_k^{\geq M + 1}}{\partial a_i}}
            &\leq 4k \frac{D}{M^{s - 1} \expr{s - 1}}, \quad k \leq 2M + 1 \\
        \sup_{a \in W} \sum_{i = 1}^\infty \abs{\frac{\partial N_k^{\geq M + 1}}{\partial a_i}}
            &\leq 6k \frac{D}{M^{s - 1} \expr{s - 1}}, \quad k > 2M + 1 \\
        \sup_{a \in W} \sum_{i = 1}^\infty \abs{\frac{\partial N_i^{\geq M + 1}}{\partial a_k}} &= 0, \quad k \leq M \\
        \sup_{a \in W} \sum_{i = 1}^\infty \abs{\frac{\partial N_i^{\geq M + 1}}{\partial a_k}} &\leq \frac{2D}{s - 1} \expr{\frac{1}{M^{s - 2}} + \frac{k}{M^{s - 1}} + \frac{1}{k^{s - 2}}}, \quad M < k \leq 2M + 1 \\
        \sup_{a \in W} \sum_{i = 1}^\infty \abs{\frac{\partial N_i^{\geq M + 1}}{\partial a_k}} &\leq \frac{2D}{s - 1} \expr{\frac{1}{M^{s - 2}} + \frac{2k}{M^{s - 1}} + \frac{1}{k^{s - 2}}} , \quad k > 2M + 1  \\
    \end{align*}
\end{lemma}
\begin{proof}

For $k \leq 2M + 1$ we have

\begin{align*}
    \sup_{a \in W} \sum_{i = 1}^\infty \abs{\frac{\partial N_k^{\geq M + 1}}{\partial a_i}} &\leq
        \sup_{a \in W} 2k \sum_{i = M + 1}^\infty \abs{a_{i + k}} + \sum_{i = M + k + 1}^\infty \abs{a_{i - k}}  \\
        &\leq \sup_{a \in W} 4k \sum_{i = M + 1}^\infty \abs{a_i}\\
        &\leq 4k \frac{D}{M^{s - 1} \expr{s - 1}}% &\leq 2k \expr{ 2\sum_{i = 1}^M A_i - e(k) A_k + \frac{2D}{M^{s - 1} \expr{s - 1}}}
\end{align*}

and for $k > 2M + 1$ except the terms above we also have a term

\begin{align*}
    \sup_{a \in W} 2k\sum_{i = M + 1}^{k - M - 1} \abs{a_{k - i}}  \leq \sup_{a \in W} 2k\sum_{i = M + 1}^{\infty} \abs{a_{i}}  \leq  2k \frac{D}{M^{s - 1} \expr{s - 1}}%
\end{align*}

% \begin{align*}
%     \sup_{a \in W} \sum_{i = 1}^\infty \abs{\frac{\partial N_k^{\geq M + 1}}{\partial a_i}} &\leq
%     \sup_{a \in W} 2k \expr{\sum_{i = M + 1}^{k - M - 1} \abs{a_{k - i}} + \sum_{i = l}^{k - 1} \abs{a_{k + i}} + \abs{a_{2k}} + \sum_{i = k+1}^\infty \expr{\abs{a_{i - k}} + \abs{a_{i + k}}}} \\
%         &\leq \sup_{a \in W} 4k \sum_{i = l}^\infty \abs{a_i}\\
%         &\leq 4k \frac{D}{\expr{M - 1}^{s - 1} \expr{s - 1}}% &\leq 2k \expr{ 2\sum_{i = 1}^M A_i - e(k) A_k + \frac{2D}{M^{s - 1} \expr{s - 1}}}
% \end{align*}

\par
Let $M < k \leq 2M + 1$. Then we have

\begin{align*}
    \sup_{a \in W} \sum_{i = 1}^\infty \abs{\frac{\partial N_i^{\geq M + 1}}{\partial a_k}} &\leq
    \sup_{a \in W} \sum_{i = M + k + 1}^\infty 2i \abs{a_{i - k}} + \sum_{i = 1}^\infty  2i \abs{a_{k + i}} \\
        &= \sup_{a \in W} \sum_{i = M + 1}^\infty \expr{2i + 2k} \abs{a_i} + \sum_{i = 1}^\infty 2i \abs{a_{k + i}}\\
        &\leq \frac{2D}{s - 1} \expr{\frac{1}{M^{s - 2}} + \frac{k}{M^{s - 1}} + \frac{1}{k^{s - 2}}}.
\end{align*}

Now let $k > 2M + 1$. Then

\begin{align*}
    \sup_{a \in W} \sum_{i = 1}^\infty \abs{\frac{\partial N_i^{\geq M + 1}}{\partial a_k}} &\leq
    \sup_{a \in W} \sum_{i = M + 1}^{k - M - 1} 2i \abs{a_{k - i}} + \sum_{i = M + k + 1}^\infty 2i \abs{a_{i - k}} + \sum_{i = 1}^\infty  2i \abs{a_{k + i}}.
\end{align*}

Second and third term are as above, so we only need to bound

\begin{align*}
    \sum_{i = M + 1}^{k - M - 1} 2i \abs{a_{k - i}} &= \sup_{a \in W} \sum_{i = M + 1}^{k - i} 2\expr{k - i} \abs{a_i} \\
        &\leq 2k \sup_{a \in W} \sum_{i = M + 1}^\infty \abs{a_i} \leq \frac{2kD}{M^{s-1}\expr{s - 1}}.
\end{align*}
\end{proof}

\begin{cor} \label{cor_nonlinear_derivative}
    For $k > 0$ there exist $\eta, \xi > 0$ such that
    \begin{align*}
        \sup_{a \in S_C} \sum_{i = 1}^\infty \abs{\frac{\partial N_k^{\geq 4}}{\partial a_i}\expr{a}}
            &\leq \eta k C^3, \\
        \sup_{a \in S_C} \sum_{i = 1}^\infty \abs{\frac{\partial N_i^{\geq 4}}{\partial a_k}\expr{a}} &\leq \xi k C^3.  \\
    \end{align*}
\end{cor}
We are again going back to the set $W$. 
Denote by $\tilde N_k := N_k - N_k^{\geq M + 1}$.
Since for $k > 2M$ terms $a_1, \dots a_M$ do not mingle with each other neither
in $\IS F_k$ nor in $\FS F_k$, we have the following easy lemma.

\begin{lemma}
    For $k > 2M$ we have
    \begin{align*}
        \tilde N_k\expr{a} = -2k \sum_{i = 1}^M a_i a_{k - i} + 2k \sum_{i = 1}^M a_i a_{k + i}.
    \end{align*}
    Moreover
    \begin{align*}
        \sup_{a \in W} \abs{\tilde N_k\expr{a}} &\leq \frac{D\expr{2^{s} + 2}}{k^{s - 2} \expr{s - 1}} \sum_{i = 1}^M A_i \\
        \sup_{a \in W} \sum_{i = 1}^\infty \abs{\frac{\partial \tilde N_k}{\partial a_i}\expr{a}}
            &\leq \frac{D\expr{2^{s} + 2}}{k^{s - 2} \expr{s - 1}} + 4k \sum_{i = 1}^M A_i\\
        \sup_{a \in W} \sum_{i = 2M + 1}^\infty \abs{\frac{\partial \tilde N_i}{\partial a_l}\expr{a}}
            &\leq \frac{D}{s - 1}\expr{\frac{2}{\expr{2M - l - 1}^{s - 2}} + \frac{1}{\expr{2M + l - 1}^{s - 2}}}, \quad l \leq M \\
        \sup_{a \in W} \sum_{i = 2M + 1}^\infty \abs{\frac{\partial \tilde N_i}{\partial a_l}\expr{a}}
            &\leq 4l \sum_{i = 1}^M A_i, \quad l > M.
    \end{align*}
\end{lemma}
\begin{proof}
    Let $k > 2M$. We have
    \begin{align*}
        \sup_{a \in W} \abs{\tilde N_k\expr{a}}
            &\leq 2k \sup_{a \in W} \sum_{i = 1}^M \abs{a_i} \expr{\abs{a_{k - i}} + \abs{a_{k + i}}} \\
            &\leq \sup_{a \in W} 2k \expr{\sum_{i = 1}^M \abs{a_{k - i}} + \abs{a_{k + i}}} \expr{\sum_{i = 1}^M A_i}.
    \end{align*}
    We have
    \begin{align*}
        \sup_{a \in W} \sum_{i = 1}^M \abs{a_{k - i}} \leq \sup_{a \in W} \sum_{i = k - M}^{k - 1} \abs{a_i} 
            \leq \sup_{a \in W}  \sum_{i = \lfloor \frac{k}{2} \rfloor + 1}^\infty \abs{a_i} \leq \frac{2^{s - 1}D}{k^{s - 1}\expr{s - 1}}
    \end{align*}
    and
    \begin{align*}
        \sup_{a \in W} \sum_{i = 1}^M \abs{a_{k + i}} \leq \sup_{a \in W} \sum_{i = k + 1}^{M + k} \abs{a_i} 
            \leq \sup_{a \in W}  \sum_{i = k + 1}^\infty \abs{a_i} \leq \frac{D}{k^{s - 1}\expr{s - 1}}.
    \end{align*}
    
    % \begin{align*}
    %     \sup_{a \in W} \abs{\tilde N_k\expr{a}}
    %         &\leq 2k \sup_{a \in W} \sum_{i = 1}^M \abs{a_i} \expr{\abs{a_{k - i}} + \abs{a_{k + i}}} \\
    %         &\leq D \frac{2^{s + 1} + 2}{k^{s - 1}} \sum_{i = 1}^M A_i,
    % \end{align*}
    % because $k - i > \frac{k}{2}$ and $k + i > k$ for $i = 1,\dots,M$.
    \par
    Now observe that in $\sum_{i = 1}^\infty \abs{\frac{\partial \tilde N_k}{\partial a_i}\expr{a}}$
    each term $a_i, i < M,$ will arise twice -- when differentiating with respect to $a_{k \pm i}$.
    Thus
    \begin{align*}
        \sup_{a \in W} \sum_{i = 1}^\infty \abs{\frac{\partial \tilde N_k}{\partial a_i}\expr{a}}
            &\leq \sup_{a \in W} \expr{2k\sum_{i = 1}^M \expr{\abs{a_{k - i}} + \abs{a_{k + i}}} + 2\sum_{i = 1}^M \abs{a_i}} \\
            &\leq \frac{D\expr{2^{s} + 2}}{k^{s - 2}} + 4k \sum_{i = 1}^M A_i.
    \end{align*}
    \par
    Let $l \leq M$. Then 
    % by Lemma \ref{lemma_sum_b_translated} we have
    \begin{align*}
        \sup_{a \in W} \sum_{i = 2M + 1}^\infty \abs{\frac{\partial \tilde N_i}{\partial a_l}\expr{a}}
            &\leq \sum_{i = 2M + 1}^\infty 2l\expr{\abs{a_{i - l}} + \abs{a_{i + l}}} \leq \frac{D}{s - 1}\expr{\frac{2}{\expr{2M - l - 1}^{s - 2}} + \frac{1}{\expr{2M + l - 1}^{s - 2}}}.
    \end{align*}
    Finally, let $l > M$. Then
    \begin{align*}
        \sup_{a \in W} \sum_{i = 2M + 1}^\infty \abs{\frac{\partial \tilde N_i}{\partial a_l}\expr{a}}
            &\leq 2\sum_{i = 1}^M \expr{\expr{i + l} + \expr{l - i}} \abs{a_{i}} \leq 4l \sum_{i = 1}^M A_i.
    \end{align*}
    % above we use unoptimal bound, i + l arises from searching for k such that k - i = l,
    % but we see that in many cases k - M > l, so this term should be ignored in bound,
    % but it would in turn require additiona shaeningans
\end{proof}

\begin{cor} \label{cor_tilde_n}
    There exist constants $\alpha, \beta > 0$ such that for $k > 2M$ and $l > 0$
    we have
    \begin{align*}
        \sup_{a \in S_C} \abs{\tilde N_k\expr{a}} &\leq \frac{\alpha C^4}{k^{s - 2}} \\
        \sup_{a \in S_C} \sum_{i = 1}^\infty \abs{\frac{\partial \tilde N_k}{\partial a_i}\expr{a}}
            &\leq \beta k C\\
        \sup_{a \in S_C} \sum_{i = 2M + 1}^\infty \abs{\frac{\partial \tilde N_i}{\partial a_l}\expr{a}}
            &\leq \beta l C.
    \end{align*}
\end{cor}

\subsection{Proof of Lemma \ref{lemma_snd_normal_form}} \label{sec_app_normal_form}
\begin{definition}
    Let $W$ be self-consistent bounds. We call function $r:W \to \mathbb{R}$ a polynomial series on $W$ if there exists a set of multiindexes
    $I$ such that $r(a) = \sum_{\iota \in I} d_\iota a^\iota$.
\end{definition}

The following simple remark will prove to be very useful.
\begin{remark} \label{rem_R_linearity}
    $R_k^2$ is linear with respect to the sequence of polynomial series $\expr{f_1, f_2, \dots}$,
    i.e. for two sequences of polynomial series $\mathrm{f}_1, \mathrm{f}_2$ we have
    $$
    R_k^2\expr{p, \lambda_1, \dots, \mathrm{f}_1 + \mathrm{f}_2} =
        R_k^2\expr{p, \lambda_1, \dots, \mathrm{f}_1} + R_k^2\expr{p, \lambda_1, \dots, \mathrm{f}_2}
    $$
\end{remark}

Let $W \subset H$ be self-consistent bounds. For a sequence of polynomial
series $\mathrm{r} = \expr{r_1, r_2, \dots}$, where $r_1, r_2, \dots: W \to \mathbb{R}$,
we introduce the following notations
% \begin{align*}
%     \HS r_k\expr{a} &:= \sum_{i \neq k}^\infty \abs{\frac{\partial r_k}{\partial a_i}\expr{a}}, \\
%     \VS_{k} \mathrm{r}\expr{a} &:= \sum_{i \neq k}^\infty \abs{\frac{\partial r_i}{\partial a_k}\expr{a}}.
% \end{align*}
\begin{align*}
    \HS r_k\expr{a} &:= \sum_{i = 1}^\infty \abs{\frac{\partial r_k}{\partial a_i}\expr{a}}, \\
    \HS r_k &:= \sup_{a \in W} \HS r_k(a), \\
    \VS_{k} \mathrm{r}\expr{a} &:= \sum_{i = 1}^\infty \abs{\frac{\partial r_i}{\partial a_k}\expr{a}}, \\
    \VS_k \mathrm{r} &:= \sup_{a \in W} \VS_k \mathrm{r}(a).
\end{align*}
Observe that the following holds.

\begin{remark} \label{rem_vs_hs}
    Let $r_1, r_2, \dots$ be a sequence of polynomial series.
    In variables $b$ given by the inverse $a_k\expr{b} = b_k + p\expr{b}$ we have for the family of functions $\hat r_k: W \ni b \mapsto r_k\expr{a\expr{b}} \in \mathbb{R}$
    \begin{align*}
        \HS \hat{r_l}\expr{b} &\leq \HS r_l\expr{a\expr{b}} + \sum_{i} \abs{\frac{\partial p}{\partial a_i}\expr{b}\frac{\partial r_l}{\partial a_k}\expr{a\expr{b}}}, \\
        \VS_l \hat{\mathrm{r}} &\leq \VS_l \mathrm{r} + \sum_{i} \abs{\frac{\partial p}{\partial a_l}\expr{b}\frac{\partial r_i}{\partial a_k}\expr{a\expr{b}}}.
    \end{align*}
\end{remark}

The following bounds on $g$ are also easily seen to be true.

\begin{remark} \label{rem_g}
    We have
    \begin{align}
        \abs{g(x)} &\leq 2x^2 = O(x^2), \quad \abs{x} < \frac{1}{2}, \\
        \abs{g'(x)} &\leq 8x^2 + 4\abs{x} = O(\abs{x}), \quad \abs{x} < \frac{1}{2}.
    \end{align}
\end{remark}
% \begin{proof}
%     Let
%     $g_n(x) = \sum_{i = n}^\infty (-1)^i x^i$.
%     Then $\abs{g_n(x)} \leq \frac{\abs{x}^n}{1 - x} < 2\abs{x}^n$
%     and $\abs{g_n'(x)} \leq \sum_{i = n}^\infty i\abs{x}^{i - 1} \leq  (\frac{\abs{x}^n}{1 - \abs{x}})' \leq
%     \frac{n x^{n - 1}\expr{1 - x} - x^n}{(1 - x)^2} \leq 4( (n - 1) \abs{x^{n}} + nx^{n - 1})$.
% \end{proof}

From now on we always implicitly assume that $C$ is small enough,
so that the arguments of $g$ are always have absolute value smaller
than $\frac{1}{2}$.
Because of Remark \ref{rem_g} the following lemma is easy to verify.

\begin{lemma} \label{lemma_ignoring_R2}
    Fix $k, l > 0$. Assume that we have $p\expr{a} = c a_k^\iota a_1, \iota = 1, 2$
    and consider change of variables given by the inverse $a_l = b_l + p\expr{b}$.
    Assume that for polynomial series $r_1, r_k, r_l: W \to R$ we have $\sup_{a \in W} \abs{r_i(b)} = O(C^{k_1})$
    and $\HS r_i = O(C^{k_2})$, $i = 1, k, l$. Then
    \begin{align*}
        \sup_{b \in W} \abs{R^2_l\expr{p, \lambda_1, \lambda_k, \lambda_l, r_1, r_k, r_l}\expr{b}} &= O(C^{k_1}), \\
        \HS R^2_l\expr{p, \lambda_1, \lambda_k, \lambda_l, r_1, r_k, r_l} &= O(C^{k_2}).
    \end{align*}
\end{lemma}

We also need the simple lemma below.

\begin{lemma} \label{lemma_tilde_nonproblematic}
    Let $\mathrm{r} := \expr{r_1, r_2, \dots}$ of polynomial series. Fix $k > 0$.
    Assume that for all $i \in \mathbb{N}$ we have
    \begin{itemize}
        \item[(i)] every monomial $m$ in $r_i$ satisfies $m(a) = a_k^{\eta} a_1^{\zeta} \tilde m(a)$, where $\eta, \zeta \in \set{0, 1, 2, \dots}$ and $\tilde m$ does not
            depend on neither $a_k$ nor $a_1$.
        \item[(ii)] there exist $\alpha, p > 0$ such that setting $\tilde r = \sum \tilde m_j$
            we have $\sup_{a \in S_C} \abs{\tilde{r}_i\expr{a}} \leq \alpha \frac{C^3}{i^{p}}$,
        \item[(iii)] there exists $\beta > 0$ such that $\VS_i \mathrm{r} \leq \beta i C$,
        \item[(iv)] there exists $\gamma > 0$ such that $\HS r_i \leq \gamma i C$.
    \end{itemize}
    Then after a change of the variables given by the inverse $a_k = b_k + c a_1^\iota b_k, \iota = 1, 2,$
    or $a_k = b_k + c a_1^2$ the sequence of functions $\mathrm{\hat{r}} :=
    \expr{\hat{r}_1, \hat{r}_2, \dots}$ given by $\hat{r}: S_C \ni b \mapsto r\expr{a\expr{b}} \in \mathbb{R}$
    satisfies the conditions (i-iv).
\end{lemma}

Now we are ready to prove Lemma \ref{lemma_snd_normal_form}.

\begin{proof}
    We analyze (\ref{eqn_ks_fourier_with_mu}) with fixed $\mu >  \textcolor{black}{\frac{1}{2}}$.
    To prove the theorem, we need to remove the term $2 a_1 a_2$ in $F_1$ and
    $-a_1^2$ in $F_2$. After changes of variables we
    will moreover need to remove the terms $b_1^2 b_2$ from the transformed $F_1$
    and $b_1^3$ from the transformed $F_3$.
    In summary, only variables which we are gonna transform are the first three
    and it is easy to see that due to
    Corollaries \ref{cor_ks_bounds}, \ref{cor_nonlinear_derivative} and to Remark
    \ref{rem_R_linearity} and Lemma \ref{lemma_ignoring_R2} we can from now on ignore
    all nonlinear terms which do not depend
    on any of those variables.
    Moreover, by Corollary \ref{cor_tilde_n} and Lemma \ref{lemma_tilde_nonproblematic}
    we can ignore $\tilde N_7, \tilde N_8, \dots$.
    \par
    Let us write the equation without all terms we can dismiss by the considerations
    above.
    \begin{align*}
        a_1' &= \lambda_1 a_1 + 2 a_1 a_2 + 2 a_2 a_3 + 2 a_3 a_4 \\
        a_2' &= \lambda_2 a_2 - 2 a_1^2 + 4 a_1 a_3 + 4 a_2 a_4 + 4 a_3 a_5, \\
        a_3' &= \lambda_3 a_3 - 6 a_1 a_2 + 6 a_1 a_4 + 6 a_2 a_5 + 6 a_3 a_6, \\
        a_k' &= \lambda_k a_k + \tilde N_k\expr{a}, \quad k = 4, 5, 6.
    \end{align*}
    We have
    \begin{align*}
        \tilde N_4\expr{a} := - 8 a_1 a_3 - 4 a_2^2 - 8 a_1 a_5 - 8 a_2 a_5 - 8 a_2 a_6, \\
        \tilde N_5\expr{a} := - 10 a_1 a_4 - 10 a_2 a_3 + 10 a_1 a_6 + 10 a_2 a_7 + 10 a_3 a_8, \\
        \tilde N_6\expr{a} := - 12 a_1 a_5 - 12 a_2 a_4 - 6 a_3^2 + 12 a_1 a_7 + 12 a_2 a_8 + 12 a_3 a_9. \\
    \end{align*}
    We can directly check term by term that those can also be ignored. We thus see
    that we can limit ourselves to the equations for $a_1', a_2', a_3'$. In those
    we can again use Lemma \ref{lemma_ignoring_R2} to exclude most terms, getting
    finally
    \begin{align*}
        a_1' &= \lambda_1 a_1 + 2 a_1 a_2, \\
        a_2' &= \lambda_2 a_2 - 2 a_1^2, \\
        a_3' &= \lambda_3 a_3 - 6 a_1 a_2.
    \end{align*}
    % The term $-6 a_1 a_2$ in the equation for $a_3'$ is non-resonant we see that after
    % change of variables given by the inverse $a_3 = b_3 + c a_1 a_2$ this term
    % vanishes and moreover by analyzing (\ref{eqn_gen_normal_form}) all new arising
    % terms can be ignored by the considerations as above.
    % \par
    We now proceed to removing $2 a_1 a_2$. It can be removed by the transformation
    given by the inverse $a_1 = b_1 + \frac{1}{\lambda_2}m\expr{b_1, a_2}$, where
    $m(b_1, a_2) = 2 b_1 a_2$. Again studying (\ref{eqn_gen_normal_form}) we see
    that in the arising terms the only one which cannot be ignored is
    $\overline m(b_1, a_2) = -\frac{4}{\lambda_2} b_1^2 a_2$, which can be removed
    without introducing any problematic terms and
    $-\frac{1}{\lambda_2} (-2b_1^2) \frac{\partial m}{\partial a_2}\expr{b_1, a_2} = \frac{4}{\lambda_2} b_1^3$
    (term from $R^2_1$), which we leave.
    \par
    In remaining transformations we do not care about exact coefficient for the theoretical
    purposes \textcolor{black}{and it is easy to verify that they do not introduce any
    resonances for $\mu > \frac{1}{2}$}, so we do not derive them here.
    \par
    It remains to remove the term $-2b_1^2$ from the transformed $F_2$.
    Only problem we get after removal of this term is that
    in the transformed $F_3$ the term $-b_1^3$ (after
    substituting $a_2 = b_2 - d b_1^2$ into the $b_1 a_2$ term),
    so we remove it too
    and we see that after this coordinate change
    we finally get the form (\ref{eqn_ks_normal_form}) with required bounds.
\end{proof}

\subsection{Computer assisted calculation of the normal forms.}
\label{sec_comp_assisted_calculation_of_normal_forms}

Let $C$ be a fixed non-negative constant interval. By \emph{bound with the constant $C$} we understand
a triple $\expr{d, n, C},$ where $d$ is a non-negative interval and $n \in \mathbb{N}$.
Consider the bounds $\alpha_1 := \expr{d_1, n_1, C}, \alpha_2 := \expr{d_2, n_2, C}$, where $n_1 \leq n_2$.
We define the operations on bounds in the following way.
\begin{align*}
    \alpha_1 + \alpha_2 &:= \expr{d_1 + C^{n_2 - n_1}d_2, n_1, C}, \\
    \alpha_1 \alpha_2 &:= \expr{d_1 d_2, n_1 + n_2, C}, \\
    \alpha_1^m &:= \expr{d_1^m, n_1 m, C}, m \in \mathbb{N}, \\
    e\expr{\alpha_1} &:= d_1 C^{n_1}.
\end{align*}
We call the function $e$ the \emph{evaluation}. We also abuse the notation
and for the interval $x$ and the bound $\alpha$ by $x < \alpha$ we mean $x < e\expr{\alpha}$.
\par
We do not define the $'-'$ operator on the bounds. Thus if we say that $q\expr{\alpha_1, \dots, \alpha_N}$
is a polynomial on bounds, we implicitly mean a polynomial with positive coefficients
and without a free term.
\par
Now let $\alpha := \expr{\alpha_1, \alpha_2, \dots, \alpha_N, \alpha_\infty}$ be bounds with the constant $C$.
We associate with them a set of the form
\begin{align} \label{eqn_S_form}
S_{\alpha} := \interval{-e\expr{\alpha_1}, e\expr{\alpha_1}} \times \dots \times \interval{-e\expr{\alpha_N}, e\expr{\alpha_N}}
    \times \prod_{i = N + 1}^\infty \interval{-\frac{e\expr{\alpha_\infty}}{i^s}, \frac{e\expr{\alpha_\infty}}{i^s}}.
\end{align}

% For a polynomial $q$ when we write $q\expr{\beta_1, \dots, \beta_N; N, s}$
% we mean that coefficients of the polynomial depend on $N$ and $s$.

\begin{definition}
Let $f$ be a polynomial series.
Let polynomials $\overline f, \overline{\frac{\partial f}{\partial a_i}}, i = 1, \dots, N,
\overline{\HS f}$ be such for every sequence of bounds $\alpha$ we have
\begin{align*}
    \max_{a \in S_\alpha} \abs{f(a)} &< \overline f\expr{\alpha}, \\
    \max_{a \in S_\alpha} \abs{\frac{\partial f}{\partial a_i}\expr{a}} &< \overline{\frac{\partial f}{\partial a_i}}\expr{\alpha}, \\
    \max_{a \in S_\alpha} \abs{\HS f(a)} &< \overline{\HS f}\expr{\alpha}.
\end{align*}

The tuple $\expr{f, \bar{f}, \overline{\frac{\partial f}{\partial a_1}}, \dots, \overline{\frac{\partial f}{\partial a_1}}, \overline{\HS f}}$
is called \emph{function with bounds}.
\par
In the sequel we will simply write that $f$ is a function with bounds meaning the definition above.
% \begin{itemize}
%     \item $\max_{a \in S} \abs{f(a)} < \overline f\expr{\alpha_1, \dots, \alpha_N, \alpha},$
%     \item for $i \leq N$ we have $\max_{a \in S} \abs{\frac{\partial f}{\partial a_i}}\expr{a} < \xi\expr{s}\frac{\partial \overline f}{\partial a_i}\expr{\alpha_1, \dots, \alpha_N, \alpha},$
% \end{itemize}

\end{definition}

\begin{remark}
    Any polynomial $p$ of the first $N$ variables is a function with bounds. Indeed,
    if we take a polynomial $\overline p$ whose coefficients are absolute values of
    coefficients of $p$, then we can simply take
    $\overline{\frac{\partial p}{\partial a_i}} := \frac{\partial \overline p}{\partial a_i}$
    and $\overline{\HS p} := \HS \overline p$.
\end{remark}

Consider a change of variables given by the inverse
$a_k\expr{b} = b_k + p\expr{b}, k \leq N$, where $p$ is a polynomial of $N$ first variables and without a free term.
We also denote the associated coordinates on bounds $\beta := \expr{\alpha_1, \dots, \beta_i, \dots, \alpha_N}$
and $\alpha\expr{\beta} := \expr{\alpha_1, \dots, \beta_k + \overline p \expr{\beta}, \alpha_{k + 1}, \dots}$.
\par
The following lemma allows us to find bounds when we change coordinates.
Let $f$ be a function with bounds. Our goal is to find bounds for $f$ in the new coordinates.
% namely for the function $\hat f: S \ni b \mapsto f\expr{a\expr{b}} \in \mathbb{R}$.

\begin{lemma}
    Let $a_k\expr{b} = b_k + p\expr{b}, k \leq N$.
    Let $f$ be a function with bounds and for
    any sequence of bounds $\alpha$ define $\hat f: S_{\alpha} \ni b \mapsto f\expr{a\expr{b}} \in \mathbb{R}$.
    Define also for a sequence of bounds $\beta$
    \begin{align*}
        \overline{\hat{f}}\expr{\beta} &:= \overline f\expr{\alpha\expr{\beta}} \\
        \overline{\frac{\partial \hat{f}}{\partial b_i}}\expr{\beta} &:= \frac{\partial \overline p}{\partial b_i}\expr{\beta} \overline{\frac{\partial f}{\partial a_k}}\expr{\alpha\expr{\beta}} + \overline{\frac{\partial f}{\partial a_i}}\expr{\alpha\expr{\beta}}, \\
        \overline{\HS \hat f}\expr{\beta} &:= \HS \overline p \expr{\beta} \overline{\frac{\partial f}{\partial a_k}}\expr{\alpha\expr{\beta}}  + \overline{\HS f}\expr{\alpha\expr{\beta}}.
    \end{align*}
    Then $\hat{f}$ with the polynomials defined above is a function with bounds.
\end{lemma}

Now we introduce an algebra on the functions with bounds.

\begin{lemma}
    Let $f_1, f_2$ be functions with bounds. Define
    \begin{align*}
        \overline{f_1 + f_2} &:= \overline{f_1} + \overline{f_2}, \\
        \overline{\frac{\partial{f_1 + f_2}}{\partial a_i}} &:= \overline{\frac{\partial{f_1}}{\partial a_i}} + \overline{\frac{\partial{f_2}}{\partial a_i}}, \\
        \overline{\HS\expr{f_1 + f_2}} &:= \overline{\HS f_1} + \overline{\HS f_2}, \\
    \end{align*}
    Then $f_1 + f_2$ with the polynomials defined above is a function with bounds.
\end{lemma}

\begin{lemma}
    Let $f_1, f_2$ be functions with bounds. Define
    \begin{align*}
        \overline{f_1 f_2} &:= \overline{f_1} \cdot \overline{f_2}, \\
        \overline{\frac{\partial{f_1 f_2}}{\partial a_i}} &:= \overline{\frac{\partial{f_1}}{\partial a_i}}\overline{f_2} + \overline{f_1}\overline{\frac{\partial{f_2}}{\partial a_i}}, \\
        \overline{\HS\expr{f_1 f_2}} &:= \overline{\HS f_1} \cdot \overline{f_2} + \overline{f_1} \cdot \overline{\HS f_2}, \\
    \end{align*}
    Then $f_1 f_2$ with the polynomials defined above is a function with bounds.
\end{lemma}

The remaining thing we need is to represent the application of the function $g$
(given by \ref{g_def}) is a function with bounds. By Remark \ref{rem_g} the following
lemma is true.
\begin{lemma}
    Let $f$ be a function with bounds.
    Let $\hat g: g \circ f$. We define
    \begin{align*}
        \overline{\hat g} &:= 2 \overline{f}^2, \\
        \overline{\frac{\partial \hat g}{\partial a_i}} &:= \overline{\frac{\partial f}{\partial a_i}} \expr{8\overline{f}^2 + 4\overline{f}}, \\
        \overline{\HS \hat g} &:= \overline{\HS f} \expr{8\overline{f}^2 + 4\overline{f}}.
    \end{align*}
    Then $\hat g$ with the polynomials defined above is a function with bounds.
\end{lemma}
Let us note that it $\hat{g}$ so defined is a function with bounds only on sets on which
$\overline{f}$ has value less than $\frac{1}{2}$; we need to remember this condition
when implementing our algorithm.
\begin{definition}
Consider a family of the functions with bounds $f_1, f_2, \dots$.\\
We denote $\mathrm{f} = \expr{f_1, f_2, \dots}$. Assume
there exist a function with bounds $f$, polynomials with positive coefficients
$\overline{\VS_i \mathrm{f}}, i = 1, \dots, N,$ \\ 
$\overline{\VS \mathrm{f}}$,
 functions $\eta, \zeta, \xi: \mathbb{N} \to \mathbb{R}$
          such that on the sets of form (\ref{eqn_S_form}) we have
\begin{align*}
\overline{f}_i\expr{\alpha} &< \eta\expr{i} \overline{f}\expr{\alpha}, i > N \\
\overline{\HS f_i}\expr{\alpha} &< \zeta\expr{i} \overline{\HS f}\expr{\alpha}, i > N \\
\VS_i \mathrm{f} &< \overline{\VS_i \mathrm{f}}\expr{\alpha}, i \leq N, \\
\VS_i \mathrm{f} &< \xi(i)\overline{\VS \mathrm{f}}\expr{\alpha}, i > N. 
\end{align*}
We call $\expr{\mathrm{f}, f, \VS_1 \mathrm{f}, \dots, \VS_N \mathrm{f}, \VS \mathrm{f}}$ \emph{the sequence with bounds}.
\end{definition}

\begin{lemma}
    \label{lemma_bounded_sequence}
    Let $a_k\expr{b} = b_k + p\expr{b}, k \leq N$, where $p$ is a polynomial of the variables
    $a_1, \dots, a_N$. Let $\expr{\mathrm{f}, f, \VS_1 \mathrm{f}, \dots, \VS_N \mathrm{f}, \VS \mathrm{f}}$
    be a sequence with bounds.
    Define on any $S$ of form (\ref{eqn_S_form}) $\hat f_i: S \ni b \mapsto f_i\expr{a\expr{b}} \in \mathbb{R}$,
    $\hat f: S \ni b \mapsto f_i\expr{a\expr{b}} \in \mathbb{R}$
    and denote $\hat{\mathrm{f}} := \expr{\hat f_1, \hat f_2, \dots}$.
    We have
    \begin{align*}
       \overline{\hat f_i} \expr{\beta} &< \eta\expr{i} \overline{\hat f}\expr{\beta}, i > N, \\
       \overline{\HS \hat f_i}\expr{\beta} &< \zeta\expr{i} \overline{\HS \hat f}\expr{\beta}, i > N, \\
       \VS_i \mathrm{\hat{f}} &< \overline{VS_i \mathrm{\hat f}} :=  \overline{\frac{\partial p}{\partial a_i}}\expr{\beta} \overline{\VS_k \mathrm{f}}\expr{\alpha\expr{\beta}} + \overline{\VS_i \mathrm{f}}\expr{\alpha\expr{\beta}}, i \leq N, \\
       \VS_i \mathrm{\hat{f}} &< \xi(i) \overline{VS \mathrm{\hat f}} := \xi(i)\overline{\VS \mathrm{f}}\expr{\alpha\expr{\beta}}, i > N. \\
    \end{align*}
    Consequently, $\hat{\mathrm{f}}$ is a sequence with bounds.
\end{lemma}

In our bounds for the KS equation on the set $S_C$ we simply have $\eta(i) = \frac{1}{i^{s - 2}}$
and $\xi(i) = \zeta(i) = i$.
\par
We are now ready to describe how we approach computation of the normal forms on the
computer. Consider the following equation.

\begin{align} \label{eqn_normal_form_computation}
    \begin{split}
    a_1' &= \lambda_1 a_1 + p_1\expr{a} + f_1\expr{a}, \\
        &\dots \\
    a_N' &= \lambda_N a_n + p_N\expr{a} + f_N\expr{a}, \\
    a_i' &= \lambda_i a_i + f_i\expr{a},
    \end{split}
\end{align}

where $\expr{f_1, \dots, f_N, \dots}$ is a sequence with bounds.
Now assume that $p_i = \sum_j m_j^i, i \leq N,$, where $m_j^i$ are non-resonant
and of the same formal order.
The arithmetic on the bounds given in this section together with
(\ref{eqn_gen_normal_form}) allows us to state that after a
change of the variables given by the inverse $a_k = b_k + \sum_j c_j m_j^k(b)$,
where $c_j$ are given by (\ref{eqn_removal_cond}). By Lemma \ref{lemma_bounded_sequence},
this gives us again equation of the form (\ref{eqn_normal_form_computation}), but
with $p_k$ removed. This allows us to do a computer-assisted proof of the Theorem \ref{thm_bif_gen}.

% \begin{lemma}
%     Let $a_k\expr{b} = b_k + p\expr{b}, k \leq N$.
%     Define $\hat f_i: S \ni b \mapsto f_i\expr{a\expr{b}} \in \mathbb{R}$
%     and $\hat f: S \ni b \mapsto f\expr{a\expr{b}} \in \mathbb{R}$.
%     Then
%     $\sup_{b \in S} \abs{\hat f_i\expr{b}} < \eta\expr{k} \overline{\hat f}\expr{\beta}, i > N$,.
% \end{lemma}

% \subsection{Transformations}
%
% % Consider an equation
% %
% % \begin{align*}
% %     a_1' &= \lambda_1 a_1 + p_1\expr{a_1, \dots, a_N} + f_1\expr{a}, \\
% %     &\ldots \\
% %     a_N' &= \lambda_1 a_1 + p_N\expr{a_1, \dots, a_N} + f_N\expr{a}, \\
% %     a_k' &= \lambda_k a_k + f_k\expr{a}, k > N,
% % \end{align*}
% %
% % where $p_1, \dots, p_N$ are polynomials.
%
% Consider an equation
%
% \begin{align*}
%     a_i' &= \lambda_i a_i + f_i\expr{a}, i \in \mathbb{N}.
% \end{align*}
%
% Fix $N > 0$ and $k \leq N$. We will consider transformation given by the inverse
% $a_i = b_i + p\expr{b}$, where $p$ is a polynomial of first $N$ variables.
% \par
% Assume that we have
% $$
% \max_{a \in S} \abs{f_i\expr{a}} \leq \xi\expr{s}\overline f_i
% $$

\section{Computer assisted proof of the bifurcation in the KS equation}
\label{sec_appendix_data_from_the_proof_of_bifurcation}

% \subsection{Comments on how to satisfy the conditions (\ref{eqn_cond1} -- \ref{eqn_cond14})}

Fix $M > 0$.
As in Section \ref{sec_ks_bounds_normal_forms}, we represent the nonlinear part
$N_k = \tilde{N}_k + N_k^\geq{M+1}$.
For the first $2M$ variables we hold all terms of $\tilde{N}$ explicitly,
and for the remaining variables we use the bounds derived in the mentioned section.
We then use the approach described in Section \ref{sec_comp_assisted_calculation_of_normal_forms}
to get to the form (\ref{eqn_gen_bif_form}) and to obtain $\alpha, \beta$ such that we can set
$\alpha_k = \alpha, \beta_k = \beta$ for all $k > 2M$.

\subsection{Bifurcation when $\mu = 1$}

By verifying the assumptions of Theorem \ref{thm_bif_gen}, we were able to prove the following.

\begin{theorem} \label{thm_comp_assisted_bif1}
     Pitchfork bifurcation occurs in the KS equation on the
    interval $$\interval{\lambda_1\expr{1.01}, \lambda_1(0.99)}.$$
\end{theorem}

To pass to the normal form we have used changes of the variables as
in the proof of Lemma \ref{lemma_snd_normal_form}; it is also clear
from this lemma that we have required monotonicity of $c(\mu)$ and $C(\mu)$.
In verifying the conditions (\ref{eqn_cond1} -- \ref{eqn_cond14}) (of course
the conditions considering the unstable directions are trivially satisfied
because there are no unstable directions here) we had
$\zeta = 1.2, \omega = 3, s = 6, p = 4, K = \sqrt{2}$ (we proved that (\ref{eqn_cond2})
is satisfied with this $K$
in the proof of Lemma \ref{lemma_big_isolation}),  $l = 6, \gamma_- = 0.05, \gamma_+  = 1.03$.
The value of $C$ at the end of the bifurcation was $C(0.99) = 0.172047$.
Using rigorous numerics we have found the following bounds for the values of $h_i$s in the normal form
for any $\mu \in [0.99, 1)$.

\begin{center}
    \begin{tabular}{ |c|c|c| }
    \hline
        $k$  & $\frac{\alpha_k}{k^{4}}$ & order of $C(\mu)$ \\
    \hline
     $1$ & $0.0369053$ & $5$\\
     $2$ & $0.364407$ & $4$\\
     $3$ & $0.133062$ & $4$\\
     $4$ & $0.490892$ & $4$\\
     $5$ & $0.0141437$ & $4$\\
     $6$ & $0.0012905$ & $4$\\
     $7$ & $0.000470282$ & $4$\\
     $8$ & $0.000254643$ & $4$\\
     $9$ & $0.0025041$ & $4$\\
     \hline
\end{tabular}
\end{center}

Now we present our bound for the derivatives of $h_i$. We present the numbers
as sum of the derivatives in rows + sum of the derivatives in columns.

\begin{center}
    \begin{tabular}{ |c|c|c| }
    \hline
        $k$  & $k \beta_k$ & order of $C(\mu)$ \\
    \hline
     $1$ & $0.471591 + 3.03665$ & $3$\\
     $2$ & $5.18774 + 7.74087$ & $1$\\
     $3$ & $14.69 + 14.8565$ & $1$\\
     $4$ & $19.9531 + 19.9139$ & $1$\\
     $5$ & $24.8912 + 49.7795$ & $1$\\
     $6$ & $29.868 + 59.7215$ & $1$\\
     $7$ & $34.8568 + 68.9155$ & $1$\\
     $8$ & $40.5246 + 57.1589$ & $1$\\
     $9$ & $44.802 + 66.0368$ & $1$\\
     \hline
\end{tabular}
\end{center}

We also evaluated largest absolute value of an expression which was the argument of $g$
verify that our change of variables is valid on the given parameter range. The result was $0.00122761$,
so we are very far away from the boundary of validity of the change of the variables.

\subsection{Bifurcation when $\mu = 0.25$}

\begin{theorem} \label{thm_comp_assisted_bif2}
     Pitchfork bifurcation occurs in the KS equation on the
    interval $$\interval{ \lambda_2(0.26), \lambda_2\expr{0.25 - 0.0002}}.$$
\end{theorem}

We have removed the terms in the given order (we skip constants)
\begin{itemize}
    \item $a_2 a_4$ from $a_2'$,
    \item $a_2^2$ from $a_4'$,
    \item $a_2^2 a_6$ from $a_2'$,
    \item $a_2^3$ from $a_6'$.
\end{itemize}
We had $\zeta = 1.2, \omega = 3, s = 6, p = 4, K = \frac{1 - 8\mu_+}{1 - 16\mu_+} = 1.00053$ (it is easy to verify (\ref{eqn_cond2})
is satisfied with this $K$),  $l = 0.066, \gamma_- = \frac{1}{15}, \gamma_+  = \frac{10}{11}$.
The value of $C$ at the end of the bifurcation was $C(\mu_+) = 0.0979273$.
Using rigorous numerics we have found the following bounds for the values of $h_i$s in the normal form
on the considered range

\begin{center}
    \begin{tabular}{ |c|c|c| }
    \hline
        $k$  & $\frac{\alpha_k}{k^{4}}$ & order of $C(\mu)$ \\
    \hline
     $1$ & $0.0408299$ & $4$\\
     $2$ & $0.00872972$ & $5$\\
     $3$ & $0.114076$ & $4$\\
     $4$ & $0.173405$ & $4$\\
     $5$ & $0.018412$ & $4$\\
     $6$ & $0.00611738$ & $4$\\
     $7$ & $0.00134561$ & $4$\\
     $8$ & $0.231311$ & $4$\\
     $9$ & $0.000211203$ & $4$\\
     $10$ & $0.00151919$ & $4$\\
     $11$ & $6.70185e-05$ & $4$\\
     $12$ & $4.38124e-05$ & $4$\\
     $13$ & $0.000560199$ & $4$\\
     \hline
\end{tabular}
\end{center}

Now we present our bound for the derivatives of $h_i$. We present the numbers
as sum of the derivatives in rows + sum of the derivatives in columns.

\begin{center}
    \begin{tabular}{ |c|c|c| }
    \hline
        $k$  & $k \beta_k$ & order of $C(\mu)$ \\
    \hline
     $1$ & $4.84811 + 9.79133$ & $1$\\
     $2$ & $0.337498 + 1.51915$ & $3$\\
     $3$ & $14.5452 + 14.593$ & $1$\\
     $4$ & $9.79877 + 14.6903$ & $1$\\
     $5$ & $24.2421 + 24.2423$ & $1$\\
     $6$ & $28.9501 + 29.0443$ & $1$\\
     $7$ & $33.9389 + 67.8762$ & $1$\\
     $8$ & $38.7946 + 77.5731$ & $1$\\
     $9$ & $43.6378 + 86.9637$ & $1$\\
     $10$ & $48.7229 + 96.6365$ & $1$\\
     $11$ & $53.8555 + 75.1066$ & $1$\\
     $12$ & $87.8383 + 82.3756$ & $1$\\
     $13$ & $63.0293 + 149.656$ & $1$\\
     \hline
\end{tabular}
\end{center}

We also evaluated largest absolute value of an expression which was the argument of $g$
verify that our change of variables is valid on the given parameter range. The result was $1.54422e-09$. 
\section{Verifying assumptions of Theorem \ref{thm:attractfp} and (\ref{eq_cc-diag-dom})} \label{sec_cc_and_ln_ver}

In this sections we quote needed bounds from \cite{ZA}.
\par
Let $\mu > 0$ be fixed. We denote
\begin{align*}
 N_k(a) = \FS F_k(a) + \IS F_k(a).
\end{align*}

The formal first derivatives of $F$ are given by
\begin{align*}
 \frac{ \partial N_i }{\partial a_j} &=  2i a_{i+j} ,  \quad \mbox{for $i=j$} \\
 \frac{ \partial N_i }{\partial a_j} &= -2i a_{i-j} + 2i a_{i+j},  \quad \mbox{for $j < i$} \\
 \frac{ \partial N_i }{\partial a_j} &= 2i a_{j-i} + 2i a_{i+j},   \quad \mbox{for $j > i$} \\
  \frac{ \partial F_i }{\partial a_j} &= i^2(1 - \mu i^2)\delta_{ij} + 2i
  \sum_{k \geq 1} ( - \delta_{k,i-j} +  \delta_{k,i+j} + \delta_{k,j-i} )a_k, \\
  \frac{\partial F_i}{\partial \mu} &= -i^4 a_i.
\end{align*}

\begin{lemma}\cite{ZA}[Lemma 5.1]
\label{lem:CC} Let $A:H \to H$ be a linear coordinate change of
the form
\begin{align*}
  A : X_m \oplus Y_m &\to X_m \oplus Y_m \\
  A(x \oplus y) &= Ax \oplus y.
\end{align*}
Let ${\tilde F}= A \circ F \circ A^{-1}$ (${\tilde F}$ is $F$
expressed in new coordinates).
\begin{align*}
 \frac{\partial {\tilde F}_i}{\partial x_j} &= \sum_{k,l=1}^m A_{ik}  \frac{\partial F_k}{\partial x_l} A^{-1}_{lj}
   \quad  \mbox{for $i \leq m$ and $j \leq m$} \\
 \frac{\partial {\tilde F}_i}{\partial x_j} &= \sum_{k \leq m} A_{ik}  \frac{\partial F_k}{\partial x_j}
   \quad \mbox{for $i \leq m$ and $j > m$} \\
 \frac{\partial {\tilde F}_i}{\partial x_j} &= \sum_{l \leq m}  \frac{\partial F_i}{\partial x_l} A^{-1}_{lj}
   \quad \mbox{for $i > m$ and $j \leq m$} \\
 \frac{\partial {\tilde F}_i}{\partial x_j} &=  \frac{\partial F_i}{\partial x_j}
   \quad \mbox{for $i > m$ and $j > m$} \\
  \frac{\partial \tilde{F}_k}{\partial \mu} &= \sum_{i=1}^m A_{ki}
  \frac{\partial F_i}{\partial \mu}, \quad  \mbox{if $k \leq m$} \\
  \frac{\partial \tilde{F}_k}{\partial \mu} &= \frac{\partial F_k}{\partial \mu}, \quad
    \mbox{if $k > m$}.
\end{align*}
\end{lemma}

Consider now the KS equation and we assume that $V=W \oplus T$ is the
self-consistent bounds for a fixed point. Let the numbers $m<M$ as
in conditions {\bf C1,C2,C3} and we assume that $a_k^{\pm}=\pm
\frac{C}{k^s}$ for $k > M$ (as in \cite{ZM}).

Let $A \in \mathbb{R}^{m \times m}$ be a coordinate change around
an approximate fixed point in $X_m$ for $m$-dimensional Galerkin
projection of (\ref{eqn_ks_fourier_with_mu}).  This matrix induces a coordinate
change in $H$. It is optimal to choose $A$ so that the
$m$-dimensional Galerkin projection of $F$ is very close to the
diagonal matrix (or to the block diagonal one when the complex
eigenvalues are present).

We will use the new coordinates in $H$. We also change the scalar product so
that the new coordinates are orthogonal.

To make notation more uniform we  set
$A_{ij}=\delta_{ij}$ if $i > m$ or $j>m$ and  $a_k=0$ for $k \leq
0 $.

We denote
\begin{equation*}
S(l,V) := \sum_{k,\ k \geq l} \max_{a \in W} |a_k|
\end{equation*}

We will estimate $S(l)$ using the following lemma.
\begin{lemma}\cite{ZA}[Lemma 5.2]
\label{lem:Sestm} Assume that $|a_k(V)| \leq \frac{C}{k^s}$ for $k
> M$, $s
> 1$, then
\begin{align*}
  S(l) &< \sum_{k=l}^M |a_k(V)|  + \frac{C}{(s-1)M^{s-1}}, \quad \mbox{for $l \leq M$} \\
  S(l) &< \frac{C}{(s-1)(l-1)^{s-1}}, \quad \mbox{for $l > M$}.
\end{align*}
\qed
\end{lemma}

We set
\begin{align*}
  {\overline S}(l) &:= \sum_{k=l}^M |a_k(V)|  + \frac{C}{(s-1)M^{s-1}}, \quad \mbox{for $l \leq M$} \\
  {\overline S}(l) &:=  \frac{C}{(s-1)(l-1)^{s-1}}, \quad \mbox{for $l > M$}.
\end{align*}

\subsection{How to verify assumptions of Theorem \ref{thm:attractfp}}

We use formulas derived in \cite{ZA}.
\par
Let $A$ be coordinate change as in \ref{lem:CC}. We wish to provide upper bound for

\begin{equation}
  l_i := \sup_{x \in V} \frac{\partial {\tilde F}_i}{\partial x_i}(x) + \sum_{j \neq i} \abs{\frac{\partial {\tilde F}_i}{\partial x_j}(x)}.
\end{equation}

We denote

\begin{equation}
  S_{ND}(i) := \sum_{j \neq i} \abs{\frac{\partial {\tilde F}_i}{\partial x_j}}.
\end{equation}

\begin{lemma}\cite[Lemma 5.3]{ZA}
If $1 \leq i \leq m$, then
\begin{align*}
S_{ND}(i) &\leq \sum_{j \leq M, j \neq i}   \sup_{x \in W} \abs{\frac{\partial \IS F_i}{\partial x_j}(x)}
+ \\
&+ \sum_{k\leq m} 2k \abs{A_{i,k}} \expr{S(M+1-k) + S(M+1+k)}
\end{align*}
\end{lemma}

\begin{lemma}\cite[Lemma 5.4]{ZA}
If $m < i \leq M$, then
\begin{align*}
S_{ND}(i) &\leq
\sum_{j \leq M, j \neq i} \sup_{x \in W} \left| \frac{\partial \IS F_i}{\partial x_j}(x)  \right|
+ \\
&+ 2i\left( S(M+1-i) + S(M+1+i) \right)
\end{align*}
\end{lemma}

\begin{lemma}\cite[Lemma 5.5]{ZA}
Assume that $M < i $. Then
\begin{align*}
  S_{ND}(i) &\leq \overline{S}_{ND}(i):= 2i \|A^{-1}\|_\infty \cdot (S(i-m) +
      S(i+1))  2iS(i+m+1)  + 4iS(1) \\
  l_i &\leq \overline{l}_i := i^2(1-\nu i^2) + \overline{S}_{ND}(i)
\end{align*}
\end{lemma}

\begin{lemma}\cite[Lemma 5.6]{ZA}
If for some $n > M$ holds
\begin{align*}
    \overline{l}_n < 0,
\end{align*}
 then
\begin{equation*}
    0 > \overline{l}_i > \overline{l}_j, \quad i < j, i \geq n
\end{equation*}

\end{lemma}

\subsection{How to verify (\ref{eq_cc-diag-dom})}

We would like to derive the formula for
\begin{equation*}
  \Gamma_{i}:=  2 \inf_{x \in W}\left|\frac{\partial \tilde{F}_i}{\partial x_i}(x) \right| - \sum_{j, j \neq i}
     \sup_{x \in W} \left| Q_{jj}\frac{\partial \tilde{F}_j}{\partial x_i}(x) +
     Q_{ii} \frac{\partial \tilde{F}_i}{\partial x_j}(x)  \right|
\end{equation*}

We will do analogous computations as in \cite{ZA}[Section 5].

We denote
\begin{align*}
K(l,n) &=\sum_{j,\ j \geq l} j \max_{a \in W} |a_{j+n}|, \\
 S_{ND}(i)&= \sum_{j, j \neq i}
     \sup_{x \in W} \left| Q_{jj}\frac{\partial \tilde{F}_j}{\partial x_i}(x) +
     Q_{ii} \frac{\partial \tilde{F}_i}{\partial x_j}(x)  \right|.
\end{align*}

We will need the following estimate.
\begin{lemma}
\label{lem:estm-sum}
\begin{equation}
  \int_{l}^\infty \frac{dx}{x^s}=\frac{1}{(s-1)l^{s-1}}  <  \sum_{k=l}^\infty \frac{1}{k^s} < \int_{l-1}^\infty \frac{dx}{x^s}=\frac{1}{(s-1)(l-1)^{s-1}}
\end{equation}
\end{lemma}

\begin{lemma}
\label{lem:Kestm} Assume that $|a_k(V)| \leq \frac{C}{k^s}$ for $k
> M$, $s > 2$, then
\begin{align*}
  K(l,n) &\leq \sum_{j,\ l \leq j  \leq  M-n}  j \max_{a \in W}
  |a_{j+n}| + \\
    &+ C\left(  \frac{1}{(s-2)(r+n-1)^{s-2}} -  \frac{1}{(s-1)(r+n)^{s-1}}
    \right),
\end{align*}
where $r=r(l,M,n)=\max(l,M-n+1)$
\end{lemma}
\begin{proof}
\begin{align*}
  K(l,n,M) =  \sum_{j,\ l \leq j  \leq  M-n}  j \max_{a \in W} |a_{j+n}| +
  \sum_{j \geq \max(l,M-n+1)} \frac{j C}{(j+n)^s}.
\end{align*}
For $r$ such that $r +n > M$ we have
\begin{align*}
  \sum_{j,\ j \geq r} \frac{j}{(j+n)^s} &=  \sum_{j,\ j \geq r} \frac{j +n}{(j+n)^s}
   - \sum_{j,\ j \geq r} \frac{n}{(j+n)^s} \leq \\
   &\leq \frac{1}{(s-2)(r+n-1)^{s-2}} -  \frac{1}{(s-1)(r+n)^{s-1}}.
\end{align*}
 \end{proof}

\begin{lemma}
\label{lem:cc-offd-i-small} If $1 \leq i \leq M$, then
\begin{align*}
  S_{ND}(i) &\leq \sum_{j \leq M, j \neq i}   \sup_{x \in W} \left| Q_{jj}\frac{\partial \tilde{F}_j}{\partial x_i}(x) +
     Q_{ii} \frac{\partial \tilde{F}_i}{\partial x_j}(x)  \right|
     + \\
     &+ \sum_{k\leq M} 2k |A_{i,k}| \left( S(M+1-k) + S(M+1+k)
     \right)+ \\
     &+ 2 \sum_{k \leq M} |A^{-1}_{ki}|  ( K(M+1,-k) +  K(M+1,k))
\end{align*}
\end{lemma}
\begin{proof} We have
\begin{align*}
   \sum_{j > M}   \sup_{x \in W} \left| Q_{jj}\frac{\partial \tilde{F}_j}{\partial x_i}(x)
   +     Q_{ii} \frac{\partial \tilde{F}_i}{\partial x_j}(x)  \right|
     &\leq  \sum_{j > M} \sup_{x \in W} \left| \frac{\partial \tilde{F}_j}{\partial x_i}(x)\right|
     + \\
     &+\sum_{j > M} \sup_{x \in W} \left| \frac{\partial \tilde{F}_i}{\partial x_j}(x)  \right|,
\end{align*}

\begin{align*}
    \sum_{j > M}  \left| \frac{\partial \tilde{F}_i}{\partial x_j}\right|  &\leq
 \sum_{j> M}  \sum_{k\leq M} |A_{i,k}| \left| \frac{ \partial F_k}{\partial x_j} \right|
    = \\
 &= \sum_{k\leq M}  |A_{i,k}| \sum_{j> M}   \left| \frac{ \partial F_k}{\partial x_j} \right| \leq \\
 &\leq \sum_{k\leq M} 2k |A_{i,k}| \sum_{j> M} \left( |a_{j-k}| + |a_{j+k}|\right) \leq \\
 &\leq \sum_{k\leq M} 2k |A_{i,k}| \left( S(M+1-k) + S(M+1+k) \right),
\end{align*}

\begin{align*}
  \sum_{j > M}  \left| \frac{\partial \tilde{F}_j}{\partial x_i}\right|
  \leq \sum_{j>M} \sum_{k \leq M} \left| \frac{\partial F_j}{\partial
  x_k}\right| \cdot |A^{-1}_{ki}| = \sum_{k \leq M}
  |A^{-1}_{ki}| \sum_{j >M} \left| \frac{\partial F_j}{\partial
  x_k}(x)\right|.
\end{align*}

Observe that for $k \leq M$ we have
\begin{align*}
\sum_{j >M} \left| \frac{\partial F_j}{\partial
  x_k}\right| \leq \sum_{j >M} 2j (|a_{j-k}| + |a_{j+k}|) \leq
  2 K(M+1,-k) + 2 K(M+1,k).
\end{align*}

 \end{proof}

Taking into account that $A_{ij}$ differs from $\delta_{ij}$ only
for $i,j \leq m$  we obtain from the lemma above the following two
lemmas
\begin{lemma}
\label{lem:cc-offd-i<=m} If $i \leq m$, then
\begin{align*}
  S_{ND}(i) &\leq \sum_{j \leq M, j \neq i}   \sup_{x \in W} \left| Q_{jj}\frac{\partial \tilde{F}_j}{\partial x_i}(x) +
     Q_{ii} \frac{\partial \tilde{F}_i}{\partial x_j}(x)  \right|
     + \\
     &+ \sum_{k\leq m} 2k |A_{i,k}| \left( S(M+1-k) + S(M+1+k)
     \right)+ \\
     &+ 2 \sum_{k \leq m} |A^{-1}_{ki}|  ( K(M+1,-k) +  K(M+1,k)).
\end{align*}
\end{lemma}

\begin{lemma}
\label{lem:cc-offd-i-medium} If $m < i \leq M$, then
\begin{align*}
  S_{ND}(i) &\leq \sum_{j \leq M, j \neq i}   \sup_{x \in W} \left| Q_{jj}\frac{\partial \tilde{F}_j}{\partial x_i}(x) +
     Q_{ii} \frac{\partial \tilde{F}_i}{\partial x_j}(x)  \right|
     + \\
     &+ 2i\left( S(M+1-i) + S(M+1+i)
     \right)+ \\
     &+ 2( K(M+1,-i) +  K(M+1,i)).
\end{align*}
\end{lemma}

\begin{lemma}
\label{lem:cc-offd-i-fartail} Let  matrix $\tilde{A} \in
\mathbb{R}^{m \times m}$ be given by
\begin{equation}
   \tilde{A}_{jk}= |A_{jk}| k.
\end{equation}

Assume that $M < i $. Then
\begin{align*}
  S_{ND}(i) &\leq \overline{S}_{ND}(i):= \left( 2\|\tilde{A}\|_1 +  2i \|A^{-1}\|_\infty \right) \cdot (S(i-m) +
      S(i+1)) + \\
    &+ 2iS(i+m+1)  + (8i-2)S(1) +\\
     &+2 \left(K(m+1,i) + K(1,0)
       \right),\\
  \Gamma_i &\geq \overline{\Gamma}_i:= 2|i^2(1-\mu i^2)| - \overline{S}_{ND}(i)
\end{align*}
\end{lemma}
\begin{proof}
It is easy to see that
\begin{align*}
\Gamma_{i} &=  2 \inf_{x \in W}\left|\frac{\partial
\tilde{F}_i}{\partial x_i}(x) \right| - \sum_{j, j \neq i}
     \sup_{x \in W} \left| Q_{jj}\frac{\partial \tilde{F}_j}{\partial x_i}(x) +
     Q_{ii} \frac{\partial \tilde{F}_i}{\partial x_j}(x)  \right|
     \geq \\
     &\geq 2 |i^2(1 - \mu i^2)| - \sum_{j, j \leq m }
     \sup_{x \in W} \left| Q_{jj}\frac{\partial \tilde{F}_j}{\partial x_i}(x) +
     Q_{ii} \frac{\partial \tilde{F}_i}{\partial x_j}(x)  \right|
     +\\
    &-   \sum_{j, j > m}  \sup_{x \in W} \left| \frac{\partial N_j}{\partial x_i}(x)\right|
      -    \sum_{j, j > m}   \sup_{x \in W} \left| \frac{\partial N_i}{\partial
      x_j}(x)    \right|
\end{align*}

 We have
\begin{align*}
  \sum_{1 \leq j \leq m}  \sup_{x \in W} \left| Q_{jj}\frac{\partial \tilde{F}_j}{\partial x_i}(x) +
     Q_{ii} \frac{\partial \tilde{F}_i}{\partial x_j}(x)  \right|
     &\leq
      \sum_{1 \leq j \leq m} \sum_{1 \leq k \leq m} |A_{jk}| \cdot \left| \frac{\partial N_k}{\partial x_i}\right| +\\ &+
     \sum_{1 \leq j \leq m} \sum_{1 \leq k \leq m}  \left| \frac{\partial N_i}{\partial
     x_k}\right|\cdot \left|A^{-1}_{kj} \right|
\end{align*}

For the first  term we obtain
\begin{align*}
 \sum_{1 \leq j \leq m} \sum_{1 \leq k \leq m} |A_{jk}| \cdot \left| \frac{\partial N_k}{\partial
 x_i}(V)\right| &\leq
  \sum_{1 \leq j \leq m} \sum_{1 \leq k \leq m}
 |A_{jk}| \left(2 k (|a_{i-k}(V)|  + |a_{i+k}(V)|) \right) \leq \\
 &\leq 2  \sum_{1 \leq j \leq m} \sum_{1 \leq k \leq m}
 |A_{jk}| k |a_{i-k}(V)|  +  2  \sum_{1 \leq j \leq m} \sum_{1 \leq k \leq m}
 |A_{jk}| k |a_{i+k}(V)| = \\
&= 2 \left\|\tilde{A} \cdot (|a_{i-1}(V)|,\dots,|a_{i-m}(V)|)
 \right\|_1 + \\
 &+ 2 \left\|\tilde{A} \cdot (|a_{i+1}(V)|,\dots,|a_{i+m}(V)|)
 \right\|_1 \leq \\
 &\leq 2 \|\tilde{A}\|_1 (S(i-m) + S(i+1)).
\end{align*}

Now we consider the second term
\begin{align*}
   \sum_{1 \leq j \leq m} \sum_{1 \leq k \leq m}  \left| \frac{\partial N_i}{\partial
     x_k}(V)\right|\cdot \left|A^{-1}_{kj} \right| &\leq
 \sum_{j=1}^m  \sum_{k=1}^m 2i(|a_{i-k}(V)| + |a_{i+k}(V)|)|A^{-1}_{kj}| = \\
 &= 2i \sum_{j=1}^m  \sum_{k=1}^m|{A^{-1}}^T_{j,k}| \cdot |a_{i-k}(V)|  +
  2i \sum_{j=1}^m  \sum_{k=1}^m|{A^{-1}}^T_{j,k}| \cdot |a_{i+k}(V)| = \\
  &= 2i \left \| {A^{-1}}^T \cdot (|a_{i-1}(V)|,\dots,|a_{i-m}(V) |)
  \right\|_1 + \\
   &+ 2i \left \| {A^{-1}}^T \cdot (|a_{i+1}(V)|,\dots,|a_{i+m}(V) |)
  \right\|_1  \leq \\
    &\leq 2i \|{A^{-1}}^T\|_1 \cdot (S(i-m) + S(i+1)) = \\
    &= 2i \|A^{-1}\|_\infty \cdot (S(i-m) + S(i+1)).
\end{align*}
Hence we obtained the following estimate
\begin{align}\label{eq_est-main-i>m}
    \begin{split}
    \sum_{1 \leq j \leq m}  &\sup_{x \in W} \left| Q_{jj}\frac{\partial \tilde{F}_j}{\partial x_i}(x) +
     Q_{ii} \frac{\partial \tilde{F}_i}{\partial x_j}(x)  \right| \leq \\
     &\leq
      \left( 2\|\tilde{A}\|_1 +  2i \|A^{-1}\|_\infty \right) \cdot (S(i-m) +
      S(i+1)).
  \end{split}
\end{align}

Now we will study the infinite sums.

First we will show  that
\begin{align}
   \sum_{j=m+1}^\infty  \sup_{x \in W} \left| \frac{\partial N_i}{\partial x_j}(x) \right|
     \leq 2i (S(i+m+1) + 2 S(1)). \label{eq_row-esti>m}
\end{align}

Indeed, we have
\begin{align*}
   \sum_{j=m+1}^\infty  \sup_{x \in W} \left| \frac{\partial N_i}{\partial x_j}(x) \right| &\leq
   \sum_{m <j < i}  2i(|a_{i-j}(V)| + |a_{i+j}(V)|) +  \\
    &+ 2i |a_{2i}(V)| + \sum_{j > i} 2i(|a_{j-i}(V)| + \\ &+ |a_{i+j}(V)|) \leq \\
   &\leq 2i \left( \sum_{j>m} |a_{i+j}(V)| + \sum_{m < j < i} |a_{i-j}(V)| + \sum_{j >i} |a_{j-i}(V)|  \right)< \\
   &< 2i \left( S(i+m+1) + 2S(1) \right).
\end{align*}

Now we will prove  that
\begin{equation}
  \sum_{j, j > m} \sup_{x \in W}   \left| \frac{\partial N_j}{\partial
  x_i}(x) \right| \leq   2 \left(K(m+1,i) + (2i-1)S(1)  + K(1,0)
     \right). \label{eq_col-esti>m}
\end{equation}

We have
\begin{align*}
  \sum_{j, j > m} \sup_{x \in W}   \left| \frac{\partial N_j}{\partial
  x_i}(x) \right| &\leq \sum_{m<j<i} 2j \left(|a_{i-j}(V)| + |a_{i+j}(V)|
  \right) + \\ &+ 2i |a_{2i}(V)| + \sum_{i<j} 2j \left(|a_{j-i}(V)| + |a_{i+j}(V)|
  \right) = \\
  &= 2 \left( \sum_{j,j>m} j|a_{j+i}(V)| + \sum_{j=m+1}^{i-1} j |a_{i-j}(V)| +
     \sum_{j=i+1}^\infty j|a_{j-i}(V)| \right) \leq \\
     &\leq 2 \left(K(m+1,i) + (i-1)S(1) + i S(1) + K(1,0)
     \right)= \\
     &= 2 \left(K(m+1,i) + (2i-1)S(1)  + K(1,0)
     \right).
\end{align*}

From
(\ref{eq_est-main-i>m},\ref{eq_row-esti>m},\ref{eq_col-esti>m}) it
follows that
\begin{align*}
  \Gamma_{i} &\geq  2|i^2(1-\mu i^2)| + \\
   &- \left( 2\|\tilde{A}\|_1 +  2i \|A^{-1}\|_\infty \right) \cdot (S(i-m) +
      S(i+1)) + \\
      &- 2i (S(i+m+1) + 2 S(1)) + \\
       &- 2 \left(K(m+1,i) + (2i-1)S(1)  + K(1,0)
     \right) = \\
      &= 2|i^2(1-\mu i^2)|  -  \left( 2\|\tilde{A}\|_1 +  2i \|A^{-1}\|_\infty \right) \cdot (S(i-m) +
      S(i+1)) + \\
       &- 2iS(i+m+1)  - (8i-2)S(1) -2 \left(K(m+1,i) + K(1,0) \right)
\end{align*}

\end{proof}

\begin{definition}
We define a function $f:\mathbb{N} \to \mathbb{R}_+$ by
\begin{align*}
  f(i) &= \frac{1}{i}(\overline{S}_{ND}(i) + 2 S(1))= \\
        &= \left( 2\|\tilde{A}\|_1/i +  2 \|A^{-1}\|_\infty \right) \cdot (S(i-m) +
      S(i+1)) + \\
     &+ 2S(i+m+1)  + 8 S(1) +\\
     &+ \frac{2}{i} \left(K(m+1,i) + K(1,0) \right).
\end{align*}
\end{definition}

The following lemma gives us the criterion for
$\Gamma_i$ to be positive for  $i$ large enough.
\begin{lemma}
\label{lem:large-i}  If for some $n > M$, holds
\begin{align*}
 2(\mu n^3 - n)   > f(n)
\end{align*}
 then
\begin{equation}
   {\overline \Gamma}_i > {\overline \Gamma}_j > 0 \qquad  \mbox{for } \qquad i > j, \quad
   j \geq n
\end{equation}
\end{lemma}
\begin{proof} From Lemma \ref{lem:cc-offd-i-fartail} it follows
that
\begin{displaymath}
    {\overline \Gamma}_i =i\left(2(\mu i^3 - i)
     - f(i) \right) + 2 S(1),
\end{displaymath}
where $f(i)$ is a positive decreasing function.

 It is easy to see that the function $i \mapsto (\mu i^3 - i)$ is
increasing and positive for $i \geq n$. Therefore
$\overline{\Gamma}_i$ is increasing and positive for $i  \geq n$.
\end{proof}

\section{Data from the proof of the heteroclinic connection} \label{sec_data_heteroclinic}

\subsection{Case of $\mu = 0.99$}

We have verified the existence of an isolating
cuboid $N$ of the source point $0$. The number of main modes was $M = 16$.
First three modes of $N$ (the dominating directions)
are equal to $$\set{[-0.0639609, 0.0474541],[-0.00153622, 0.00137737],[-0.00134996, 0.00100156]}$$
and the tail is given by $\frac{C}{k^s} = \frac{3.85104}{k^{20}}$.
We have verified the cone conditions on $N$ by verifying the condition
(\ref{eq_cc-diag-dom}) with $\lambda \approx 0.000197702$.
\par
The approximate target point is $\set{0.173564,-0.00508437,7.43631e-05}$.
We have found its basin of attraction $R$, whose first three modes are
equal to $$\set{[0.16661, 0.181094],
[-0.00559088, -0.0046115],
[-3.3707e-06, 0.000146987]}$$ and the tail is given by
$\frac{C}{k^s} = \frac{6.53907e-05}{k^{12}}$. On $R$
the logarithmic norm is less than $-6.53835e-05$.
\par
We integrate the $a_1$-right boundary of the set $N$, i.e.
the set $N_1^+ \oplus \pi_{a_2, a_3, \dots} N$.
We have used integration algorithm with $m = 8$ main modes and
with $M = 16$ modes held explicitly (the near tail is modes $9, 10, \dots, 16$).
We have used the time step $h = 0.0002$. After $1500000$ time steps the
bounds we have obtained for the first three modes are
\begin{align*}\big\{&\interval{0.16956, 0.16956} + \interval{-1.33044e-06, 1.33044e-06},\\
&\interval{-0.00485229, -0.00485229} + \interval{-7.62099e-08, 7.62099e-08},\\
&\interval{6.93307e-05, 6.93307e-05} + \interval{-1.63319e-09, 1.63319e-09}\big\}\end{align*}
and for the tail are
$\frac{C}{k^s} = \frac{1.63471e-14}{k^{16}}$.

\subsection{Case of $\mu = 0.75$}

We have verified the existence of an isolating
cuboid $N$ of the source point $0$. The number of main modes was $M = 16$.
First three modes of $N$ (the dominating directions)
are equal to $$\set{[-0.273971, 0.260674],[-0.0402309, 0.0410396],[-0.0310602, 0.028675]}$$
and the tail is given by $\frac{C}{k^s} = \frac{4.26881e+08}{k^{20}}$.
We have verified the cone conditions on $N$ by verifying the condition
(\ref{eq_cc-diag-dom}) with $\lambda \approx 0.00368637$.
\par
The approximate target point is $\set{0.712361, -0.123239, 0.0101786}$.
We have found its basin of attraction $R$, whose first three modes are
equal to $$\set{[0.658428, 0.764085],
[-0.146629, -0.0983076],
[0.00373825, 0.015761]}$$ and the tail is given by
$\frac{C}{k^s} = \frac{202.882}{k^{12}}$. On $R$
the logarithmic norm is less than $-0.0035869$.
\par
We integrate the $a_1$-right boundary of the set $N$, i.e.
the set $N_1^+ \oplus \pi_{a_2, a_3, \dots} N$.
We have used integration algorithm with $m = 8$ main modes and
with $M = 16$ modes held explicitly (the near tail is modes $9, 10, \dots, 16$).
We have used time step $h = 0.0005$. After $20000$ time steps the
bounds we have obtained for the first three modes are
\begin{align*}\big\{&\interval{0.710814, 0.710814} + \interval{-4.67305e-05, 4.67305e-05},\\
&\interval{-0.122663, -0.122663} + \interval{-1.67197e-05, 1.67197e-05},\\
&\interval{0.0100529, 0.0100529} + \interval{-2.04451e-06, 2.04451e-06}\big\}\end{align*} and for the tail are
$\frac{C}{k^s} = \frac{0.000374526}{k^{14}}.$

\end{document}